\documentclass[ijoc,nonblindrev]{informs3} 

\SingleSpacedXI

\usepackage{subfig}
\usepackage{float}
\usepackage{booktabs}
\usepackage{appendix}
\usepackage{lineno,hyperref}

\usepackage{float}
\usepackage{graphicx}
\usepackage{subfig}
\usepackage{mathrsfs}
\usepackage{graphics}
\usepackage{float}
\usepackage{amsmath}
\usepackage{amsfonts}
\usepackage{mathrsfs}
\usepackage{amsmath}
\usepackage{longtable}
\usepackage{optidef}
\usepackage[subnum]{cases}

\usepackage{optidef}
\usepackage{graphicx}
\usepackage{float}
\usepackage{graphicx}
\usepackage{subfig}
\usepackage{mathrsfs}
\usepackage{float}
\usepackage[table]{xcolor}
\usepackage{amsmath}
\usepackage{amsfonts}
\usepackage{booktabs}
\usepackage{natbib}
\usepackage{bbm}
\usepackage{algorithm}
\usepackage[noend]{algpseudocode}
\usepackage{multirow}
\algdef{SE}[DOWHILE]{Do}{doWhile}{\algorithmicdo}[1]{\algorithmicwhile\ #1}
\usepackage{pgfplots}
\usepackage[letterpaper, portrait, margin=1in]{geometry}

\usepackage{booktabs}
\usepackage{tabularx}

\newcommand{\norm}[1]{\left\lVert#1\right\rVert}

\usepackage{newtxtext} 

\newcommand{\Acal}{\mathcal{A}}
\newcommand{\Pcal}{\mathcal{P}}

\newcommand{\Fcal}{\mathcal{F}}
\newcommand{\Rcal}{\mathcal{R}}
\newcommand{\Tcal}{\mathcal{T}}
\newcommand{\Wcal}{\mathcal{W}}
\newcommand{\Ccal}{\mathcal{C}}
\newcommand{\Scal}{\mathcal{S}}
\newcommand{\Ical}{\mathcal{I}}
\newcommand{\Lcal}{\mathcal{L}}

\newcommand{\tmin}{t^{\text{min}}}

\hypersetup{
    colorlinks,
    linkcolor=blue,
    citecolor=blue,
    urlcolor=blue
}
\usepackage{natbib}
 \bibpunct[, ]{(}{)}{,}{a}{}{,}%
 \def\bibfont{\small}%

\TheoremsNumberedThrough     

\EquationsNumberedThrough    

\begin{document}


\RUNAUTHOR{Mo et al.} 

\RUNTITLE{Individual Path Recommendation}

\TITLE{Individual Path Recommendation Under Public Transit Service Disruptions Considering Behavior Uncertainty}

\ARTICLEAUTHORS{%
\AUTHOR{Baichuan Mo$^{1,2,*}$, Haris N. Koutsopoulos$^3$, Zuo-Jun Max Shen$^4$, Jinhua Zhao$^5$}
\AFF{$^{*}$ Corresponding author}
\AFF{$^{1}$Department of Civil Engineering, Tsinghua University, Beijing, China, 100084}
\AFF{$^{2}$Department of Civil and Environmental Engineering, Massachusetts Institute of Technology, Cambridge, MA 02139}
\AFF{$^{3}$Department of Civil and Environmental Engineering, Northeastern University, Boston, MA 02115}
\AFF{$^{4}$Department of Industrial Engineering and Operations Research, University of California, Berkeley, Berkeley, CA 94720}
\AFF{$^{5}$Department of Urban Studies and Planning, Massachusetts Institute of Technology, Cambridge, MA 20139}
} 

\ABSTRACT{%
\small Public transit passengers need guidance during service disruptions. This study proposes an individual-based path (IPR) recommendation model. The model decides which paths to recommend for each passenger with the objective of minimizing system travel time and respecting passengers' path choice preferences. We assume the recommendations could affect passengers' path choice probabilities, but their actual choices are uncertain. This behavior uncertainty makes the problem a stochastic optimization with decision-dependent distributions. We propose a single-point approximation method to eliminate the expectation operator by introducing two new concepts: $\epsilon$-feasibility and $\Gamma$-concentration, which control the mean and variance of path flows in the optimization problem. The approximation yields a tractable single-stage mixed integer linear formulation, which can be solved efficiently with Benders decomposition. The approximation gap is approved to be bounded from the above. Additional theoretical analysis shows that $\epsilon$-feasibility and $\Gamma$-concentration are strongly connected to expectation and chance constraints in a typical stochastic optimization formulation, respectively. The model is implemented in a real-world case study using data from an urban rail disruption in the Chicago Transit Authority system, and a synthetic case study with varied network sizes and incident locations. In the real-world case study, results show that the proposed IPR model reduces the average travel times in the system by 6.6\% compared to the status quo and by 4.2\% compared to a capacity-based benchmark model. In the synthetic case study, the proposed model shows 15.0\% to 1.8\% lower system travel time compared to the capacity-based benchmark, depending on the network sizes and demand situations.   
}%


\maketitle

%


\section{Introduction}\label{intro}

With aging systems and near-capacity operations, service disruptions often occur in urban public transit (PT) systems. For example, the Chicago Transit Authority (CTA) reported an average of 75.4 abnormal events per day in 2019, with at least 1.1 events each day resulting in delays of over 20 minutes \citep{mo2022toward}. These incidents may result in passenger delays, cancellation of trips, and economic losses \citep{cox2011transportation}. During a significant disruption where the service is interrupted for a relatively long period of time (e.g., 1 hour), affected passengers usually need to find an alternative path or use other travel modes (such as transfer to another bus route). However, due to a lack of knowledge of the system (especially during incidents), the routes chosen by passengers may not be optimal or even cause more congestion \citep{mo2021impact}. For example, during a rail disruption, most of the passengers may choose bus routes that are parallel to the interrupted rail line as an alternative. However, given limited bus capacity, parallel bus lines may become oversaturated and passengers have to wait for a long time to board due to being denied boarding (or left behind).  

One of the strategies to better guide passengers is to provide path recommendations so that passenger flows are re-distributed in a better way and the system travel times are minimized. This can be seen as solving an optimal passenger flow distribution (or assignment) problem over a public transit network. However, different from the typical flow redistribution problem, there are several unique characteristics and challenges for the path recommendation problem under PT service disruptions.  
\begin{itemize}
    \item Passengers may have different preferences on different alternative paths. This heterogeneity suggests that we cannot treat a group of passengers simply as flows. Individualization is needed in the path recommendation design.
    \item Passengers may not follow the recommendation. When providing a specific path recommendation to a passenger, their actual path choices are uncertain (though the recommendation may change their preferences). This behavior uncertainty brings challenges to the recommendation system design and has not been considered in the path recommendation literature. In the context of individualization, the behavior uncertainty is also individual-specific, which requires a more granular modeling approach.   
\end{itemize}

To tackle these challenges, this study proposes an individual-based path recommendation (IPR) model to reduce system congestion during public transit disruptions considering passenger behavior uncertainty. The decision variable in this study is $x_{p,r}$, a binary variable indicating whether path $r$ is recommended to passenger $p$ or not. However, passengers may not follow the recommendation. Their behavior uncertainty is modeled using a conditional path choice probability distribution given their received path recommendation. The original IPR formulation yields a stochastic optimization with decision-dependent distributions \citep{goel2006class,drusvyatskiy2023stochastic}. We propose a single-point approximation method to eliminate the expectation operator by introducing two new concepts: $\epsilon$-feasibility and $\Gamma$-concentration. The former describes the relationship between path flows and recommendations from the ``expectation'' perspective. While the latter constrains the variance of the path flows under recommendations. The approximation yields a tractable single-stage mixed integer linear formulation, which can be solved efficiently with Benders decomposition in large-scale scenarios. The approximation gap is proved to be bounded from the above. Additional theoretical analysis shows that $\epsilon$-feasibility and $\Gamma$-concentration are strongly connected to expectation and chance constraints in a typical stochastic optimization formulation, respectively. The proposed framework is implemented in a real-world case study using data from an urban rail disruption in the CTA system, and a synthetic case study with varied network sizes and incident locations. 

The main contributions of the paper are threefold: 
\begin{itemize}
    \item The paper proposes a framework with prior path utility and posterior path choice distribution given recommendations to model behavior uncertainty, yielding a stochastic optimization problem with decision-dependent distributions. We propose a single-point approximation method to eliminate the expectation operator by introducing two new concepts: $\epsilon$-feasibility and $\Gamma$-concentration, which control the mean and variance of path flows in the optimization problem. 
    \item The proposed single-point approximation yields a tractable single-stage mixed integer linear formulation, which can be solved efficiently with Benders decomposition.
    We prove that, with $\epsilon$-feasibility and $\Gamma$-concentration, the gap of the proposed approximation is bounded from the above. Additional theoretical analysis shows that $\epsilon$-feasibility and $\Gamma$-concentration are strongly connected to expectation and chance constraints in a typical stochastic optimization formulation, respectively.
\end{itemize} 

The remainder of the paper is organized as follows. A literature review is discussed in Section \ref{sec_liter}. In Section \ref{sec_prob_desc}, we describe the problem conceptually and analytically. Section \ref{sec_form} develops the formulation for the individual path recommendation problem and the modeling of the behavior uncertainty. Section \ref{sec_theory} provides all theoretical analysis of the proposed method. In Section \ref{sec_case_study}, we apply the proposed model on the CTA system as a case study and analyze the results. Finally, we conclude the paper and summarize the main findings in Section \ref{sec_conclusion}.

\section{Literature review}\label{sec_liter}
In this section, we go through previous studies on path recommendations and behavior uncertainties. Specifically, we show how our paper is connected with stochastic optimization with decision-dependent distributions (also known as performative predictions).
\subsection{Personalized recommendations in the transportation field}
Personalized recommendation systems in transportation have garnered significant attention in recent years, aiming to enhance user experience by tailoring travel options to individual preferences \citep{thiengburanathum2016overview, quijano2020recommender}. For example, \citet{zhu2020personalized} explored the use of personalized incentives to promote sustainable travel behaviors, demonstrating that targeted rewards can effectively encourage eco-friendly transportation choices. \citet{song2018personalized} introduced a personalized menu optimization approach with a preference updater, applying it to a Boston case study to improve user satisfaction in transportation services. \citet{wu2022personalized} proposed a conceptual framework for personalized travel mode recommendation in a multimodal system. Similarly, \citet{liu2019hydra} developed a context-aware, multi-modal transportation recommendation system that integrates various travel modes, considering user context to provide optimal travel suggestions. \citet{danaf2019online} explored online discrete choice models, applying them to personalized transportation recommendations to predict user preferences more accurately. Additionally, \citet{lim2021origin} presented an origin-aware next-destination recommendation model utilizing personalized preference attention, enhancing the accuracy of destination predictions based on user history. Collectively, these studies underscore the importance of personalization in transportation systems, leveraging user data and advanced modeling techniques to offer tailored travel recommendations.

\subsection{Path recommendations during disruptions}
Path recommendation is a popular topic in the transportation fields \citep{zhang2024survey}. The main goal is to provide a suitable route for a user based on their travel constraints and preferences, with applications in navigation, deliveries \citep{liu2020integrating}, travel planning, etc. Methods for path recommendation include search-based methods \citep{stentz1994optimal, geisberger2008contraction, delling2011customizable, dai2015personalized,wang2014r3}, probability-based \citep{karaman2011sampling, chen2016learning, qu2019profitable}, constraint-based \citep{qu2014cost, reza2017optimal, lai2018urban, cheng2019taxic}, deep learning-based \citep{huang2020multi, wang2021personalized, bhumika2022marrs, wang2023asnn}, and reinforcement learning-based \citep{ji2020spatio, xia2022efficient, chen2022curl}. 

Most previous studies on path recommendation under incidents are like designing a ``trip planner''. That is, the main objective is to find available routes or the shortest path given an OD pair when the network is interrupted by incidents. For example, \citet{bruglieri2015real} designed a trip planner to find the fastest path in the public transit network during service disruptions based on real-time mobility information. \citet{bohmova2013robust} developed a routing algorithm in urban public transportation to find reliable journeys that are robust for system delays. \citet{roelofsen2018assessing} provided a framework for generating and assessing alternative routes in case of disruptions in urban public transport systems. \citet{fang2024travellm} use the large language models to recommend alternative routes in public transit during network disruptions. To the best of the authors' knowledge, most of the previous studies did not consider minimizing system travel time as the goal for path recommendations. The only exceptions are \citet{mo2023robust} and \citet{dai2024guidance}, where they aim to provide station-based path recommendations to guide passengers in urban rail systems in order to minimize system travel time. Compared to typical path recommendation studies where the routes are selected independently for each person, the recommendation with system travel time minimization needs to co-design the strategy of all individuals simultaneously, which is more challenging.

Providing path recommendations during disruptions is similar to the topic of passenger evacuation under emergencies. The objective of evacuation is usually to minimize the total evacuation time, where the decisions need to consider interaction among individual choices. For example, \citet{abdelgawad2012large} developed an evacuation model with routing and scheduling of subway and bus transit to alleviate congestion during the evacuation of busy urban areas. \citet{wang2019optimization} proposed an optimal bus bridging design method under operational disruptions on a single metro line. \citet{tan2020evacuating} proposes an evacuation model with urban bus networks as alternatives in the case of common metro service disruptions by jointly designing the bus lines and frequencies. However, although these passenger evacuation papers focus on minimizing the system travel time, there are several differences from this paper. First, in our paper, the service disruption is not as severe as the emergency situation. We assume the service will recover after a period of time and passengers are allowed to wait. They do not necessarily need to cancel trips or follow the evacuation plan, as assumed in previous evacuation studies. Second, in this article, we assume the supply-side adjustments (e.g., bus rerouting) are given and are not considered decision variables in the optimization. We focus on providing information to the passengers to better utilize the existing resources/capacities of the system, while many evacuation studies focus on supply-side adjustment (e.g., rerouting and rescheduling). Third, this paper considers passenger heterogeneity and focuses on individual-level path recommendations, while previous evacuation papers simply model passengers as flows. Besides, we also assume that passengers may not follow the recommendation (i.e., behavior uncertainty), which has not been considered in any evacuation paper before.

\subsection{Behavior uncertainty and compliance}
Behavior uncertainty is a well-known challenge in transportation modeling \citep{mahmassani1984uncertainty}. It is one of the major reasons for traffic flow instability and the difficulties in predicting traffic conditions \citep{han2024reducing}. In the context of path recommendation, behavior uncertainties cause noncompliance (i.e., passengers may not follow the recommendations), which brings more challenges in evaluating system travel time and other objective functions that depend on passenger flows. For example, \citet{wong2023understanding},  using survey data, revealed there is an evacuation-reluctant group of people who are not willing to follow mandatory evacuation orders in California wildfire. \citet{wang2022non} considered passengers' noncompliance behavior in a transit-based evacuation pick-up point assignment problem.

Typically, passenger's behavior is modeled using various econometrics approaches \citep{ben1985discrete, train2009discrete, mo2021impacts} or machine learning models \citep{mirchevska2013behavior, wang2020deep}. These models output the probability distribution for the passenger's possible behavior. At the aggregate level, numerous studies have taken demand uncertainty into consideration for decision-making. The applications include ride-sharing \citep{guo2021robust}, bus scheduling \citep{guo2024robust}, transit route planning \citep{yoon2020contextual}, and supply chain management \citep{jung2004simulation}. 

However, at the individual level, the number of studies is limited. The main reason is that individual-level decision-making is usually discrete. So it is challenging to use typical robust optimization to address discrete uncertain variables \citep{subramanyam2021robust}. In terms of stochastic optimization, the number of possible scenarios increases exponentially with the number of individuals in the system. Typical ways to address a large number of scenario combinations include 1) Monte Carlo simulation \citep{niederreiter1992quasi, homem2014monte, bartl2022monte}, 2) scenario reduction by clustering \citep{hu2019new, bertsimas2023optimization, keutchayan2023problem}, and 3) better scenario generation methods, such as optimal quantization \citep{bardou2009optimal} and moment matching \citep{hoyland2003heuristic, mehrotra2013generating}. Though these methods help to deal with a large number of combinations of scenarios, most of them only apply to ``exogenous'' uncertainties (i.e., the distribution of the uncertain parameters is given). In this study, we assume an individual's path choice probabilities are affected by our path recommendation decisions (e.g., they are more likely to choose paths that we recommend). This is known as stochastic optimization with decision-dependent distributions \citep{goel2006class, hellemo2018decision, drusvyatskiy2023stochastic}, or preformative predictions in the machine learning field \citep{mendler2020stochastic, perdomo2020performative}. Specifically, these problems can be formulated as:
\begin{align}\label{eq_so_dep}
    \min_{x \in \mathcal{X}} \mathbb{E}_{\zeta \sim \mathbb{P}(x)} \left[ f(x, \zeta) \right],
\end{align}
where $\zeta$ is the random parameter with distribution $\mathbb{P}$ depending on $x$ and $f(x, \zeta)$ is the objective function. Typical methods for solving these problems include 1) performative stochastic gradient descent (for continuous problems) and 2) applying sample average approximation (SAA) in branch-and-cut (for integer problems). However, $f(x, \zeta)$ in the IPR problem represents the system travel time, which usually has no analytical formulations due to the hard-capacity constraints in transit network loading processes \citep{mo2023robust, mo2023ex}. Therefore, the traditional ways are hard to adopt.

\section{Problem description}\label{sec_prob_desc}
In this section, we describe how the individual path recommendation problem is defined, including network and terminology, decision variables and objectives, as well as the scope and time point for decision-making. 
\subsection{Conceptual description}
\begin{figure}[htb]
\centering
\includegraphics[width = 1.0 \linewidth]{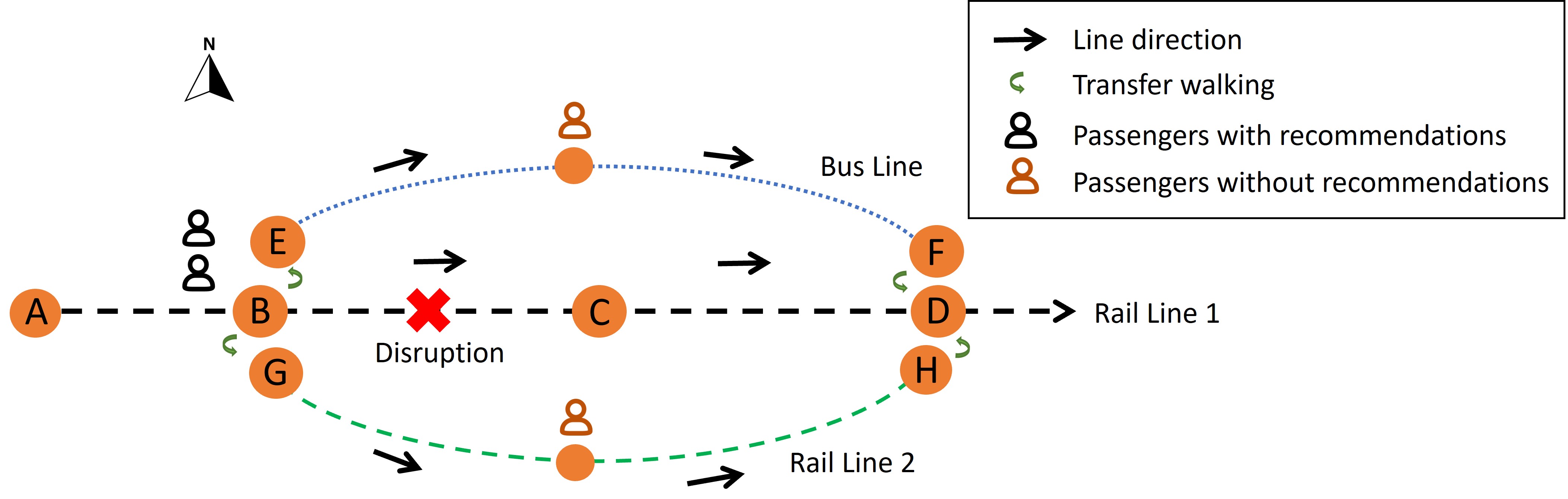}
\caption{Example of the individual path recommendation problem}
\label{fig_example}
\end{figure}

Consider a public transit system with both rail (or subway) and bus services. A public transit network consists of stations (bus or rail stations) and lines (bus or rail lines). A line is defined based on how vehicles (bus or train) are operated. A station is where passengers can board or alight. It is associated with three attributes: station ID, line ID, and line direction. Technically, a station in this study represents a boarding/alighting platform. Passengers can transfer from a station in Line 1 to another station in Line 2. We call the first station of a line where the vehicle departs from the garage the ``terminal'' station. We call a station with multiple lines crossing it as a ``transfer'' station. Figure \ref{fig_example} shows an example of the transit network. In this example, Station A is associated with Rail Line 1 in the east direction. It is also a terminal station. Station B is a transfer station, which can be connected to Station E on the bus line or Station G on Rail Line 2. 

Consider a service disruption in a public transit system. During the disruption, some stations in the incident line (or the whole line) are blocked. Passengers in the blocked vehicles are usually offloaded to the nearest platforms. To respond to the incident, some operating changes are made, such as dispatching shuttle buses, rerouting existing services, short-turning in the incident line, headway adjustment, etc. Assume that all information about the operating changes is available. These changes define a new PT service network and available path sets. Our objective is to develop an individual-based path recommendation model that, when an incident happens, provides a recommended path to every passenger who uses their phones, websites, or electronic boards at stations to enter their origin, destination, and departure time. The recommendation considers the individual's preferences and behavioral histories. Hence, passengers with the same origin, destination, and departure time may get different recommended paths. The overall system aims to minimize the total travel time for all passengers, including passengers in nearby rail or bus lines without incidents (note that these passengers may experience additional crowding due to transfer passengers from the incident line). 

In this example of Figure \ref{fig_example}, Rail Line 1 has an incident between Stations B and C and cannot provide service for a period of time. Both of the two passengers at Station B want to go to Station D. Assuming that they request path recommendations. The alternative paths include using the bus line (blue dashed line), using Rail Line 2 (green dashed line), or waiting for the system to recover (i.e., still using Rail Line 1). Note that using either the bus line or Rail Line 2 will take away capacity from passengers who originally used these two services (i.e., the orange passengers in the figure). Hence, the model should consider the total travel time of all four passengers in the system to design recommendation strategies. Moreover, as mentioned in the introduction, behavior uncertainty needs to be considered. In this example, if we recommend a passenger use a bus line, he/she may not follow the recommendation and choose Rail Line 2 instead.

\subsection{Analytical description}\label{sec_ana_des}
Let us divide the analysis period into several time intervals with equal length $\tau$ (e.g., $\tau  = 5$ min). Let $t$ be the integer time index. $t=1$ is the start of the incident and $t \leq 0$ indicates the time before the incident. Let $\mathcal{P}$ be the set of passengers that will receive path recommendations. We assume $\mathcal{P}$ is known as we can obtain passengers' requests before running the model. Given the revised operation during the incident, let $\mathcal{R}_p$ be the feasible path set for each passenger $p \in \mathcal{P}$. Note that $\mathcal{R}_p$ includes all feasible services that are provided by the PT operator. A path $r \in \mathcal{R}_p$  may be waiting for the system to recover (i.e., using the incident line), or transfer to nearby bus lines, using shuttle services, etc. We do not consider non-PT modes such as TNC or driving for the following reasons: 1) This study aims to design a path recommendation system used by PT operators. The major audience should be all PT users. Considering non-PT modes needs the supply information of all other travel modes, and involves non-PT users (such as the impact of traffic congestion on drivers), which is beyond the scope of this study. Future research may consider a multi-modal path recommendation system. 2) Passengers using non-PT modes can be simply treated as demand reduction for the PT system. So their impact on the PT system can still be captured. 

Given a passenger $p \in \mathcal{P}$, we aim to determine $x_{p,r}$ for each $p$, where  $x_{p,r}$ indicates whether path $r \in \mathcal{R}_p$ is recommended to passenger $p$ or not. Assume only one path is recommended to each passenger, we have
\begin{align}
    \sum_{r \in \mathcal{R}_p} x_{p,r} = 1 \quad \forall p \in \Pcal.
\end{align}
Note that we can relax this assumption by designing the recommendation as a ``composition'' including multiple paths or travel times. This generalization is discussed in Appendix \ref{sec_rec_com}. 

The set $\mathcal{P}$ includes passengers with different origins, destinations, and departure times. If an incident ends at $t^{\text{end}}$, the recommendation should consider a time horizon after $t^{\text{end}}$ because there is remaining congestion in the system. Hence, we provide recommendations until time $T^D > t^{\text{end}}$ (e.g., $T^D$ can be one hour after $t^{\text{end}}$). Therefore, the departure times for passenger $p \in \mathcal{P}$ range from $[1,T^D]$ ($T^D$ and $t^{\text{end}}$ are both time indices).

The recommendation model will be solved at $t = 1$ and will generate the recommendation strategies $\boldsymbol{x} = (x_{p,r})_{p \in \mathcal{P}, r \in \mathcal{R}_p}$ for passengers who depart at time $t\in[1, T^D]$. In reality, the model can be implemented in a rolling horizon manner. Specifically, at each time interval $t \geq 1$, we first update the demand and supply information in the system, including new demand estimates, new to-be-recommended passenger set $\mathcal{P}$, newly available path sets $\Rcal_p$, new service routes and frequencies, new incident duration estimates, new onboard passenger estimates, etc. Based on this information, we solve the model to obtain recommendations for passengers with departure time in $[t, T^D]$. But we only implement the recommendation strategies for passengers who depart at the current time $t$. More discussions can be found in Appendix \ref{append_roll}. 

\begin{figure}[htb]
\centering
\includegraphics[width = 1 \linewidth]{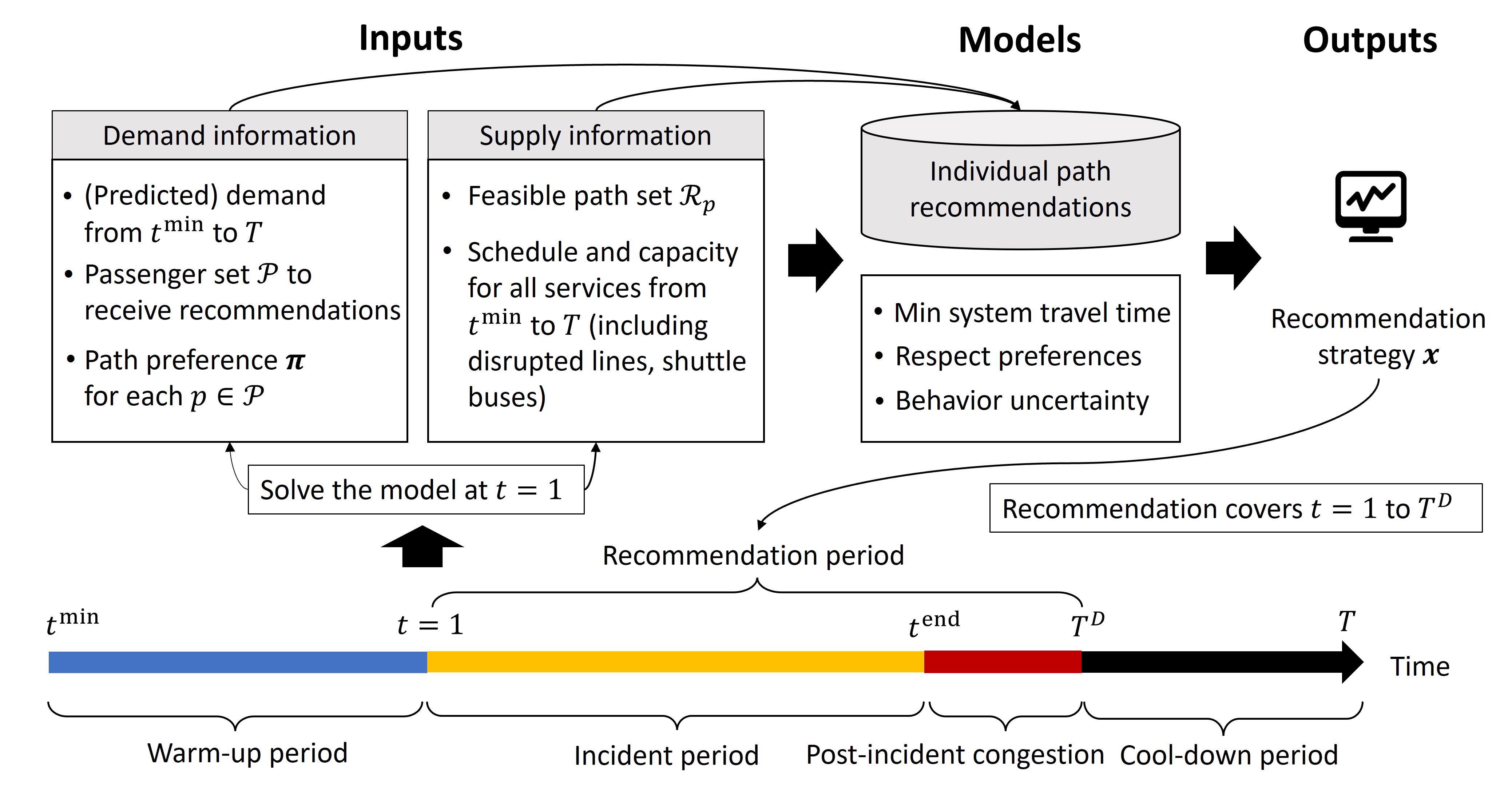}
\caption{Problem description and model framework}
\label{fig_frame}
\end{figure}

Therefore, in the following formulation, we only focus on solving the model at $t = 1$, which is the start of the incident. The whole analysis period includes warm-up and cool-down periods to better estimate the system states (e.g., vehicle loads, passenger travel times, etc.). Therefore, the analysis period is defined as $[t^{\text{min}},T]$, where $t^{\text{min}} < 1$ (time before the incident) and $T > T^D$. For example, $t^{\text{min}}$ and $T$ can be one hour before and after $t=1$ and $T^D$, respectively. And we define all time intervals in the analysis period as $\Tcal = \{\tmin,\tmin+1,..., T\}$. The overall path recommendation framework can be summarized in Figure \ref{fig_frame}. 

It is worth noting that the supply information is assumed to be known in this study. An ideal response strategy to disruptions should decide supply (e.g. rescheduling) and demand (e.g., path recommendation) strategies simultaneously. However, considering both in one decision-making problem is beyond the scope of this study and deserves a separate study. We admit this is a limitation of the paper. 

In reality, only considering path recommendations also has practical implementations. For example, some PT operators may manually decide the supply-side changes (such as dispatching shuttle buses and rerouting) by themselves, as they believe human decisions can incorporate better domain knowledge and have better interpretability. Our model could take their supply-side decision as input and only solve the path recommendation problem. Moreover, in the context of algorithm-driven supply change decisions, our model can be a follow-up subproblem to decide the demand-side strategies, thus making the supply-side and demand-supply strategy design two subproblems. Note that since the model can also be implemented with the rolling horizon method (Appendix \ref{append_roll}), the supply-side information may vary over time according to the newest updates from the operators. 



\section{Formulations}\label{sec_form}
In this section, we elaborate on the detailed formulation of the individual path recommendation model. Section \ref{sec_behavior_unc} describes how passengers' behavior uncertainties (i.e., non-compliance to recommendation) are modeled based on a random utility maximization framework. Section \ref{sec_ind_path_model} provides the overall formulation of the individual path recommendation model as a stochastic optimization with decision-dependent distributions. Section \ref{sec_optimal_flow} presents a linear programming model to calculate the system travel time. Section \ref{sec_single} elaborates on the proposed single-point approximation with $\epsilon$-feasibility and $\Gamma$-concentration. Section \ref{sec_bender} shows how the individual path recommendation model can be solved efficiently using Benders decomposition. 


\subsection{Behavior uncertainty}\label{sec_behavior_unc}
\begin{figure}[htb]
\centering
\includegraphics[width = 0.6 \linewidth]{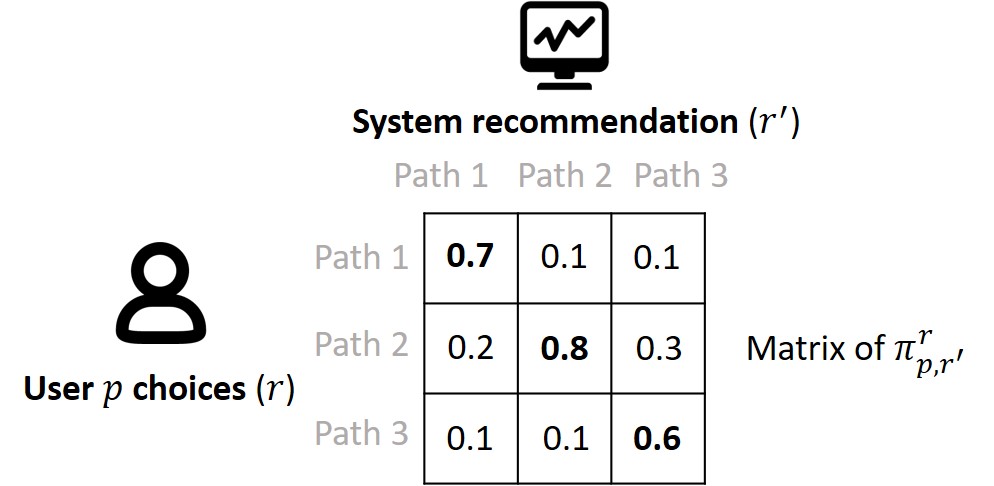}
\caption{Example of conditional path choice probability}
\label{fig_bh_un}
\end{figure}

Consider a passenger $p$ with a path set $\mathcal{R}_p$. Their inherent utility of using path $r$ is denoted as $V_{p}^{r}$. If path $r'$ was recommended, the impact of the recommendation on the utility of path $r$ is denoted as $I_{p,r'}^r$. Hence, his/her overall utility of using path $r$ can be represented as
\begin{align}
    U_{p}^r = V_{p}^{r} + \sum_{r' \in \mathcal{R}_p}x_{p,r'} \cdot I_{p,r'}^r + \xi_{p}^{r} \quad \forall r \in \mathcal{R}_p, \;p \in \mathcal{P} ,
\end{align}
where $\xi_{p}^{r}$ is the random error. $x_{p,r'} = 1$ if passenger $p$ is recommended path $r'$, otherwise $x_{p,r'} = 0$. Let $ \pi_{p,r'}^r$ be the conditional probability that passenger $p$ chooses path $r$ given that the recommended path is $r'$. Assuming a utility-maximizing behavior, we have
\begin{align}
    \pi_{p,r'}^r = \mathbb{P}( V_{p}^{r} + I_{p,r'}^r + \xi_{p}^{r} \geq V_{p}^{r''} +  I_{p,r'}^{r''} + \xi_{p}^{r''}, \; \forall r'' \in \mathcal{R}_p) .
\end{align}
Different assumptions for the distribution of $\xi_{p}^{r}$ can lead to different expressions. For example, if $\xi_{p}^{r}$ are i.i.d. Gumbel distributed, the choice probability reduces to multinomial logit model \citep{train2009discrete, mo2024robust} and we have 
\begin{align}
    \pi_{p,r'}^r = \frac{\exp(V_{p}^{r} +  I_{p,r'}^r) }{\sum_{r''\in \mathcal{R}_p}  \exp(V_{p}^{r''} + I_{p,r'}^{r''} )}.
    \label{eq_pi_prob}
\end{align}
The values of $ V_{p}^{r}$ and $I_{p,r'}^r$ can be calibrated using data from individual-level surveys or smart card data, which deserves separate research. When developing the individual path recommendation model, we assume $\pi_{p,r'}^r$ is known. Figure \ref{fig_bh_un} shows an example of the conditional probability matrix. The specific values assume that paths with recommendations are more likely to be chosen.  

The conditional probability $\pi_{p,r'}^r$ captures the individual's inherent preference for different paths as well as the response to the recommendation system. It varies across individuals and reflects their behavioral uncertainties. This study focuses on design path recommendation systems based on the value of $\pi_{p,r'}^r$.

\subsection{Individual path recommendation}\label{sec_ind_path_model}
Let $\mathbbm{1}_{p,r'}^r$ be the indicator random variable representing whether passenger $p$ \textbf{actually} chooses path $r$ or not given that he/she is recommended path $r'$. By definition, $\mathbbm{1}_{p,r'}^r$ is a Bernoulli random variable with $ \mathbb{E}[\mathbbm{1}_{p,r'}^r] = \pi_{p,r'}^r$ and $ \text{Var}[\mathbbm{1}_{p,r'}^r] = \pi_{p,r'}^r\cdot (1-\pi_{p,r'}^r)$. Consider an OD pair $(u,v)$ and departure time $t$, where $u, v$ are two stations in the transit network. Let $\mathcal{R}^{u,v}$ be the set of feasible paths for OD pair $(u,v)$. Define $Q_t^{u,v,r}$ as the number of passengers in $\Pcal$ with OD pair $(u,v)$ and departure time $t$, who use path $r \in \mathcal{R}^{u,v}$. We have 
\begin{align}
    Q_t^{u,v,r}(\boldsymbol{x}) = \sum_{p\in \mathcal{P}^{u,v}_{t}} \sum_{r' \in \mathcal{R}^{u,v}}x_{p,r'} \cdot \mathbbm{1}_{p,r'}^r ,
\end{align}
where $\mathcal{P}^{u,v}_{t} \subseteq \mathcal{P}$ is the set of passengers with OD pair $(u,v)$ arriving at the system at time interval $t$ that receive path recommendations. Note that $Q_t^{u,v,r}(\boldsymbol{x})$ is also a random variable and it is binomial distributed. Define $\boldsymbol{Q}(\boldsymbol{x}) := (Q_{t}^{u,v,r})_{ t\in\mathcal{T},(u,v,r)\in\mathcal{F}}$, where $\Fcal$ is the set of all triplets $(u,v,r)$ in the system. Denote its distribution as $\mathbb{P}_{\boldsymbol{Q}(\boldsymbol{x})}$, which depends on the recommendation $\boldsymbol{x}$. The objective of this study is to provide a recommendation strategy that minimizes the system travel time (STT), which can be formulated as:
\begin{align}\label{eq_opr}
   (\text{S-IPR}) \quad \min_{\boldsymbol{x}\in\mathcal{X}} \quad & \mathbb{E}_{\boldsymbol{Q}(\boldsymbol{x})} \left[STT\left(\boldsymbol{Q}(\boldsymbol{x})\right)\right],
\end{align}
where $\mathcal{X} := \{\boldsymbol{x}: \sum_{r\in R_p} x_{p,r} = 1, x_{p,r} \in \{0,1\},  \forall p \in \mathcal{P}, r\in \mathcal{R}_p \}$ is the feasible recommendation set. $STT(\boldsymbol{Q}(\boldsymbol{x}))$ is a pseudo function that returns the system travel time of the public transit network given the path flow $\boldsymbol{Q}(\boldsymbol{x})$ (e.g., a simulation or transit network loading model). 

Besides the total system travel time, many recommendation systems also aim to respect passengers' preferences. That is, if possible, a path with high inherent utility $V_{p}^{r}$ should be recommended. Hence the following term is added to the objective function. 
\begin{align}
    \max \sum_{p\in\mathcal{P}} \sum_{r\in\mathcal{R}_p} x_{p,r} \cdot V_{p}^{r} \Longleftrightarrow \min \sum_{p\in\mathcal{P}} \sum_{r\in\mathcal{R}_p} -x_{p,r} \cdot V_{p}^{r}.
\end{align}
The final S-IPR model considering individual preference can be extended as:
\begin{align}\label{eq_ipr}
    \text{(S-IPR-P)} \quad \min_{\boldsymbol{x}\in\mathcal{X}} \quad &  \Psi \cdot \sum_{p\in\mathcal{P}} \sum_{r\in\mathcal{R}_p} -x_{p,r} \cdot V_{p}^{r} +  \mathbb{E}_{\boldsymbol{Q}(\boldsymbol{x})} \left[STT\left(\boldsymbol{Q}(\boldsymbol{x})\right)\right],
\end{align}
where $\Psi$ is a hyperparameter to adjust the scale and balance the trade-off between system efficiency and passenger preferences.

The S-IPR-P problem (\ref{eq_ipr}) is a stochastic optimization with decision-dependent uncertainties \citep{goel2006class, drusvyatskiy2023stochastic}, which is also known as ``performative predictions''. This is because the expectation of the system travel time is taken over the distribution of the path flows $\boldsymbol{Q}(\boldsymbol{x})$, and the distribution of $\boldsymbol{Q}(\boldsymbol{x})$ depends on the decision variable $\boldsymbol{x}$. For normal stochastic optimization, the expectation of the objective function is usually taken over some exogenous parameters (with a fixed distribution that is invariant to decision variables).

The typical ways to solve it are performative stochastic gradient descent (for continuous problems) or applying sample average approximation in branch-and-cut (for integer problems). However, $STT(\cdot)$ usually has no analytical formulation due to the hard-capacity constraints in transit systems \citep{mo2023robust, mo2023ex}. Therefore, the traditional ways are hard to adopt. 

In this study, we propose a new idea to approximate the performative stochastic optimization using a deterministic formulation with two new constraints: ``$\epsilon$-feasibility'' and ``$\Gamma$-concentration''. The idea is that, since the distribution of $\boldsymbol{Q}(\boldsymbol{x})$ is flexible and depending on $\boldsymbol{x}$, we aim to restrict its distribution such that $\mathbb{E}_{\boldsymbol{Q}(\boldsymbol{x})} \left[STT\left(\boldsymbol{Q}(\boldsymbol{x})\right)\right]$ can be approximated by some deterministic values.

\subsection{Solving system travel time as a linear programming}\label{sec_optimal_flow}
One of the main challenges in solving the S-IPR-P problem (\ref{eq_ipr}) is that $STT\left(\cdot\right)$ has no analytical formulation. It is usually obtained by transit network loading or simulation models \citep{mo2020capacity}. \citet{bertsimas2020joint} shows that the network loading process can be formulated as linear programming (LP). In this section, we adapt their formulations for the transit network loading in a public transit system with service disruptions. 

Consider an OD pair $(u,v)$ and departure time $t$. Define $f_t^{u,v,r}$ as the number of passengers \textbf{not in} $\Pcal$ with OD pair $(u,v)$ and departure time $t$, who use path $r \in \mathcal{R}^{u,v}$. Specifically, $Q_t^{u,v,r}(\boldsymbol{x})$ represents the passenger flows that receive recommendations while $f_t^{u,v,r}$ represents those who do not. Hence, the total path flow in $r\in \mathcal{R}^{u,v}$ is $Q_t^{u,v,r} + f_t^{u,v,r}$. Let $d_t^{u,v}$ be the total demand of OD pair $(u,v)$ at time $t$, we have
\begin{align}\label{eq_demand}
  \sum_{r\in\mathcal{R}^{u,v}}  Q_t^{u,v,r}(\boldsymbol{x}) + f_t^{u,v,r} = d_t^{u,v} \quad \forall (u,v) \in \Wcal, t \in \Tcal ,
\end{align}
where $\Wcal$ is the set of all OD pairs. 



\begin{figure}[htb]
\centering
\includegraphics[width = 0.9 \linewidth]{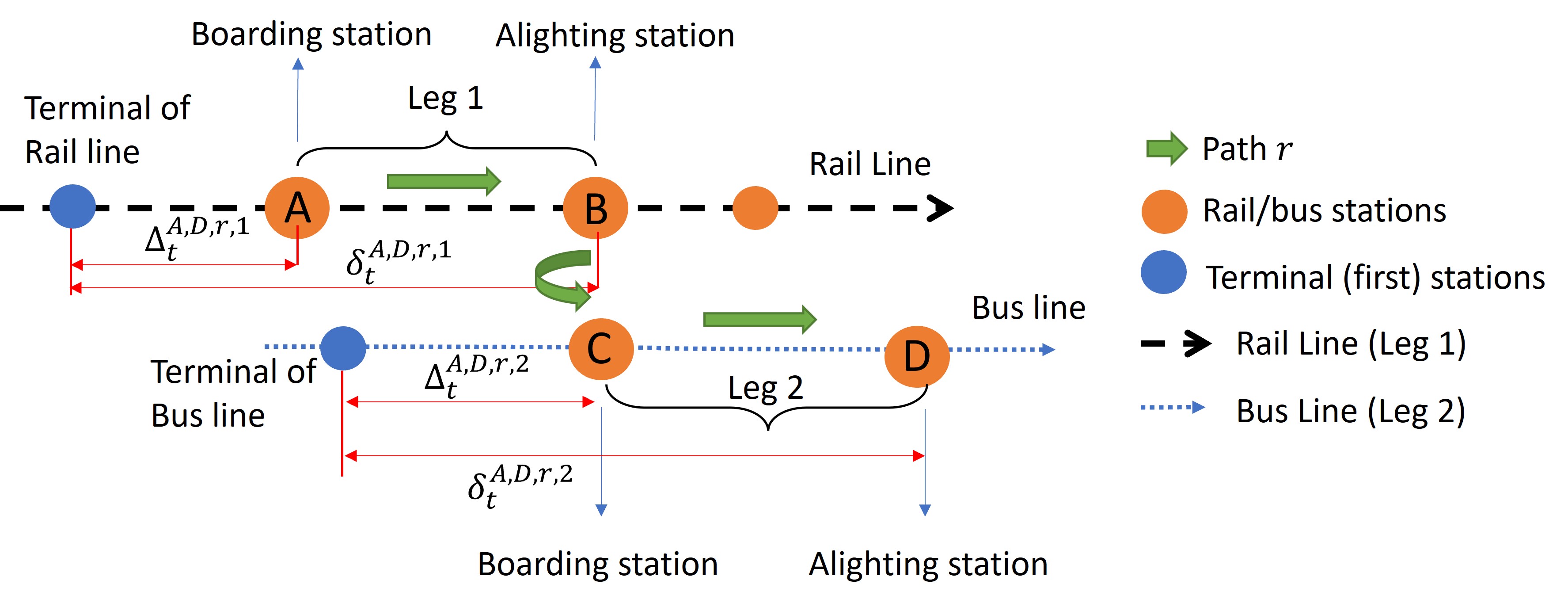}
\caption{Definition of paths and legs}
\label{fig_path}
\end{figure}

Consider a path $r$ for OD pair $(u,v)$. A path may include multiple legs, where each leg is associated with the service in a rail or a bus line. For example, the path $r$ in Figure \ref{fig_path} with origin at station A and destination at station D (indicated by green arrows) has two legs: the first one (A to B) in the rail line and the second in the bus line (C to D). Every leg has a boarding and an alighting station. For example, Leg 1 (resp. 2) in this example has boarding station A (resp. C) and alighting station B (resp. D). Let $\Ical^{u,v,r} = \{1,...,|\Ical^{u,v,r}|\}$ be the set of legs for path $r$. We use a four-element tuple $(u,v,r,i)$ to represent a leg $i$ of path $r$ for OD pair $(u,v)$, where $i\in\Ical^{u,v,r}$. In this example, leg 1 is $(A,D,r,1)$ and leg 2 is $(A,D,r,2)$. 

Let $\Delta^{u,v,r,i}_t$ (resp. $\delta^{u,v,r,i}_t$) be the travel time between the \textbf{terminal} and the \textbf{boarding} (resp. \textbf{alighting}) station of leg $(u,v,r,i)$ for a vehicle \textbf{departing} from the terminal at time $t$, where a terminal is the first departure station of a line. As shown in the example (Figure \ref{fig_path}), $\Delta^{A,D,r,1}_t$ (resp. $\delta^{A,D,r,1}_t$) represent the travel time from the terminal to station A (resp. station B) on the rail line. $\Delta^{A,D,r,2}_t$ (resp. $\delta^{A,D,r,2}_t$) represent the travel time from the terminal to station C (resp. station D) on the bus line.  

Hence, if a vehicle departs at time $t$, its arrival time at the boarding (resp. alighting) station of leg $(u,v,r,i)$ is $t+\Delta^{u,v,r,i}_t$ (resp. $t+\delta^{u,v,r,i}_t$). Then, $\delta^{u,v,r,i}_t - \Delta^{u,v,r,i}_t$ represents the total in-vehicle time of leg $(u,v,r,i)$ for the vehicle. 

Define $z^{u,v,r,i}_{t}$ (decision variable) as the total number of onboard passengers in leg $(u,v,r,i)$ who board a vehicle that \textbf{had departed} from the terminal at time $t$. There are three types of constraints for the network flow description: 1) existing flow constraints, 2) vehicle capacity constraints, and 3) flow conservation constraints.

\textbf{Existing flows constraints:} Although the path recommendations start at time $t=1$, there are passengers that already boarded the vehicles. Ignoring these existing flows may lead to an overestimation of the system's available capacity. To capture the existing onboard flows at $t=1$, we define the set of all onboard flow indices at time $t=1$ as
\begin{align}
    \Omega_1 = \{(u,v,r,i,t):  t+\Delta^{u,v,r,i}_t \leq 1 \leq t+\delta^{u,v,r,i}_t \},
\end{align}
where $ \Omega_1 $ represents the indices of passenger flows who enter the system before the incident starts and have not left the system at the time of the incident (i.e., $t=1$).
The existing flow constraints can be expressed as
\begin{align}
    z_{t}^{u,v,r,i} = \hat{z}_{t}^{u,v,r,i} \quad \forall (u,v,r,i,t) \in \Omega_1,
    \label{eq_existing}
\end{align}
where $\hat{z}_{t}^{u,v,r,i}$ are constants that capture the existing onboard flows when the incident happens. These flows can be directly obtained from a simulation model or real-time passenger counting data.

\textbf{Capacity constraints:} Transit vehicles have limited capacity. Consider a vehicle departing at time $t$ on line $l$ (referred to as vehicle $(l,t)$). We denote its total number of onboard passengers at time $t'$ as $O_{l,t,t'}$. Specifically, $O_{l,t,t'}$ can be expressed as
\begin{align}
O_{l,t,t'}(\boldsymbol{z}) = \sum_{\{(u,v,r,i) \in \texttt{OBLegs}(l,t,t')\}} z_{t}^{u,v,r,i} \quad  \forall l \in \Lcal, \forall t \in \Tcal, t' = t,t+1,...,T_{l,t},
\end{align}
where $T_{l,t}$ is the time index that vehicle $(l,t)$ arrives at the last station of line $l$. $\boldsymbol{z}$ is the decision variable vector defined as $\boldsymbol{z} = (z_{t}^{u,v,r,i})_{ t\in \Tcal,  (u,v,r) \in \mathcal{F}, i \in \Ical^{u,v,r}}$.  $\texttt{OBLegs}(l,t,t')$ is the set of legs for flows that could onboard vehicle $(l,t)$ at time $t'$, defined as
\begin{align}
    \texttt{OBLegs}(l,t,t') = \{
    & (u,v,r,i): \text{Leg } (u,v,r,i) \text{ on line $l$}, \text{and }  t+\Delta_{t}^{u,v,r,i} \leq t' \leq t +\delta_{t}^{u,v,r,i}\}.
    \label{eq_onboard_set}
\end{align}
Equation (\ref{eq_onboard_set}) implies that the legs can be selected based on its associated vehicle arrival time at the boarding and alighting stations. Then the capacity constraint is: 
\begin{align}\label{eq_cap_const}
   O_{l,t,t'}(\boldsymbol{z}) \leq K_{l,t} \quad \forall l \in \Lcal ,t\in\Tcal,t' = t,t+1,...,T_{l,t},
\end{align}
where $K_{l,t}$ is the capacity of the vehicle $(l,t)$. $\Lcal$ is the set of all lines. 

\textbf{Flow conservation constraint:} There are two different flow conservation constraints: 1) flow conservation at origin stations and 2) at transfer stations. To ensure the origin flow conservation, the cumulative number of arrival passengers should be larger than the cumulative number of boarding passengers at an origin at any time (note that we use ``larger than'' because passengers left behind due to capacity constraints are allowed). This indicates that not all arrival passengers can board due to potentially being left behind because of capacity constraints. 

The number of arriving passengers (i.e., demand) for path $(u,v,r)$ at time $t$ is $d_t^{u,v}$. And the number of boarding passengers at the origin station (i.e., $u$) at time $t$ is $z_{t'}^{u,v,r,1}$ (i.e., the first leg) with $t' + \Delta^{u,v,r,1}_{t'} = t$. $t'$ is the vehicle departure time from the terminal. $t' + \Delta^{u,v,r,1}_{t'}$ is the time when the vehicle arrives at the boarding station of leg $(u,v,r,1)$. Therefore, the origin flow conservation constraint can be written as:
\begin{align}\label{eq_flow_cons_1}
    \sum_{\{t': t^{\text{min}}\leq t' + \Delta^{u,v,r,1}_{t'} \leq t \}} z_{t'}^{u,v,r,1} \leq   \sum_{t' = t^{\text{min}}}^t Q_t^{u,v,r}(\boldsymbol{x}) + f_t^{u,v,r} \quad \forall (u,v,r)\in\Fcal,t \in \Tcal.
\end{align}

Now consider the flow conservation at a transfer station. All arrival passengers at a transfer station of a path are the onboard passengers from the last leg. Therefore, we use a similar way to define the transfer flow conservation: the cumulative number of onboard passengers from the last leg should be larger than the cumulative number of boarding passengers at the transfer station. And the number of boarding passengers at the transfer station is simply $z_{t' }^{u,v,r,i}$ with $i\in\Ical^{(u,v,r)}\setminus\{1\}$. Hence, flow conservation constraints at a transfer station are:
\begin{align}\label{eq_flow_cons_2}
\sum_{\{t': t^{\text{min}} \leq t' + \Delta_{t'}^{u,v,r,i}  \leq t \}} z_{t' }^{u,v,r,i}  \leq \sum_{\{t': t^{\text{min}} \leq t' + \delta_{t'}^{u,v,r,i-1}  \leq t\}} z_{t'}^{u,v,r,i-1}\quad \forall (u,v,r)\in\Fcal,i\in\Ical^{(u,v,r)}\setminus\{1\},t \in \Tcal .
\end{align}
Note that $z_{t' }^{u,v,r,i}$ is defined as the onboard passengers for vehicles \textbf{departing} at time $t'$. Therefore, $t' + \delta_{t'}^{u,v,r,i-1}$ is the alighting time for passengers at leg $i-1$ (which is also the transfer demand arrival time at leg $i$ as we assume transfer walk time is within a time interval $\tau$ and is negligible). $t' + \Delta_{t'}^{u,v,r,i}$ is the boarding time for passengers at leg $i$.

The objective is to minimize the total travel time for all passengers in the system. Total travel time can be decomposed into waiting time and in-vehicle time.

\textbf{In-vehicle time:} Total in-vehicle time is simply the onboard flow multiplied by the travel time on each leg:
\begin{align}
IVT(Q_t^{u,v,r}(\boldsymbol{x}), \boldsymbol{z}) = \sum_{(u,v,r)\in\Fcal}\sum_{i\in\Ical^{u,v,r} }\sum_{t\in\Tcal} z^{u,v,r,i}_t \cdot T^{\text{IVT}}_{u,v,r,i,t}
\end{align}
where $T^{\text{IVT}}_{u,v,r,i,t}$ is the in-vehicle time of leg $(u,v,r,i)$ of the vehicle departing at time $t$. 

\textbf{Waiting time:} There are two causes of waiting time: 1) waiting time because of vehicle headways, and 2) waiting time resulting from being left behind. During a specific time interval $t$, all left behind passengers would have a waiting time of $\tau$. All boarding passengers, assuming uniform arrival, have an average waiting time that is half of the time interval (i.e., $\frac{\tau}{2}$).
Therefore, the total waiting time for passengers at station $s$ and time $t$ can be formulated as 
\begin{align}
    WT_{s,t} = \tau(AD_{s,t} + XD_{s,t} - BD_{s,t}) +  \frac{\tau}{2}(BD_{s,t+1} - BD_{s,t}),
\end{align}
where $AD_{s,t}$ represents the \textbf{cumulative arriving demand} at station $s$ \textbf{up to} time $t$, $XD_{s,t}$ represents the \textbf{cumulative transferring demand} at station $s$ \textbf{up to} time $t$, and $BD_{s,t}$ represents the \textbf{cumulative boarded demand} at station $s$ \textbf{up to} time $t$. Hence, $(BD_{s,t+1} - BD_{s,t})$ represents the total number of boarding passengers at time $t$ and station $s$, and $(AD_{s,t} + XD_{s,t} - BD_{s,t})$ represents the total number of left behind passengers at station $s$ and time $t$. Finally, the total system waiting time is 
\begin{align}
WT(\boldsymbol{Q}(\boldsymbol{x}),\boldsymbol{z}) = \sum_{s \in \mathcal{S}} \sum_{t = 1}^T WT_{s,t}
\end{align}
The cumulative arriving demand $AD_{s,t}$ is simply all arriving passengers with origin $s$ up to time $t$:
\begin{align}
AD_{s,t} = \sum_{\{(u,v,r): u = s\}} \sum_{t' = t^{\text{min}}}^t \left(Q_t^{u,v,r}(\boldsymbol{x}) + f_t^{u,v,r} \right) \quad \forall s\in \Scal, t \in \Tcal,
\end{align}
where $\Scal$ is the set of all stations.

The cumulative transferring demand is all passengers alighting at station $s$ from their previous leg $i-1$ for their next leg $i$:
\begin{align}
XD_{s,t} = \sum_{\{(u,v,r,i) \in \texttt{Trans}(s)\}} \sum_{\{ t': t^{\text{min}}\leq t' +\delta^{u,v,r,i-1}_{t'}\leq t\}} z_{t'}^{u,v,r,i-1} \quad \forall t \in \mathcal{T},
\end{align}
where $\texttt{Trans}(s)$ is the set of legs that transfer at station $s$.  

The cumulative boarded demand is all passengers that successfully board a vehicle at station $s$ at time $t$. Define $\texttt{Board}(s)$ as the set of all legs with boarding station $s$, we have
\begin{align}
BD_{s,t} = \sum_{\{(u,v,r,i) \in \texttt{Board}(s)\}} \sum_{ \{t':  t^{\text{min}} \leq t'+\Delta^{u,v,r,i}_{t'} \leq t \}} z_{t'}^{u,v,r,i} \quad \forall t \in \mathcal{T}.
\end{align}

Taking everything into consideration, the total travel time in the system is $WT(\boldsymbol{Q}(\boldsymbol{x}),\boldsymbol{z}) + IVT(\boldsymbol{Q}(\boldsymbol{x}), \boldsymbol{z})$. Assuming everyone in the system would board the first available vehicle, the total system travel time can be obtained by minimizing $WT + IVT$:
\begin{subequations}
\label{eq_of}
\begin{align}
    (\text{LP-STT}) & \quad STT(\boldsymbol{Q}(\boldsymbol{x})) := \min_{\boldsymbol{z}} \quad  WT(\boldsymbol{Q}(\boldsymbol{x}),\boldsymbol{z}) + IVT(\boldsymbol{Q}(\boldsymbol{x}), \boldsymbol{z})  \\
    \text{s.t.} \quad
    & \text{Constraints (\ref{eq_demand}), (\ref{eq_existing}), (\ref{eq_cap_const}), (\ref{eq_flow_cons_1}), and (\ref{eq_flow_cons_2})} ,\label{const_OF1} \\
    & z^{u,v,r,i}_{t} \geq 0 \quad \forall t\in \Tcal,  (u,v,r) \in \mathcal{F}, i \in \Ical^{u,v,r}\label{const_OF2}.
\end{align}
\end{subequations}

\textbf{Incident specification:} Equation (\ref{eq_of}) is a general formulation of the optimal flow problem. Now we will introduce how the incident-specific information is incorporated into this problem. We assume the incident causes a service disruption in a specific line (if only several stations are interrupted, we can separate the line into multiple lines so that the assumption always holds). The service disruption in a line can be seen as the stop of vehicles for a period of time. The vehicle stopping can be captured by the parameters $\Delta^{u,v,r,i}_t$, $\delta^{u,v,r,i}_t$, and $K_{l,t}$. Specifically, a long stop due to an incident can be seen as an increase in travel time from the terminal to downstream stations (i.e., increase in $\Delta^{u,v,r,i}_t$ and $\delta^{u,v,r,i}_t$). Moreover, since there is no vehicle dispatching during the incident, we set $K_{l,t} = 0$ for the corresponding time and line. In this way, we can model the incident without changing the formulation.

\subsection{Eliminating stochasticity with a single-point approximation}\label{sec_single}
With the LP formulation for $STT(\cdot)$, we can reformulate the S-IPR-P problem (\ref{eq_ipr}) as
\begin{align}\label{eq_opr2}
   (\text{S-IPR-P}) \quad \min_{\boldsymbol{x}\in\mathcal{X}} \quad &  \Psi \cdot \sum_{p\in\mathcal{P}} \sum_{r\in\mathcal{R}_p} -x_{p,r} \cdot V_{p}^{r} + \mathbb{E}_{\boldsymbol{Q}(\boldsymbol{x})} \left[\min_{\boldsymbol{z}\in\mathcal{Z}(\boldsymbol{Q}(\boldsymbol{x}))} WT(\boldsymbol{Q}(\boldsymbol{x}),\boldsymbol{z}) + IVT(\boldsymbol{Q}(\boldsymbol{x}), \boldsymbol{z})\right] ,
\end{align}
where $\mathcal{Z}(\boldsymbol{Q}(\boldsymbol{x}))$ is the feasible region for the LP-STT problem (i.e., constraints (\ref{const_OF1}) and (\ref{const_OF2})). The S-IPR-P problem has an outer and inner minimization problem. However, they cannot be directly combined due to the expectation operator. In this study, we propose to remove the expectation operator using a single-point approximation. Specifically, we wish to find a path flow $\boldsymbol{q} \in \mathcal{D}(\boldsymbol{x})$ such that
\begin{align}\label{eq_approx_exp}
    \mathbb{E}_{\boldsymbol{Q}(\boldsymbol{x})} \left[\min_{\boldsymbol{z}\in\mathcal{Z}(\boldsymbol{Q}(\boldsymbol{x}))} WT(\boldsymbol{Q}(\boldsymbol{x}),\boldsymbol{z}) + IVT(\boldsymbol{Q}(\boldsymbol{x}), \boldsymbol{z})\right] \approx \min_{\boldsymbol{q} \in \mathcal{D}(\boldsymbol{x}), \boldsymbol{z}\in\mathcal{Z}(\boldsymbol{q})} WT(\boldsymbol{q},\boldsymbol{z}) + IVT(\boldsymbol{q}, \boldsymbol{z}),
\end{align}
where $\mathcal{D}(\boldsymbol{x})$ is some feasible region that is determined by the recommendation $\boldsymbol{x}$ such that the single-point approximation is reasonable. The intuition for defining $\mathcal{D}(\boldsymbol{x})$ is to make sure that under the recommendation $\boldsymbol{x}$, the distribution of $\boldsymbol{Q}(\boldsymbol{x})$ is concentrate to the mean and the final solution $\boldsymbol{q}^*$ is close to the mean. In this case, we could use a single point to represent the expectation. We define two new concepts, ``$\epsilon$-feasibility'' and ``$\Gamma$-concentration'', that matches with the intuition.

\begin{definition}[\textbf{$\epsilon$-feasible flows}] 
A flow $q_{t}^{u,v,r}$ is $\epsilon$-feasible if and only if
\begin{align}
  \left| q_{t}^{u,v,r} -  {\mu}_{t}^{u,v,r}(\boldsymbol{x}) \right| \leq \epsilon_{t}^{u,v,r}, \quad \forall  (u,v,r) \in \mathcal{F}, t = t^{\text{min}},...,T ,
  \label{eq_eps_feasible}
\end{align}
where \begin{align}
   {\mu}_{t}^{u,v,r}(\boldsymbol{x}) := \mathbb{E}\left[ Q_t^{u,v,r} \right] =  \sum_{p\in \mathcal{P}^{u,v}_{t}} \sum_{r' \in \mathcal{R}^{u,v}}x_{p,r'} \cdot \pi_{p,r'}^r,
\end{align}
and $\epsilon_{t}^{u,v,r}$ is a small positive constant.
\end{definition}

If $\boldsymbol{q}$ is $\epsilon$-feasible, it means $\boldsymbol{q}$ is close to the expectation of the actual flow $\boldsymbol{Q}(\boldsymbol{x})$ under recommendation strategy $\boldsymbol{x}$. $\epsilon$-feasibility results in a direct linear constraint on $\boldsymbol{q}$. 

\begin{definition}[\textbf{$\Gamma$-concentrated flows}] 
A flow $q_{t}^{u,v,r}$ is $\Gamma$-concentrated if and only if it is $\epsilon$-feasible and for any constant $a > \epsilon_{t}^{u,v,r}$, we have
\begin{align}
    \mathbb{P}\left[\left|  Q_t^{u,v,r}  - {q}_{t}^{u,v,r} \right| \geq a \right] \leq \left(\frac{\Gamma_{t}^{u,v,r}}{a-\epsilon_{t}^{u,v,r}}\right)^2 \quad \forall  (u,v,r) \in \mathcal{F}, t = t^{\text{min}},...,T ,
    \label{eq_gamma_concen}
\end{align}
where $\Gamma_{t}^{u,v,r}$ is a small positive constant. 
\end{definition}

If $\boldsymbol{q}$ is $\Gamma$-concentrated, it means that the probability that $Q_t^{u,v,r}$  and ${q}_{t}^{u,v,r}$ are very different (i.e., with difference greater than $a$) is bounded from above, suggesting that $Q_t^{u,v,r}$ is concentrated around $q_t^{u,v,r}$.

Notice that we can modify the definition of $\Gamma$-concentration using the following Proposition:
\begin{proposition}\label{prop_gamma_linear}
The $\Gamma$-concentration inequality (\ref{eq_gamma_concen}) holds if the variance of $Q_t^{u,v,r}$ is bounded from above by $(\Gamma_{t}^{u,v,r})^2$. Mathematically:
\begin{align}
   \text{Var}\left[ Q_t^{u,v,r}  \right]  =  \sum_{p\in \mathcal{P}^{u,v}_{t}} \sum_{r' \in \mathcal{R}^{u,v}}x_{p,r'} \cdot \pi_{p,r'}^r \cdot (1-\pi_{p,r'}^r) \leq (\Gamma_{t}^{u,v,r})^2 .
\end{align}
\end{proposition}
The variance formulation is based on the fact that $(x_{p,r'})^2 = x_{p,r'}$ and  $\text{Cov}[\mathbbm{1}_{p,r'}^r,\mathbbm{1}_{p,r''}^r] = 0$ if $r' \neq r''$. The proof of Proportion \ref{prop_gamma_linear} is a direct result of Chebyshev's inequality and the triangle inequality (see Appendix \ref{proof_prop_gamma_conc}).

Then, we can define $\mathcal{D}(\boldsymbol{x})$ as the set of path flows that are both $\epsilon$-feasible and $\Gamma$-concentrated:
\begin{align}
\mathcal{D}(\boldsymbol{x}) := \Big\{&\boldsymbol{q}\geq 0 : \notag \\
 &q_{t}^{u,v,r} \geq \sum_{p\in \mathcal{P}^{u,v}_{t}}  \sum_{r' \in \mathcal{R}^{u,v}}x_{p,r'} \cdot \pi_{p,r'}^r  - \epsilon_{t}^{u,v,r}, \quad \forall (u,v,r) \in \mathcal{F}, t\in\mathcal{T} \label{eq_feasibility1}\\
 & q_{t}^{u,v,r} \leq  \sum_{p\in \mathcal{P}^{u,v}_{t}} \sum_{r' \in \mathcal{R}^{u,v}}x_{p,r'} \cdot \pi_{p,r'}^r + \epsilon_{t}^{u,v,r}, \quad \forall (u,v,r) \in \mathcal{F}, t\in\mathcal{T} \label{eq_feasibility2}
 \\
 &\sum_{p\in \mathcal{P}^{u,v}_{t}} \sum_{r' \in \mathcal{R}^{u,v}}x_{p,r'} \cdot \pi_{p,r'}^r \cdot (1-\pi_{p,r'}^r) \leq (\Gamma_t^{u,v,r})^2,
 \quad \forall (u,v,r)\mathcal{F}, t\in\mathcal{T} \Big \} .\label{eq_conce}
\end{align}
Constraints (\ref{eq_feasibility1}), (\ref{eq_feasibility2}), and (\ref{eq_conce}) are both linear, which keeps the tractability of the formulation. 

With the approximation of (\ref{eq_approx_exp}), we can drop the expectation operator, and combine the inner and outer minimization problem. Then the approximated individual path recommendation problem considering preference (A-IPR-P) can be reformulated as
\begin{align}\label{eq_aipr}
    \text{(A-IPR-P)} \quad \min_{\boldsymbol{x}\in\mathcal{X}, \boldsymbol{q} \in \mathcal{D}(\boldsymbol{x}), \boldsymbol{z}\in\mathcal{Z}(\boldsymbol{q})} - \Psi \cdot \sum_{p\in\mathcal{P}} \sum_{r\in\mathcal{R}_p}x_{p,r} \cdot V_{p}^{r} + 
 WT(\boldsymbol{q},\boldsymbol{z}) + IVT(\boldsymbol{q}, \boldsymbol{z}) .
\end{align}
The quality of the approximation and other theoretical analysis regarding $\epsilon$-feasibility and $\Gamma$-concentration will be discussed in Section \ref{sec_theory}. For now, we just need to know that a single-point approximation is proposed to transform the original stochastic optimization with decision-dependent distributions to a tractable mixed-integer linear programming (MILP).

For a better understanding of the $\epsilon$-feasibility and $\Gamma$-concentration, let us consider a simple example where there is only one path flow $Q$ and the system travel time is $STT({Q}) = \alpha\cdot Q$. In this case, $\mathbb{E}_{Q}[STT({Q})]=\alpha\cdot\mathbb{E}[Q]=STT(\mathbb{E}[Q])]$. This shows the benefit of $\epsilon$-feasibility: if $q$ is $\epsilon$-feasible and $\epsilon=0$, we would have $\mathbb{E}_{Q}[STT({Q})] =STT(q)$, the single-point approximation becomes exact. Hence, the $\epsilon$-feasibility helps to make the approximated system travel time close to the mean. 

In terms of $\Gamma$-concentration, let us consider an extreme scenario where $Q$ is uniformly distributed in $[0,1]$ under the optimal recommendation strategy $\boldsymbol{x}^*$. In this case, the distribution of $STT(Q)$ will be uniformly in $[0,\alpha]$ (Figure \ref{fig_ill_uniform}). Then, $\boldsymbol{x}^*$ becomes meaningless because, under this recommendation, the actual system travel time has too many variations. However, if we add $\Gamma$-concentration to $Q$ (i.e, $\text{Var}[Q]\leq\Gamma^2$), we would have $\text{Var}[STT(Q)] \leq (\alpha\cdot\Gamma)^2$, the distribution of $\text{Var}[STT(Q)]$ will be concentrated (Figure \ref{fig_ill_conc}) toward the mean value, making the recommendation of minimizing the mean value meaningful.

\begin{figure}[htbp]
\centering

\subfloat[Uniform\label{fig_ill_uniform}]{
\begin{tikzpicture}
  \begin{axis}[
    ymin=0, ymax=0.2,
    xmin=0, xmax=11,
    axis lines=left,
    xlabel=$STT(Q)$, 
    ylabel=Probability density,
    xtick={0, 5, 10},
    xticklabels={$0$, $\alpha/2$, $\alpha$},
    ytick=\empty,
    height=5cm,
    width=0.45\textwidth,
  ]
    \addplot [
        domain=0:10, 
        samples=2, 
        fill=blue!20, 
        draw=none
    ] {1/10} \closedcycle;

    \addplot[domain=0:10, samples=2, very thick, blue] {1/10};
    \addplot[domain=-1:0, samples=2, dashed, gray] {0};
    \addplot[domain=10:11, samples=2, dashed, gray] {0};

    \node at (axis cs:5,0.11) [anchor=south] {No concentration constraints};
  \end{axis}
\end{tikzpicture}
}
\hfill
\subfloat[Concentrated\label{fig_ill_conc}]{
\begin{tikzpicture}
  \begin{axis}[
    ymin=0, ymax=0.5,
    xmin=0, xmax=11,
    domain=0:10, 
    samples=100,
    axis lines=left,
    xlabel=$STT(Q)$,
    ylabel=Probability density,
    xtick={0, 5, 10},
    xticklabels={$0$, $\alpha/2$, $\alpha$},
    ytick=\empty,
    height=5cm,
    width=0.45\textwidth,
  ]
    \addplot [
        domain=0:10,
        samples=100,
        fill=red!20,
        opacity=0.7
    ] {1/(sqrt(2*pi*1)) * exp(-(x - 5)^2 / (2*1))} \closedcycle;
    \addplot[very thick, red] {1/(sqrt(2*pi*1)) * exp(-(x - 5)^2 / (2*1))};
    \node at (axis cs:5,0.42) [anchor=south] {$\text{Var}[STT(Q)] \leq (\alpha\cdot\Gamma)^2$};
  \end{axis}
\end{tikzpicture}
}
\caption{Illustration for the impact $\Gamma$-concentration constraints on the system travel time distribution}
\label{fig:uniform_normal_shaded}
\end{figure}

\subsection{Solving the problem by Benders decomposition}\label{sec_bender}
The structure of (\ref{eq_aipr}) allows us to efficiently solve it by Benders decomposition (BD) \citep{benders1962partitioning}. The basic idea of BD is to decompose the problem into a master problem and a subproblem and solve these problems iteratively. The decision variables are divided into difficult variables, which in our case are the binary variables $\boldsymbol{x}$, and a set of easier variables, the continuous $\boldsymbol{q}$ and $\boldsymbol{z}$. At each iteration, the master problem determines one possible leader decision $\boldsymbol{x}$. This solution is used in the subproblem to generate optimality-cuts or feasibility-cuts, which are added to the master problem. In this study, the master problem decides the recommendation strategies, which is a MILP of a smaller scale and can be solved efficiently using existing solvers. The subproblem reduces to the LP-STT problem (\ref{eq_of}) with one more linear constraint (still linear programming). This format makes the BD an appropriate algorithm for the original problem. As this is a conventional method, we put the formulation details in Appendix \ref{append_benders}.

\section{Theoretical analysis}\label{sec_theory}
The key idea under our proposed method is using a single point to approximate the expectation of the system travel time (\ref{eq_approx_exp}). A natural question is how good the approximation is. In Section \ref{sec_stt_diff}, we show that the gap of the approximation is bounded from above with $\epsilon$-feasibility and $\Gamma$-concentration constraints. In Section \ref{sec_two_stage_SO}, we show that this approximation can be seen as a way of approximating the recourse function in a two-stage stochastic optimization. In Section \ref{sec_relation}, we provide further theoretical analysis to show that $\epsilon$-feasibility and $\Gamma$-concentration are related to expectation and chance constraints in a typical stochastic optimization formulation, respectively.

\subsection{Approximation quality analysis}\label{sec_stt_diff}
Let $(\boldsymbol{{q}^*, {z}^*, {x}^*})$ be the optimal solution of our approximated A-IPR-P formulation (\ref{eq_aipr}).
Given the optimal recommendation strategy $\boldsymbol{x}^*$, let the corresponding random variable of path flows be $\boldsymbol{Q}({\boldsymbol{x}^*})$. Let $\mathbb{P}_{\boldsymbol{Q}({\boldsymbol{x}^*})}(\cdot)$ be the probability density function of $\boldsymbol{Q}({\boldsymbol{x}^*})$. The expectation of $STT(\boldsymbol{Q}({\boldsymbol{x}^*}))$ can be expressed as
\begin{align}
    \mathbb{E}_{\boldsymbol{Q}({\boldsymbol{x}^*})}[STT(\boldsymbol{Q}({\boldsymbol{x}^*}))] = \sum_{\hat{\boldsymbol{q}}\in\mathcal{Q}({\boldsymbol{x}^*})} \left[\min_{\boldsymbol{z}\in\mathcal{Z}(\hat{\boldsymbol{q}})} WT(\hat{\boldsymbol{q}},\boldsymbol{z}) + IVT(\hat{\boldsymbol{q}}, \boldsymbol{z})\right] \cdot \mathbb{P}_{\boldsymbol{Q}({\boldsymbol{x}^*})}(\hat{\boldsymbol{q}}) ,
\end{align}
where $\mathcal{Q}({\boldsymbol{x}^*})$ is the support of the random variable $\boldsymbol{Q}({\boldsymbol{x}^*})$, defined as:
\begin{align}
    \mathcal{Q}({\boldsymbol{x}^*}) = \{\hat{\boldsymbol{q}}\geq 0: \hat{q}_{t}^{u,v,r} = \sum_{p\in \mathcal{P}^{u,v}_{t}} \sum_{r' \in \mathcal{R}^{u,v}}x_{p,r'}^* \cdot \hat{{1}}_{p,r'}^r,\; \forall \;\hat{{1}}_{p,r'}^r \in\{0,1\}, \sum_{r\in\mathcal{R}_p}\hat{{1}}_{p,r'}^r = 1,\;\forall(u,v,r)\in\mathcal{F}, t\in\mathcal{T} \}.
\end{align}
It represents the set of all possible values of network flows given recommendation $\boldsymbol{x}^*$. $\hat{{1}}_{p,r'}^r$ is a binary deterministic variable that can either be 0 or 1 (i,e., the possible realization of ${\mathbbm{1}}_{p,r'}^r$).

$\mathbb{E}_{\boldsymbol{Q}({\boldsymbol{x}^*})}[STT(\boldsymbol{Q}({\boldsymbol{x}^*}))]$ is the expected system travel time given recommendation strategy ${\boldsymbol{x}^*}$. However, our model is optimized using (\ref{eq_aipr}), where the minimization is conducted over the model-evaluated system travel time under the path flow $\boldsymbol{q}^*$:
\begin{align}
    STT({\boldsymbol{q}^*}) = \min_{\boldsymbol{z}\in\mathcal{Z}(\boldsymbol{q}^*)}WT(\boldsymbol{q}^*,\boldsymbol{z}) + IVT(\boldsymbol{q}^*, \boldsymbol{z}).
\end{align}
It is worth analyzing the relationship between the model-evaluated system travel time ($STT({\boldsymbol{q}^*})$) and the expected system travel time ($\mathbb{E}_{\boldsymbol{Q}({\boldsymbol{x}^*})}[STT(\boldsymbol{Q}({\boldsymbol{x}^*}))]$). This analysis tells us how well our proposed approach can approximate the real system performance indicator. 

\begin{lemma}\label{lemma_berge}
 $STT({{\boldsymbol{q}}}) =\min_{\boldsymbol{z}\in\mathcal{Z}(\boldsymbol{q})}WT(\boldsymbol{q},\boldsymbol{z}) + IVT(\boldsymbol{q}, \boldsymbol{z})$ is continuous in terms of ${{\boldsymbol{q}}}$ if the set of optimal flows is bounded (i.e., there are a limited number of flow patterns that permits the optimal system travel time). 
\end{lemma}
The proof is a direct implementation of Berge's Maximum Theorem  \citep{sundaram1996first} and is shown in Appendix \ref{app_proof_lemma}. Lemma \ref{lemma_berge} implies that a small change in the path flows only results in small changes in system travel time. Since the system travel time is usually bounded from above given a finite-scale transit network, small flow changes should not yield infinite changes in the system travel time. Hence, $STT({{\boldsymbol{q}}})$ should also have a bounded gradient. Combining the continuity property in Lemma \ref{lemma_berge}, we conclude that $STT({{\boldsymbol{q}}})$ is Lipschitz continuous. That is, there exists a constant $L$ such that, for any network flows $\boldsymbol{q}_1$ and $\boldsymbol{q}_2$, we have 
\begin{align}
    |STT({\boldsymbol{q}}_1) -  STT({\boldsymbol{q}}_2)| \leq L \cdot \norm{\boldsymbol{q_1} - \boldsymbol{q_2}}_{1}  .
\end{align}


\begin{proposition}\label{eq_prop_bound_stt_diff}
Let $(\boldsymbol{{q}^*, {z}^*, {x}^*})$ be the optimal solution of the A-IPR-P problem (\ref{eq_aipr}). $\boldsymbol{Q}({\boldsymbol{x}^*})$ is the random variable of path flows. The difference between the model-evaluated system travel time and the expected system travel time is bounded from above if
\begin{itemize}
    \item the set of optimal flows is bounded (to implement Lemma \ref{lemma_berge}),
    \item the network flows are bounded from above (i.e., there exist $\boldsymbol{q}^\text{Max} < \infty$ such that $\hat{\boldsymbol{q}}\leq \boldsymbol{q}^\text{Max},\;\forall\; \hat{\boldsymbol{q}}\in\mathcal{Q}({\boldsymbol{x}^*})$), and
    \item $|\boldsymbol{q}^* - \mathbb{E}[\boldsymbol{Q}({\boldsymbol{x}^*})]|\leq\boldsymbol{\epsilon}$ ($\epsilon$-feasibility) and $\text{Var}[\boldsymbol{Q}({\boldsymbol{x}^*})] \leq \boldsymbol{\Gamma}^2$ ($\Gamma$-concentration),
\end{itemize}
where $\mathbb{E}[\boldsymbol{Q}({\boldsymbol{x}^*})] := (\mathbb{E}[{Q}_i({\boldsymbol{x}^*})])_{(u,v,r)\in\mathcal{F},t\in\mathcal{T}}$, $\boldsymbol{\epsilon} := (\epsilon_t^{u,v,r})_{(u,v,r)\in\mathcal{F},t\in\mathcal{T}}$, $\boldsymbol{\Gamma} := \left(\Gamma_t^{u,v,r}\right)_{(u,v,r)\in\mathcal{F},t\in\mathcal{T}}$. The bound of the difference is determined by both $\boldsymbol{\epsilon}$ and $\boldsymbol{\Gamma}$. Mathematically,
\begin{align}
    & \big| \mathbb{E}_{\boldsymbol{Q}({\boldsymbol{x}^*})}[STT(\boldsymbol{Q}({\boldsymbol{x}^*}))]- SST(\boldsymbol{q}^*)\big| \leq 2L\cdot\norm{\boldsymbol{\epsilon}}_1  + L\cdot\big(\norm{\mathbb{E}[\boldsymbol{Q}({\boldsymbol{x}^*})]}_1  + \norm{\boldsymbol{q}^{\text{Max}}}_1 + 2\norm{\boldsymbol{\epsilon}}_{1}\big)\cdot \norm{\boldsymbol{\Gamma}}_2^2 .
\end{align}
\end{proposition}

\begin{remark}\label{remark_prop2}
Proposition \ref{eq_prop_bound_stt_diff} shows that even if the model is optimized on a realization of the system travel time (not the expectation), as long as we impose the $\epsilon$-feasibility and $\Gamma$-concentration, the model-evaluated system travel time and the expected system travel will be similar if $\boldsymbol{\epsilon}$ and $\boldsymbol{\Gamma}$ are small. It provides a quality guarantee of the approximation. 
\end{remark}

Remark \ref{remark_prop2} raises a natural question on how we set the values for $\boldsymbol{\epsilon}$ and $\boldsymbol{\Gamma}$. Ideally, we wish them to be as small as possible to have a good approximation quality. Based on the definition of $\epsilon$-feasibility (i.e., (\ref{eq_feasibility1}) and (\ref{eq_feasibility2})), $\boldsymbol{\epsilon}$ can be as small as $0$ without violating any constraints. Therefore, we should set $\boldsymbol{\epsilon}=0$ for better approximation quality. However, when $\boldsymbol{\Gamma}$ is too small, though we have a lower variance and approximation gap, the recommendation strategies will be very restricted (or even infeasible). In that case, the quality of the recommendation may not be good.  In summary, $\boldsymbol{\epsilon}$ should be set as 0. For $\boldsymbol{\Gamma}$, though small values of $\boldsymbol{\Gamma}$ are preferred for the approximation qualities. There are also trade-offs for the solution quality and problem feasibility. We may need to test different settings to select the best hyperparameters.  

As the selection of $\boldsymbol{\Gamma}$ could affect the feasibility of the problem, we propose the following method to set its value. Consider the group of passengers in $\mathcal{P}_t^{u,v}$. The largest and smallest possible variance of the total variance $\sum_{r\in\mathcal{R}^{u,v}}\text{Var}[Q_{t}^{u,v,r}]$ can be obtained as:
\begin{align}
    (\Gamma_t^{u,v,\text{Max}})^{2} :=\max_{\boldsymbol{x}\in\mathcal{X}} \sum_{r\in\mathcal{R}^{u,v}}\text{Var}[Q_{t}^{u,v,r}(\boldsymbol{x})] = \sum_{p\in \mathcal{P}^{u,v}_{t}} \max_{r' \in \mathcal{R}_p} \left( \sum_{r\in\mathcal{R}_p} x_{p,r'} \cdot \pi_{p,r'}^r \cdot (1-\pi_{p,r'}^r) \right),\\
   (\Gamma_t^{u,v,\text{Min}})^{2} :=\min_{\boldsymbol{x}\in\mathcal{X}} \sum_{r\in\mathcal{R}^{u,v}}\text{Var}[Q_{t}^{u,v,r}(\boldsymbol{x})] = \sum_{p\in \mathcal{P}^{u,v}_{t}} \min_{r' \in \mathcal{R}_p} \left( \sum_{r\in\mathcal{R}_p} x_{p,r'} \cdot \pi_{p,r'}^r \cdot (1-\pi_{p,r'}^r) \right).
\end{align}
Therefore, the values of $\boldsymbol{\Gamma}$ should satisfy:
\begin{align}
    \Gamma_t^{u,v,\text{Min}} \leq \sum_{r\in\mathcal{R}^{u,v}}{\Gamma}_t^{u,v,r} \leq  \Gamma_t^{u,v,\text{Max}}, \quad\forall (u,v)\in\mathcal{W}, t\in\mathcal{T}.
\end{align}
Note that both $\Gamma_t^{u,v,\text{Min}}$ and $\Gamma_t^{u,v,\text{Max}}$ can be pre-calculated given $\boldsymbol{\pi}$. Then, we can modify the $\Gamma$-concentration constraints (\ref{eq_conce}) to 
\begin{align} \label{eq_modify_conce}
    \sum_{r\in\mathcal{R}^{u,v}}\sum_{p\in \mathcal{P}^{u,v}_{t}} \sum_{r' \in \mathcal{R}^{u,v}}x_{p,r'} \cdot \pi_{p,r'}^r \cdot (1-\pi_{p,r'}^r) \leq \sum_{r\in\mathcal{R}^{u,v}}(\Gamma_t^{u,v,r})^2,
 \quad \forall (u,v,r)\mathcal{F}, t\in\mathcal{T}.
\end{align}
With the new constraints, we can directly set the value of $\sum_{r\in\mathcal{R}^{u,v}}\Gamma_t^{u,v,r}$ in $\left[\Gamma_t^{u,v,\text{Min}}, \Gamma_t^{u,v,\text{Max}}\right]$. It is worth noting that since $\text{Var}[Q^{u,v,r}_t]\geq 0$ for all $(u,v,r)\in\mathcal{F}$ and $t\in\mathcal{T}$, imposing an upper-bound to the summation ($\sum_{r\in\mathcal{R}^{u,v}}\text{Var}[Q_{t}^{u,v,r}]$) also implies an upper bound of each single element (just the specific value of the upper bound are different). Hence, all the derivations above still hold.

\subsection{Connections to two-stage stochastic optimization}\label{sec_two_stage_SO}
The original S-IPR-P problem (\ref{eq_opr2}) needs to solve both recommendation strategy $\boldsymbol{x}$ and onboarding flow $\boldsymbol{z}$, a possible alternative formulation is a two-stage stochastic optimization, where the first stage is to determine the recommendation and the second stage is to determine $\boldsymbol{z}$. 

A typical way to solve the two-stage stochastic optimization is to construct an approximation $\hat{STT}(\boldsymbol{x})$ for $\mathbb{E}_{\boldsymbol{Q}({\boldsymbol{x}})}[STT(\boldsymbol{Q}({\boldsymbol{x}}))]$. Then we solve 
\begin{align}\label{eq_two_stage_app1}
    \quad \min_{\boldsymbol{x}\in\mathcal{X}} \quad & \Psi \cdot \sum_{p\in\mathcal{P}} \sum_{r\in\mathcal{R}_p} -x_{p,r} \cdot V_{p}^{r} +   \hat{STT}(\boldsymbol{x}) 
\end{align}
for the first stage instead\footnote{The approximation is then updated based on the second-stage solutions}. From this perspective, we can treat $\epsilon$-feasibility and $\Gamma$-concentration as a way of constructing $\hat{STT}(\boldsymbol{x})$, that is,
\begin{align}\label{eq_two_stage_app2}
    \quad  \hat{STT}(\boldsymbol{x}) := \min_{\boldsymbol{q} \in \mathcal{D}(\boldsymbol{x}), \boldsymbol{z}\in\mathcal{Z}(\boldsymbol{q})} WT(\boldsymbol{q},\boldsymbol{z}) + IVT(\boldsymbol{q}, \boldsymbol{z}).
\end{align}
Therefore, combining (\ref{eq_two_stage_app1}) and (\ref{eq_two_stage_app2}) as a one-stage optimization problem yields our A-IPR-P formulation (\ref{eq_aipr}). 


\subsection{Connections to expectation and chance constraints}\label{sec_relation}
Consider the original S-IPR-P formulation (\ref{eq_opr2}), the inner minimization problem (i.e., $\min_{\boldsymbol{z}\in\mathcal{Z}(\boldsymbol{Q}(\boldsymbol{x}))} WT(\boldsymbol{Q}(\boldsymbol{x}),\boldsymbol{z}) + IVT(\boldsymbol{Q}(\boldsymbol{x}), \boldsymbol{z})$) has uncertainties in both objective and constraints. A typical way to solve it is to transform the original constraints into expectation and chance constraints. For the convenience of analysis, we reformulate the inner minimization in the following general format:
\begin{subequations}
\label{eq_opt_with_rd}
\begin{align}
    (\text{Inner-Minimization})  \quad\min_{\boldsymbol{z}} & \quad g(\boldsymbol{Q}(\boldsymbol{x}), \boldsymbol{z}) \\
     \text{s.t.} & \quad h_j(\boldsymbol{Q}(\boldsymbol{x}), \boldsymbol{z}) \leq 0 \quad \forall j\in\mathcal{J}
\end{align}
\end{subequations}
where $g(\boldsymbol{Q}(\boldsymbol{x}), \boldsymbol{z}) := WT(\boldsymbol{Q}(\boldsymbol{x}),\boldsymbol{z}) + IVT(\boldsymbol{Q}(\boldsymbol{x}), \boldsymbol{z})$ is the objective function. $h_j(\cdot)$ is the constraint function. $\mathcal{J}$ is the set of constraints indices. $h_j(\cdot)$ is defined such that 
$\{\boldsymbol{z}: h_j(\boldsymbol{Q}(\boldsymbol{x}), \boldsymbol{z}) \leq 0,  \forall j\in\mathcal{J}\} \Leftrightarrow \mathcal{Z}(\boldsymbol{Q}(\boldsymbol{x}))$. Then, we can also approximate the expected system travel time as
\begin{subequations}
 \label{eq_opt_with_deter}
\begin{align}
   \mathbb{E}_{\boldsymbol{Q}(\boldsymbol{x})} \left[STT\left(\boldsymbol{Q}(\boldsymbol{x})\right)\right] \approx \min_{\boldsymbol{z}} & \quad\mathbb{E}_{\boldsymbol{Q}(\boldsymbol{x})}[ g(\boldsymbol{Q}(\boldsymbol{x}), \boldsymbol{z})] \\
     \text{s.t.} & \quad \mathbb{E}_{\boldsymbol{Q}(\boldsymbol{x})}[h_i(\boldsymbol{Q}(\boldsymbol{x}), \boldsymbol{z})] \leq 0 \quad\forall j \in\mathcal{J} \qquad \text{(Expectation constraints)} \label{eq_exp_const} \\  \;\text{ or/and }\; &\mathbb{P}_{\boldsymbol{Q}(\boldsymbol{x})}[h_j(\boldsymbol{Q}(\boldsymbol{x}), \boldsymbol{z}) \leq 0] \geq {\eta} \quad \forall j \in\mathcal{J} \qquad \text{(Chance constraints)}  \label{eq_perc_const} 
\end{align}
\end{subequations}
where $\eta$ is a predefined parameter for the probability guarantee of the constraints.
However, the formulations in (\ref{eq_opt_with_deter}) are in general hard to solve except that we have the closed-form expressions for $\mathbb{E}_{\boldsymbol{Q}(\boldsymbol{x})}[\cdot]$ and $\mathbb{P}_{\boldsymbol{Q}(\boldsymbol{x})}[\cdot]$ (or using some approximation techniques for the constraints). In this section, we aim to show that $\epsilon$-feasibility and $\Gamma$-concentration are highly related to the expectation and chance constraints, respectively. Therefore, adding them to the constraints would help with the approximation of the expected system travel time.

\begin{proposition}\label{prop_eplison_lowerbound}
Define:
\begin{align}
&STT_{\text{SO}} := \min_{\boldsymbol{z}}\{\mathbb{E}[g(\boldsymbol{Q}(\boldsymbol{x}), \boldsymbol{z})]:\; \mathbb{E}[h_j(\boldsymbol{Q}(\boldsymbol{x}), \boldsymbol{z})] \leq 0, \; \forall j \in\mathcal{J}\} \label{eq_stoch_ep},\\ 
&  STT_{\text{EP}}(\boldsymbol{\epsilon}) :=\min_{\boldsymbol{z}, \boldsymbol{q}}\{g(\boldsymbol{q}, \boldsymbol{z}):\; h_j(\boldsymbol{q}, \boldsymbol{z}) \leq 0,\; \forall j \in\mathcal{J}, \;|{\boldsymbol{q}} - \mathbb{E}[\boldsymbol{Q}(\boldsymbol{x})] | \leq \boldsymbol{\epsilon} \} ,
\end{align}
where $STT_{\text{SO}}$ is the optimal solution of the stochastic optimization problem with expectation constraints. $STT_{\text{EP}}$ is the optimal solution of the proposed approach with the $\epsilon$-feasibility constraint. If $\boldsymbol{\epsilon} = 0$, and $g(\cdot)$ and $h_j(\cdot)$ are both convex functions (corresponding to the convex optimization), we have 
\begin{align}
    STT_{\text{EP}}(\boldsymbol{\epsilon} = 0) \leq STT_{\text{SO}}.
\end{align}
\end{proposition}
The proof is based on Jensen’s inequality and is shown in Appendix \ref{append_proof2}. Proposition \ref{prop_eplison_lowerbound} shows that, the $\epsilon$-feasibility constraint is highly related to the expectation constraints (\ref{eq_exp_const}). Applying $\epsilon$-feasibility constraint will give a lower bound of the original problem with the expectation constraints. Proposition \ref{prop_eplison_lowerbound} is also related to the certainty-equivalent (or mean-field) variant of a stochastic optimization problem. When $\boldsymbol{\epsilon} = 0$, we have $STT_{\text{EP}}(\boldsymbol{\epsilon} = 0) = STT\left(\mathbb{E}\left[\boldsymbol{Q}(\boldsymbol{x})\right]\right)$. That is, we directly use the mean to solve the original stochastic optimization problem.



Now we will show that $\Gamma$-concentration is highly related to the chance constraint. 
\begin{proposition}\label{prop_gamma_app}
If $h_j(\cdot)$ is Lipschitz continuous with $\ell_2$ norm, that is, there exists a positive constant $C$ such that, for all $\boldsymbol{q_1}$, $\boldsymbol{q_2}$, and $\boldsymbol{z}$, $|h_j(\boldsymbol{q_1}, \boldsymbol{z}) - h_j(\boldsymbol{q_2}, \boldsymbol{z})| \leq C \cdot \norm{\boldsymbol{q_1} - \boldsymbol{q_2}}_2$. Define the constraint set of the $\Gamma$-concentration as:
\begin{align}
    \mathcal{H}_{\Gamma}(\boldsymbol{Q}(\boldsymbol{x})) = \{\boldsymbol{z}: \mathbb{E}[h_j(\boldsymbol{Q}(\boldsymbol{x}), \boldsymbol{z})] \leq 0, \;\forall j \in\mathcal{J}, 
      \text{Var}[\boldsymbol{Q}(\boldsymbol{x})]\leq \boldsymbol{\Gamma}^2 \}.
\end{align}
And define the set of a weaker chance constraint as:
\begin{align}
    \mathcal{H}_{C}(\boldsymbol{Q}(\boldsymbol{x})) = \{\boldsymbol{z}: \mathbb{E}[h_j(\boldsymbol{Q}(\boldsymbol{x}), \boldsymbol{z})] \leq 0,\; 
      \mathbb{P}\left[h_j(\boldsymbol{Q}(\boldsymbol{x}), \boldsymbol{z})  \leq \frac{C \cdot \norm{\boldsymbol{\Gamma}}_2}{\sqrt{1-\eta}} \right] \geq \eta, \; \forall j \in\mathcal{J} \}.
\end{align}
Then we have: $\mathcal{H}_{\Gamma}(\boldsymbol{Q}(\boldsymbol{x})) \subseteq \mathcal{H}_{C}(\boldsymbol{Q}(\boldsymbol{x}))$. That is, the $\Gamma$-concentration constraints can deduce a weaker version of the chance constraints. 
\end{proposition}
The proof uses Chebyshev’s inequality and is shown in Appendix \ref{append_proof3}. Proposition \ref{prop_gamma_app} shows that $\Gamma$-concentration constraints are tighter than a weaker version of the chance constraints. When $\boldsymbol{\Gamma}$ is sufficiently small, we would have $\mathbb{P}\left[h_j(\boldsymbol{Q}(\boldsymbol{x}), \boldsymbol{z})  \leq  \frac{C \cdot \norm{\boldsymbol{\Gamma}}_2}{\sqrt{1-\eta}} \right] \approx \mathbb{P}[h_j(\boldsymbol{Q}(\boldsymbol{x}), \boldsymbol{z})  \leq 0]$. In this case, we derive the chance constraints from the $\Gamma$-concentration constraints. 

\section{Actual case study}\label{sec_case_study}
In this section, we implement the proposed algorithm in a real-world disruption case happening in the Chicago Transit Authority (CTA) urban rail system. We first introduce the background of the case study in Section \ref{sec_case_design}, followed by model parameters, benchmark models, and experiment designs. The numerical results are summarized in Section \ref{sec_results}.
\subsection{Case study design}\label{sec_case_design}

\begin{figure}[htb]
\centering
\includegraphics[width = 0.6\linewidth]{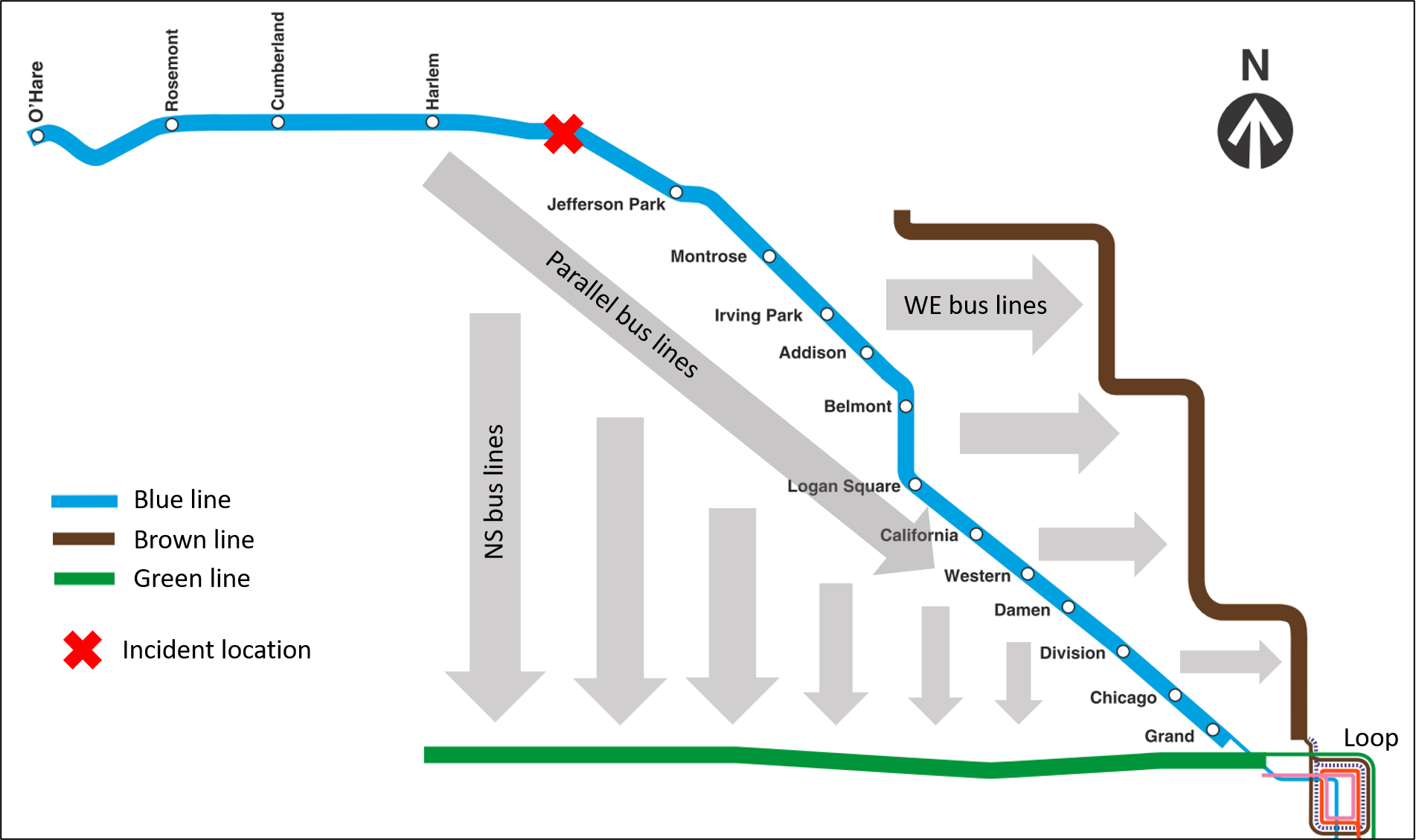}
\caption{Case study network}
\label{fig_incident}
\end{figure}

We consider an actual incident in the Blue Line of the CTA system (Figure \ref{fig_incident}). The incident started at 8:14 AM and ended at 9:13 AM on Feb 1st, 2019 due to infrastructure issues between Harlem and Jefferson Park stations (the red X in the figure) that led to a whole Blue Line suspension. During the disruption (morning hours), the destination for most of the passengers is the ``Loop'' in the CBD area in Chicago. There are four alternative paths to the Loop: 1) using the Blue Line (i.e., waiting for the system to recover), 2) using the parallel bus lines, 3) using the North-South (NS) bus lines to transfer to the Green Line, and 4) using the West-East (WE) bus lines to transfer to the Brown Line. Based on the service structure, the route sets $\Rcal^{(u,v)}$ for each OD pair $(u,v)$ can be constructed.

In the case study, we divide the time into $\tau = 5$ mins equal-length intervals, and focus on solving the problem at $t=1$ (i.e., beginning of the incident). We assume that the set of passengers to receive recommendations ($\Pcal$) consists of all passengers with their intended origins at the Blue Line and destinations in the Loop. A simulation model \citep{mo2020capacity} is used to get the system state up to time $t=1$ (i.e., the incident time 8:14 AM) and generate $ \hat{z}_{t}^{u,v,r,i}$ and $\Omega_1$. The recommendation strategy covers passengers departing between $t=1$ and $T^D=23$, approximately one hour after the end of the incident (9:13 AM). The analysis period is set as $\tmin = -13$ and $T = 34$, approximately one hour before $t=1$ and after $T^D$, providing enough buffer (warm-up and cool-down time) for passengers in $\Pcal$ to finish their trips. As demand and incident duration predictions are out of the scope of this paper, we simply use the actual demand and incident duration for all experiments. Our other work \citep{mo2023robust} proposes to use robust and stochastic optimization to address demand and incident duration uncertainty, respectively.

In terms of the generation of the synthetic conditional probability matrix $\boldsymbol{\pi}$ used for the case study. During the incident, CTA does not provide specific path recommendation information. For every individual, we assume that their actual path choices (referred to as the ``status quo'' choices) reflect their inherent preferences. Appendix \ref{app_infer_status_quo} presents the method and results of inferring passengers' status quo choices during the disruption using smart card data \citep{mo2021inferring}. Around 49\% of the passengers chose to wait while others either took the parallel buses or transferred to the rails. The basic idea is to track their tap-in records when entering the Blue Line and nearby bus routes, and compare them with their historical travel histories to get the transfer information. 

Given the status quo choices, we assume that the ``true'' passenger $p$'s inherent utility for path $r$ is given by 
\begin{equation}
V_{p}^r=\left\{
\begin{aligned}
& 1 + v_p^r& \quad\text{     if $r$ is $p$'s actual path choice} \\
& v_p^r & \quad\text{     otherwise},
\end{aligned}
\right. \quad \forall \; p\in \mathcal{P}, r\in \mathcal{R}_p \label{eq_gen_pre},
\end{equation}
where $v_p^r$ is drawn uniformly from $\mathcal{U}[0,1]$. Equation (\ref{eq_gen_pre}) indicates every path has a random utility $v_p^r$ normalized to 0 $\sim$ 1. The chosen path has an additional utility value of 1. We assume that the impact of the recommendation of $r'$ on the utility of path $r$ is 
\begin{equation}
I_{p,r'}^r=\left\{
\begin{aligned}
& \text{Drawn from }\mathcal{U}[0,5]& \quad\text{     if $r = r'$} \\
& 0 & \quad\text{     otherwise},
\end{aligned}
\right. \quad \forall \; p\in \mathcal{P}, r,r'\in \mathcal{R}_p \label{eq_gen_pre2}.
\end{equation}
Equation (\ref{eq_gen_pre2}) means that the utility of the path recommended (i.e., $r = r'$) has an additional positive impact drawn uniformly from $\mathcal{U}[0,5]$. The utilities of paths not being recommended ($r \neq r'$) do not change. Given (\ref{eq_gen_pre}) and (\ref{eq_gen_pre2}), we can generate the conditional probability $\boldsymbol{\pi}$ using (\ref{eq_pi_prob}). It is worth mentioning that the above assumptions for generating synthetic passenger prior preferences are based on two reasonable principles: 1) Passenger's actual chosen paths have a higher inherent utility. 2) Recommendations of a path can increase its probability of being chosen. 

\subsection{Parameter settings}
The convergence gap threshold for Benders decomposition is set as $1\times 10^{-8}$. The $\epsilon$-feasibility and $\Gamma$-concentration parameters are set with different values to test the performance. Specifically, we set $\epsilon_{t}^{u,v,r} = \epsilon \cdot q^{u,v,r}_t$ for $\epsilon \in \{0, 0.01, 0.03, 0.05, 0.1\}$, indicating different level of deviations on path flows. For $\Gamma$-concentration, we implement the modified constraint in (\ref{eq_modify_conce}), and set $\sum_{r\in\mathcal{R}^{u,v}}({\Gamma}_t^{u,v,r})^2 = \Gamma \cdot \left((\Gamma_t^{u,v,\text{Max}})^2 - (\Gamma_t^{u,v,\text{Min}})^2\right) + (\Gamma_t^{u,v,\text{Min}})^2$ for $\Gamma \in \{0,0.25,0.5,0.75,1\}$. Note that $\Gamma=0$ indicates that the model is forced to choose the lowest variance recommendations, while $\Gamma=1$ means that there are no constraints on the variance.

\subsection{Benchmark models}
There are two benchmark path choice scenarios we use for comparison purposes:

\textbf{Status-quo path choices}. This scenario provides the status quo situation which does not include any recommendations. It represents the worst case. In this scenario, no behavior uncertainty is considered because this is based on the actual path choices realized by passengers.

\textbf{Capacity-based path recommendations}. The capacity-based path recommendations aim to recommend passengers to different paths according to the available capacity of paths. Specifically, for a path in OD pair $(u,v)$ and time $t$, its capacity is the total available capacity of all vehicles passing through the first boarding station of the path during the time period. For example, for a path consisting of an NS bus route and the Green Line, the path capacity is the total available capacity of all buses at the boarding station of the NS bus route during time interval $t$. The available capacity can be obtained from a simulation model using historical demand as the input or using historical passenger counting data. The available capacity for the Blue line (the incident line) depends on modified operations during the incident (i.e., the service suspension is considered). When no vehicles operate in the Blue line during time interval $t$, the path capacity is zero. Let the available capacity of the path $(u,v,r)$ at time $t$ be $C^{u,v,r}_t$. The set of all capacity-based recommendation strategies is defined as:
\begin{align}
    \mathcal{X}^{\text{Cap}} = \left\{\boldsymbol{x}\in\mathcal{X}: \frac{\sum_{p\in\mathcal{P}^{u,v}_t} x_{p,r}}{|\mathcal{P}^{u,v}_t|} = \frac{C^{u,v,r}_t}{\sum_{r'\in\mathcal{R}^{u,v}} C^{u,v,r'}_t},\; \forall (u,v,r)\in\mathcal{F},t\in\mathcal{T}\right\}.
\end{align}
We randomly choose one strategy $\boldsymbol{x}^{\text{Cap}} \in\mathcal{X}^{\text{Cap}}$ to implement, which is equivalent to randomly selecting passengers to recommend paths such that the proportion of passengers being recommended with path $r$ is the same as its available capacity proportion out of all paths in that OD pair $(u,v)$ and time $t$.

\subsection{System travel time evaluation}\label{sec_gen_sys_tt}
Given a recommendation strategy $\boldsymbol{x}$, as mentioned above, the actual system travel time is a random variable because of the passenger behavior uncertainty. To obtain the mean and standard deviation of the system travel time, we generate multiple passenger choice realizations based on $\boldsymbol{\pi}$ and $\boldsymbol{x}$. For each generated passenger choice ($\hat{\mathbbm{1}}_{p,r'}^r$), the realized path flows are
\begin{align}
    \hat{q}_t^{u,v,r} = \sum_{p\in \mathcal{P}^{u,v}_{t}} \sum_{r' \in \mathcal{R}^{u,v}}x_{p,r'} \cdot \hat{\mathbbm{1}}_{p,r'}^r \quad \forall (u,v,r)\in\Fcal , t \in \Tcal.
\end{align}
The corresponding realized system travel time is $STT(\hat{\boldsymbol{q}})=\min_{\boldsymbol{z}\in\mathcal{Z}(\hat{\boldsymbol{q}})} WT(\hat{\boldsymbol{q}},\boldsymbol{z}) + IVT(\hat{\boldsymbol{q}}, \boldsymbol{z})$. This process is repeated with multiple realizations, providing the sample mean and standard deviation of the system travel time under recommendation strategy $\boldsymbol{x}$.

\subsection{Experiment design}
As this paper considers various components (such as individual path recommendations, passengers' path preferences, behavior uncertainty, etc.), it is useful to test different components separately to identify the impact of each one. Hence, we design the following test cases, each one with specific parameter settings to systematically evaluate the impacts of each component. 

\textbf{Model performance compared to benchmark models}. The most straightforward model validation is to evaluate the effect of reducing system travel time. In this test case, we set $\Psi = 0$, meaning that we ignore the passengers' preferences and focus only on minimizing system travel time. The results of this test case are discussed in Section \ref{sec_model_perform} with different values of $\Gamma$ and $\epsilon$.

\textbf{The benefit of considering behavior uncertainty}. In this test case, we evaluate the importance of incorporating behavior uncertainty in the model. The model without behavior uncertainty assumes that passengers take the recommended path. The recommendation strategy is obtained by solving (\ref{eq_ipr}) with $\pi_{p,r'}^r = 1$ if $r=r'$. Similarly, we set $\Psi = 0$. Note that, when we evaluate the recommendation strategy, the behavior uncertainty is still considered in generating the system travel time (see Section \ref{sec_gen_sys_tt}). The results of this test case are shown in Section \ref{sec_bene_unc}

\textbf{Impact of considering passenger preferences}. In all the above tests, $\Psi = 0$ is used, focusing on the system travel time. In this test case, we evaluate the model performance under different values of $\Psi$ in order to assess the impact of considering passenger preferences. The results of this test case are discussed in Section \ref{sec_impact_pax_pref}.

\subsection{Numerical results}\label{sec_results}
In this section, we summarize the experiment results of the proposed model compared with other benchmarks. The computational time and system travel time comparison are shown in Sections \ref{sec_res_compu_time} and \ref{sec_model_perform}. Sections \ref{sec_bene_unc} and \ref{sec_impact_pax_pref} discuss the impact of considering behavior uncertainties and passenger preferences. 
\subsubsection{Model convergence and computational time}\label{sec_res_compu_time}
\begin{figure}[htbp]
\centering
\includegraphics[width = 0.6\linewidth]{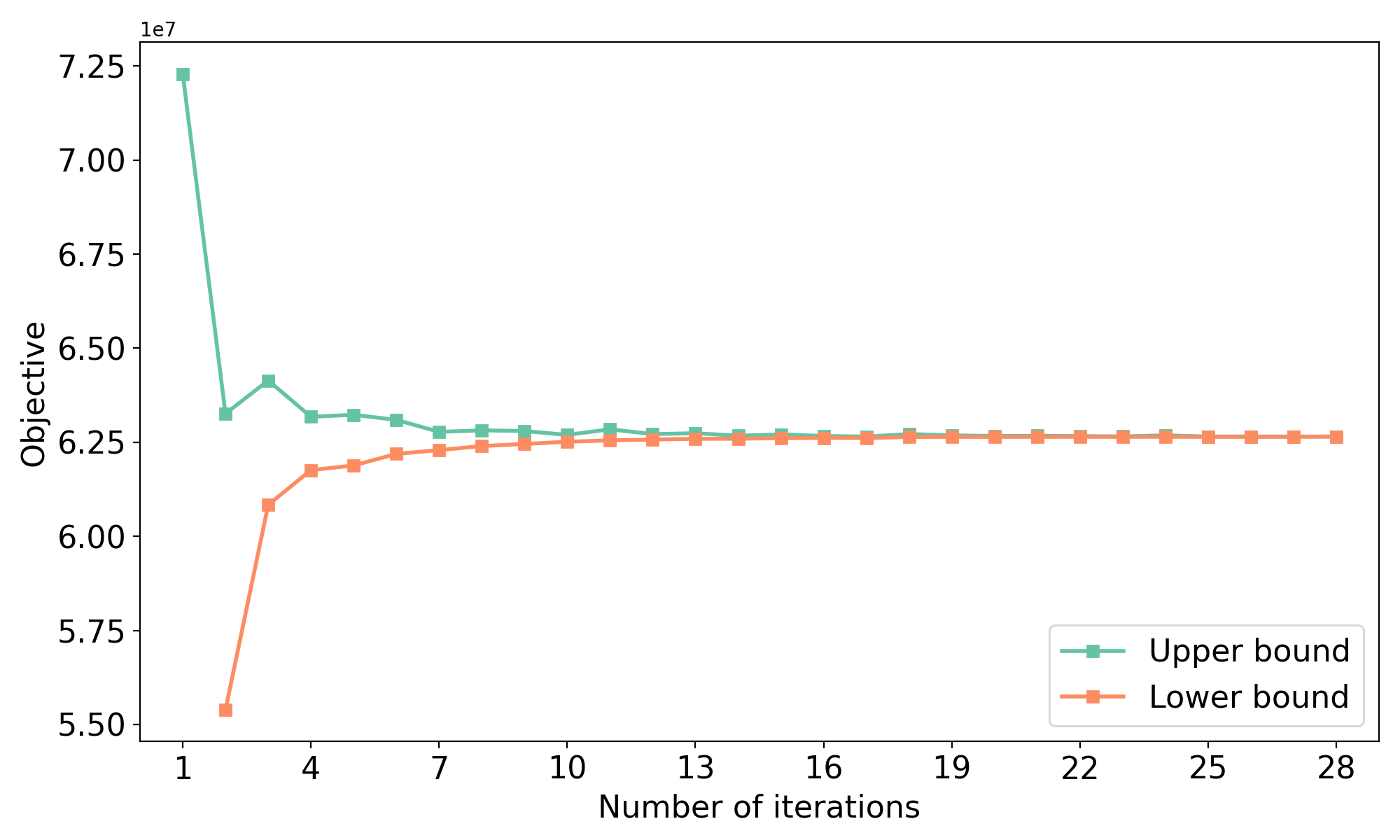}
\caption{Convergence of the Benders decomposition}
\label{fig_BD_con}
\end{figure}

Figure \ref{fig_BD_con} shows the convergence of the BD algorithm. As expected, the lower bound of the model keeps increasing, while the upper bound, after dropping significantly in early iterations, exhibits some fluctuations. The model converges after 28 iterations with a relative gap of less than 1$\times10^{-8}$. The number of optimality cuts was 28 and no feasibility cut was generated. 

Table \ref{tab_time_comp} compares the computational time of the Benders decomposition and off-the-shelf solvers. 
The BD algorithm was implemented using Julia 1.6 with the Gurobi 9.1 solver \citep{gurobi} on a personal computer with the I9-9900K CPU. The total computational time is 17.8 seconds (master problem 8.2 seconds + subproblem 9.6 seconds), which is more efficient than directly using the Mixed integer programming (MIP) solvers, including Gurobi \citep{gurobi}, CPLEX \citep{cplex2009v12}, GLPK (GNU Linear Programming Kit) \citep{makhorin2008glpk}, and CBC (Coin-or branch and cut) \citep{forrest2005cbc}.
\begin{table}[htbp]
\centering
\caption{Computational time comparison}\label{tab_time_comp}
\resizebox{0.6\textwidth}{!}
{
{\renewcommand{\arraystretch}{1} 
\begin{tabular}{@{}lcclcc@{}}
\toprule
Solver                 & CPU time (sec) & Gap     & Solver & CPU time (sec) & Gap     \\ \midrule
BD & 17.8                     & 0.000\% & Gurobi & 55.1                     & 0.000\% \\
CPLEX                  & 65.7                     & 0.000\% & CBC    & 425.4                    & 0.000\% \\
GLPK                   &       562.6                   &   0.000\%      &   &                          &         \\ \bottomrule
\end{tabular}
}
}
\end{table}

\subsubsection{Model performance compared to benchmark models}\label{sec_model_perform}
In this section, we compare the system travel time under the proposed individual path recommendations (without post-adjustment) and two benchmark models. All travel times (except for the status quo that is deterministic) are calculated based on 10 replications using the randomly sampled actual path choices based on the given recommendation (see Section \ref{sec_gen_sys_tt}). We only chose 10 replications because the results are relatively stable and the standard deviations (std.) are small.

\begin{table}[htbp]
\centering
\caption{Average travel time (ATT) comparison for different models}\label{tab_result}
{\renewcommand{\arraystretch}{1} 
\begin{tabular}{lcccc}
\toprule
\multirow{2}{*}{Models} & \multicolumn{2}{c}{ATT (all passengers)} & \multicolumn{2}{c}{ATT (incident line passengers)} \\ \cmidrule(l){2-5} 
                       & Mean (min)               & Std.  (min)                & Mean (min)                       & Std.  (min)                      \\ \midrule
Status quo             & 28.318                      & N.A.                       & 40.255                           & N.A.                            \\
Capacity-based         & 27.609 (-2.5\%)                   & 0.033                      & 33.848 (-15.9\%)                        & 0.165                           \\ \midrule
\cellcolor{gray!25}IPR model ($\epsilon=0.0, \Gamma=1.0$)            & \cellcolor{gray!25}26.454 (-6.6\%)                   & \cellcolor{gray!25}0.019                      & \cellcolor{gray!25}32.526 (-19.2\%)                        & \cellcolor{gray!25}0.204                          \\ 
IPR model ($\epsilon=0.0, \Gamma=0.75$)            & 26.461 (-6.5\%)                   & 0.019                      & 32.563 (-19.1\%)                        & 0.183                           \\ 
IPR model ($\epsilon=0.0, \Gamma=0.5$)            & 26.492 (-6.4\%)                   & 0.019                      & 32.776 (-18.6\%)                        & 0.167                           \\ 
IPR model ($\epsilon=0.0, \Gamma=0.25$)            & 26.649 (-5.9\%)                   & 0.017                      & 33.768 (-16.1\%)                        & 0.145                           \\ 
IPR model ($\epsilon=0.0, \Gamma=0.0$)            & 28.315 (-0.0\%)                   & 0.015                      & 39.950 (-0.76\%)                        & 0.126                           \\  \midrule
IPR model ($\epsilon=0.01, \Gamma=1.0$)            & 26.455 (-6.6\%)                   & 0.019                      & 32.529 (-19.2\%)                        & 0.199                           \\ 
IPR model ($\epsilon=0.03, \Gamma=1.0$)            & 26.459 (-6.6\%)                   & 0.018                      & 32.539 (-19.2\%)                        & 0.191                           \\ 
IPR model ($\epsilon=0.05, \Gamma=1.0$)            & 26.460 (-6.6\%)                   & 0.017                      & 32.555 (-19.1\%)                        & 0.187                           \\ 
IPR model ($\epsilon=0.10, \Gamma=1.0$)            & 26.467 (-6.5\%)                   & 0.017                      & 32.574 (-19.1\%)                        & 0.178                           \\ \midrule
IPR model ($\epsilon=0.01, \Gamma=0.5$)            & 26.493 (-6.4\%)                   & 0.019                      & 32.779 (-18.6\%)                        & 0.165                           \\ 
IPR model ($\epsilon=0.03, \Gamma=0.5$)            & 26.495 (-6.4\%)                   & 0.018                      & 32.786 (-18.6\%)                        & 0.163                           \\ 
IPR model ($\epsilon=0.05, \Gamma=0.5$)            & 26.498 (-6.4\%)                   & 0.017                      & 32.791 (-18.5\%)                        & 0.147                           \\ 
IPR model ($\epsilon=0.10, \Gamma=0.5$)            & 26.504 (-6.4\%)                   & 0.017                      & 32.805 (-18.5\%)                        & 0.128                           \\ \midrule
IPR model ($\epsilon=0.01, \Gamma=0.25$)            & 26.653 (-5.9\%)                   & 0.017                      & 33.771 (-16.1\%)                        & 0.147                           \\ 
IPR model ($\epsilon=0.03, \Gamma=0.25$)            & 26.655 (-5.9\%)                   & 0.017                      & 33.801 (-16.0\%)                        & 0.144                           \\ 
IPR model ($\epsilon=0.05, \Gamma=0.25$)            & 26.677 (-5.8\%)                   & 0.017                      & 33.812 (-16.0\%)                        & 0.137                           \\ 
IPR model ($\epsilon=0.10, \Gamma=0.25$)            & 26.691 (-5.7\%)                   & 0.017                      & 33.830 (-16.0\%)                        & 0.141                
          \\ 
\bottomrule
\multicolumn{5}{l}{\begin{tabular}[c]{@{}l@{}} \small Numbers in parentheses represent percentage travel time reduction compared to the status quo \\
\small The best model is highlighted in gray\end{tabular}}

\end{tabular}
}
\end{table}

Table \ref{tab_result} shows that the proposed model (IPR) significantly reduces the average travel time (ATT) in the system compared to the status quo for most of the hyperparameter settings. The best model with $\epsilon=0.0, \Gamma=1.0$ has a 6.6\% reduction in travel times of all passengers in the system. And for passengers in the incident line (i.e., passengers who received the recommendation, $\Pcal$),  the average travel time reduction is 19.2\%. Our model also outperforms the capacity-based benchmark path recommendation strategy, which reduces the travel time of all passengers by 2.5\% and incident line passengers by 15.9\%. We also compare individual-level travel time savings. Results are elaborated in Appendix \ref{append_ind_tt}.

Comparing different parameter settings, with the decreasing values in $\Gamma$, the model has a tighter constraint on the variance. We also observe decreasing trends in the standard deviations. However, the lower value in $\Gamma$ will make the recommendation strategies more restricted. Hence, the quality of the recommendations is worse with higher system travel times (for all passengers and incident line passengers). Specifically, when $\Gamma=0$, the model is forced to choose the smallest variance recommendations, resulting in significantly bad solutions (even worse than the capacity-based recommendations). In terms of $\epsilon$, as we analyzed in Section \ref{sec_stt_diff}, the best value should be 0 since it provides a better approximation gap and does not impose additional constraints on the model. The experiment results also validate this under various values of $\Gamma$.

\subsubsection{Benefits of considering behavior uncertainty}\label{sec_bene_unc}
In this section, we aim to compare the model with and without considering the behavior uncertainty. The model without behavior uncertainty assumes that all passengers follow the recommended path when designing the recommendation (but they may not in reality). 

Table \ref{tab_result2} shows the comparison of average travel time for the two models. The IPR model with behavior uncertainty (BU) uses $\epsilon=0, \Gamma=1$ as they are the best hyperparameters. As expected, considering behavior uncertainty in the path recommendation design achieves a smaller travel time for all passengers and incident line passengers. Note that, though the 0.94\% reduction (around 15 seconds saving per passenger) is relatively small, considering the large number of passengers in the system, the total travel time savings are still significant.     

\begin{table}[htb]
\centering
\caption{Average travel time (ATT) comparison with and without behavior uncertainty (BU)}\label{tab_result2}
{\renewcommand{\arraystretch}{1} 
\begin{tabular}{lcccc}
\toprule
\multirow{2}{*}{Models} & \multicolumn{2}{c}{ATT (all passengers)} & \multicolumn{2}{c}{ATT (incident line passengers)} \\ \cmidrule(l){2-5} 
& Mean (min)                & Std. (min)                  & Mean (min)              & Std. (min)                       \\ \midrule
IPR model (w.o. BU)         & 26.706            & 0.026               & 32.852         & 0.122                  \\
IPR model (w. BU)        & 26.454 (-0.94\%)        & 0.019             & 32.526 (-1.0\%)         & 0.204             \\ \bottomrule
\multicolumn{5}{l}{\begin{tabular}[c]{@{}l@{}} \small Numbers in parentheses represent percentage travel time reduction compared to the IPR \\\small  model w.o. BU \end{tabular}}
\end{tabular}
}
\end{table}

\subsubsection{Impact of respecting passenger's prior preferences}\label{sec_impact_pax_pref}
In this section, we evaluate the impact of different values of $\Psi$ in terms of respecting passengers' prior preferences. Besides the system travel time, we also evaluate the total utility, defined as the sum of the prior utilities of the recommended path:
\begin{align}
    TU(\boldsymbol{x}) = \sum_{p\in\mathcal{P}} \sum_{r\in\mathcal{R}_p} x_{p,r} \cdot V_{p,r}. 
\end{align}
Note that the maximum value of $TU(\boldsymbol{x})$ is achieved when every passenger is recommended with their preferred path (i.e., the path with the highest prior utility, $V_{p,r}$). Denote this maximum value as $TU^{\text{max}}$. The relative ratio of total utility, $\frac{TU(\boldsymbol{x})}{TU^{\text{max}}}$, represents the fraction of the total (prior) utility that the recommendation has achieved.

Another indicator is the number of passengers recommended with their preferred path (denoted as $NP(\boldsymbol{x})$). Similarly, we also define the proportion of passengers recommended with their preferred path (i.e., $\frac{NP(\boldsymbol{x})}{|\Pcal|}$, where $|\Pcal| = 5,827$ in the case study). 

\begin{figure}[htb]
\centering
\subfloat[Average travel time]{\includegraphics[width=0.5\textwidth]{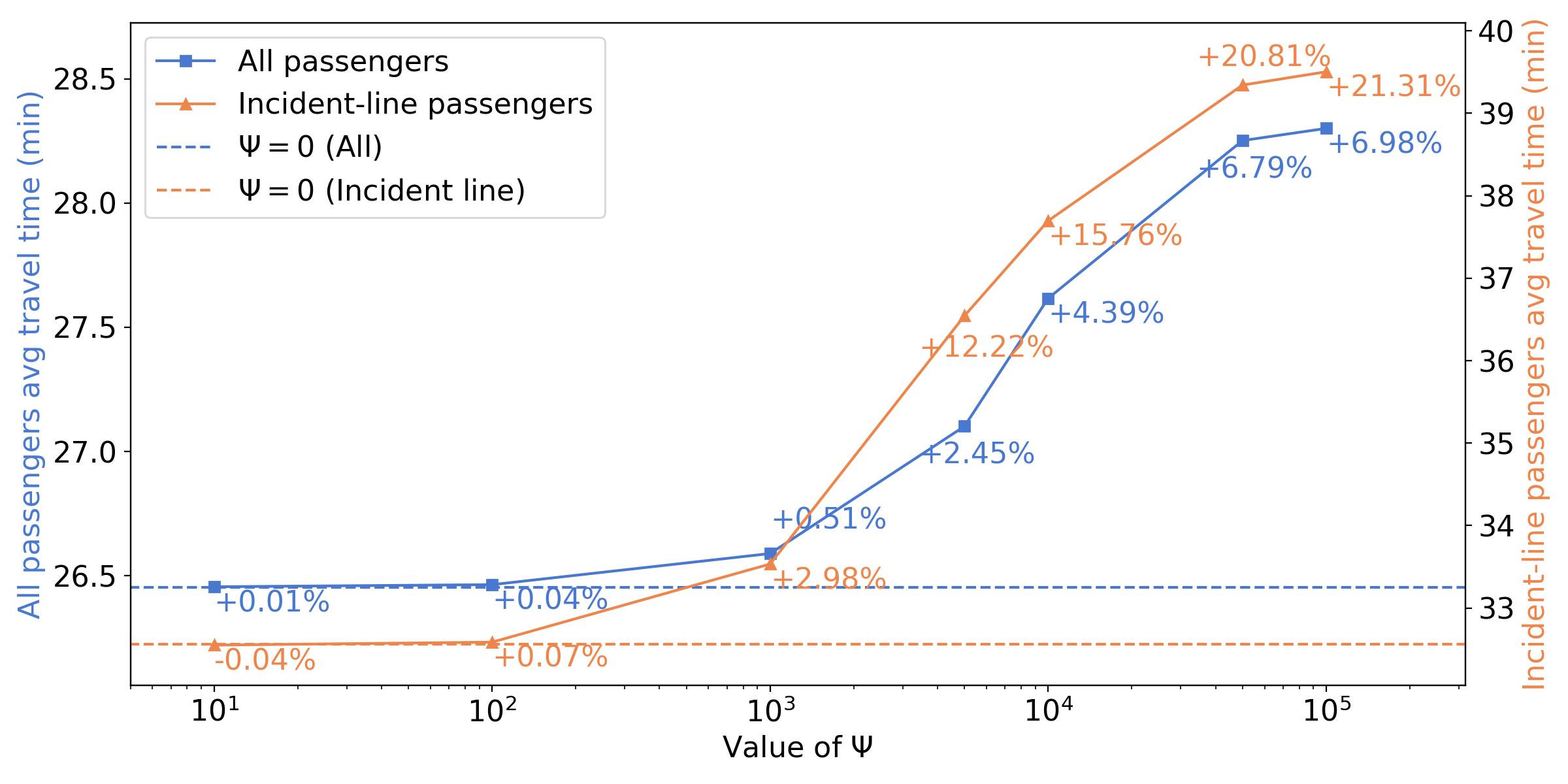}\label{fig_pax_prefer1}}
\hfil
\subfloat[Total utility and number of passengers being recommended with preferred path]{\includegraphics[width=0.5\textwidth]{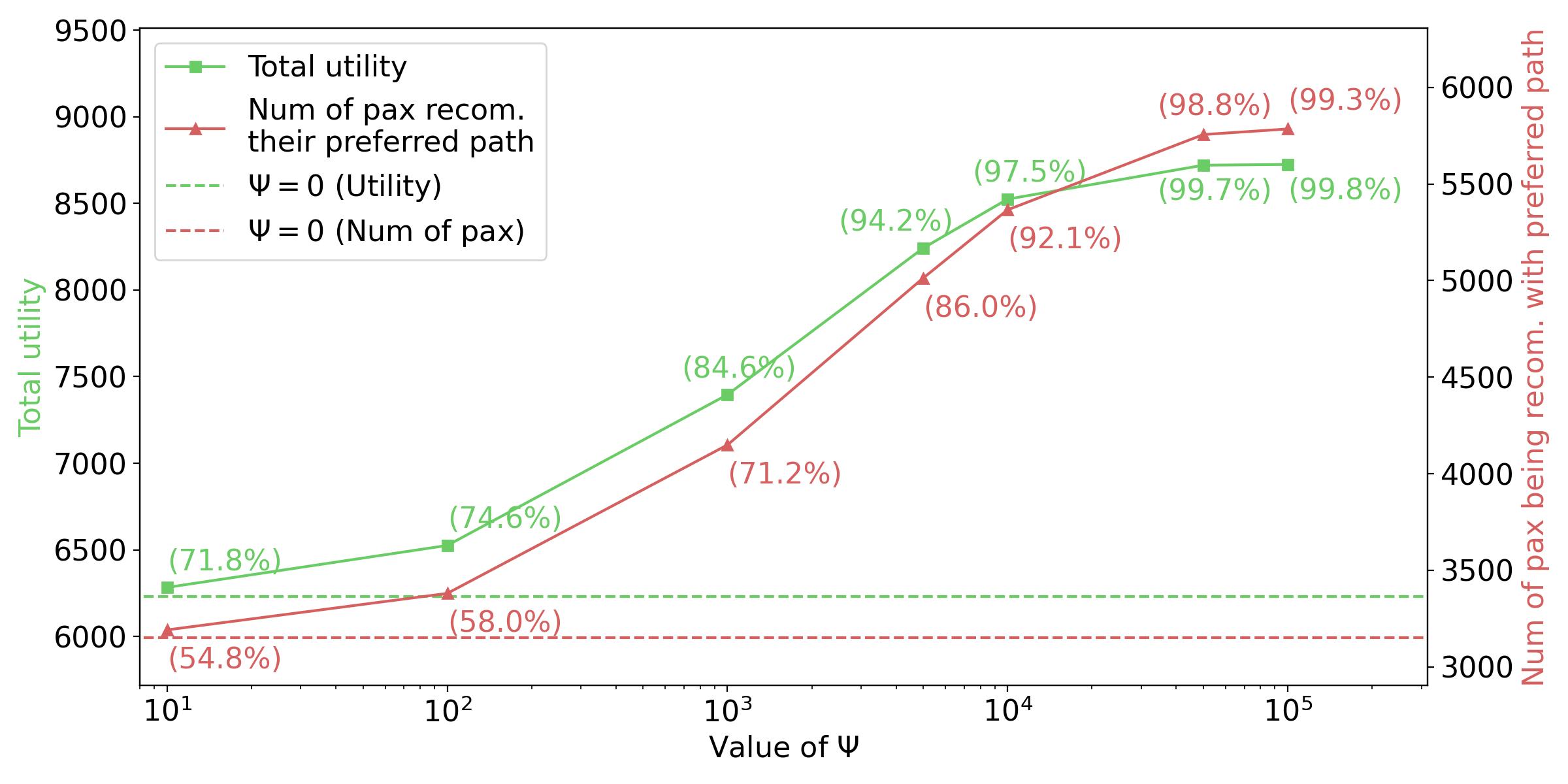}\label{fig_pax_prefer2}}
\caption{Impact of different values of $\Psi$ on results. The percentage change in Figure (a) is compared with the scenario of $\Psi = 0$. The percentage in parentheses in Figure (b) represents the relative ratio of total utility and proportion of passengers recommended with their preferred path, respectively.}
\label{fig_compare_psi}
\end{figure}

Figure \ref{fig_compare_psi} shows the results for different values of $\Psi$. All models use $\epsilon=0, \Gamma=1$. The x-axis is plotted in a log scale. In Figure \ref{fig_pax_prefer1}, the average travel time for all passengers and incident-line passengers increases with the increase of $\Psi$, which is as expected because the larger value of $\Psi$ means that the recommendation generation focuses more on satisfying passenger's inherent preferences rather than minimizing the system travel time. Similarly, in Figure \ref{fig_pax_prefer2},  as expected, both $TU(\boldsymbol{x})$ and $NP(\boldsymbol{x})$ increase with the increase in  $\Psi$. When $\Psi = 10^5$, the average travel time of the incident line passengers increased by 21.3\%, which is close to the status quo scenario. This is because we generate passengers' prior utilities based on the status quo choices. Figure \ref{fig_pax_prefer2} shows that nearly all passengers in $\Pcal$ are recommended with their preferred path when $\Psi = 10^5$.

Figure \ref{fig_compare_psi} illustrates the trade-off between respecting passengers' preferences and reducing system congestion. When the value of  $\Psi$ is relatively small (e.g., less than $10^3$), increasing $\Psi$ can effectively increase the total utility and number of passengers recommended with their preferred path. Meanwhile, the system travel time only slightly increases. But when $\Psi$ is large (e.g., greater than  $10^4$), increasing  $\Psi$ significantly increases the system travel time, but the impact on increasing the passenger's utility is limited. The reason may be that, in the system, there are some passengers whose preferred paths are not at the capacity bottlenecks. Hence, when $\Psi$ is small, the optimal solution recommends those passengers use their preferred paths without significantly impacting the system travel time. When $\Psi$ is large, passengers are recommended to use their preferred paths even if these paths are highly congested, causing a significant increase in the system travel time. The results imply that a reasonable value of $\Psi$ should be relatively small. With small $\Psi$, most of the passengers (e.g., more than 70\%) are recommended to use their preferred paths without significantly reducing the system efficiency.

\subsection{Impact of different behavior uncertainty scenarios}\label{sec_res_unc}
The model's performance may vary according to the level of passengers' uncertainties. In this section, we test different scenarios on passenger's conditional path choice probabilities. We modify the impact of the recommendation of path $r'$ on the utility of path $r$ as
\begin{equation}
I_{p,r'}^r=\left\{
\begin{aligned}
& \text{Drawn from }\mathcal{U}[I^{\text{Min}},I^{\text{Max}}]& \quad\text{     if $r = r'$} \\
& 0 & \quad\text{     otherwise},
\end{aligned}
\right. \quad \forall \; p\in \mathcal{P}, r,r'\in \mathcal{R}_p \label{eq_gen_pre_scenario}.
\end{equation}
Besides our previous testing of $I^{\text{Min}}=0, I^{\text{Max}}=5$, we added following additional scenarios:
\begin{itemize}
    \item Smaller recommendation impact: passengers are less affected by the recommendation of the system, we consider $I^{\text{Min}}=0, I^{\text{Max}}=3$ and $I^{\text{Min}}=0, I^{\text{Max}}=1$
    \item Random recommendation impact: passengers may be positively or negatively affected by the recommendation. We consider $I^{\text{Min}}=-5, I^{\text{Max}}=5$, $I^{\text{Min}}=-3, I^{\text{Max}}=3$, and $I^{\text{Min}}=-1, I^{\text{Max}}=1$.
\end{itemize}

Results are shown in Table \ref{tab_result_unc}. The IPR model with behavior uncertainty (BU) uses $\epsilon=0, \Gamma=1$. Each number is the average of 10 replications for travel time evaluations. The results highlight that incorporating behavior uncertainty in the IPR model leads to significant improvements in travel time performance, especially under scenarios with larger behavior uncertainties (e.g., $I^{\text{Min}}=-5, I^{\text{Max}}=5$). The average travel time for incident line passengers decreases from 36.67 mins to 33.58 mins (-8.4\%). The model is robust in handling both positive and negative recommendation influences.


\begin{table}[htb]
\centering
\caption{Impact of different behavior uncertainty scenarios}\label{tab_result_unc}
{
{\renewcommand{\arraystretch}{1} 
\begin{tabular}{@{}cccccc@{}}
\toprule
\multirow{2}{*}{{Scenario [$I^{\text{Min}}, I^{\text{Max}}$]}} & \multirow{2}{*}{{Model}} & \multicolumn{2}{c}{{ATT (all passengers)}} & \multicolumn{2}{c}{{ATT (incident line passengers)}} \\ \cmidrule(l){3-6} 
& & Mean (min)                & Std. (min)          & Mean (min)              & Std. (min)                       \\ \midrule
\multirow{3}{*}{{[}0, 5{]}}           & Capacity-based                     &   27.609 &0.033                      &        33.848 & 0.165                                                                \\
                                      & IPR (w.o. BU)                      &  26.706 &0.026                                                        &   32.852 &0.122                                                                       \\
                                      & \cellcolor{gray!25}IPR (w. BU)                        &  \cellcolor{gray!25}26.454 & \cellcolor{gray!25}0.019                                                              &        \cellcolor{gray!25}32.526 & \cellcolor{gray!25}0.204                                                                \\ \midrule
\multirow{3}{*}{{[}0, 3{]}}           & Capacity-based                     &   27.777 &0.034                      &        34.927 &0.181                                                                  \\
                                      & IPR (w.o. BU)                      &  27.177 &0.031                                                        &   34.475 &0.131                                                                       \\
                                      & \cellcolor{gray!25}IPR (w. BU)                        &  \cellcolor{gray!25}26.771 & \cellcolor{gray!25}0.020                                                              &        \cellcolor{gray!25}33.843 & \cellcolor{gray!25}0.211                                                                \\ \midrule
\multirow{3}{*}{{[}0, 1{]}}           & Capacity-based                     &   28.202 &0.037                      &        37.037 &0.201                                                                  \\
                                      & IPR (w.o. BU)                      &  28.049 &0.034                                                        &   36.978 &0.136                                                                       \\
                                      & \cellcolor{gray!25}IPR (w. BU)                        &  \cellcolor{gray!25}27.535 & \cellcolor{gray!25}0.022                                                              &        \cellcolor{gray!25}36.202 & \cellcolor{gray!25} 0.215                                                                \\ \midrule
\multirow{3}{*}{{[}-5, 5{]}}           & Capacity-based                     &   28.252 &0.044                      &        36.667 &0.221                                                                  \\
                                      & IPR (w.o. BU)                      &  28.001 &0.038                                                        &   36.176 &0.142                                                                       \\
                                      & \cellcolor{gray!25}IPR (w. BU)                        &  \cellcolor{gray!25}26.714 & \cellcolor{gray!25}0.026                                                              &        \cellcolor{gray!25}33.578 & \cellcolor{gray!25} 0.210                                                                \\ \midrule
\multirow{3}{*}{{[}-3, 3{]}}           & Capacity-based                     &   28.324 &0.045                      &        37.236 &0.217                                                                  \\
                                      & IPR (w.o. BU)                      &  28.123 &0.041                                                        &   36.807 &0.139                                                                       \\
                                      & \cellcolor{gray!25}IPR (w. BU)                        &  \cellcolor{gray!25}27.018 & \cellcolor{gray!25}0.028                                                              &        \cellcolor{gray!25}34.931 & \cellcolor{gray!25}0.214                                                                \\ \midrule
\multirow{3}{*}{{[}-1, 1{]}}           & Capacity-based                     &   28.418 &0.047                     &        37.864 &0.214                                                                  \\
                                      & IPR (w.o. BU)                      &  28.394 &0.042                                                        &   37.761 &0.137                                                                       \\
                                      & \cellcolor{gray!25}IPR (w. BU)                        &  \cellcolor{gray!25}27.864 & \cellcolor{gray!25}0.030                                                              &        \cellcolor{gray!25}36.857& \cellcolor{gray!25}0.222                                                                \\ \bottomrule

\multicolumn{6}{l}{\begin{tabular}[c]{@{}l@{}} \small Numbers in parentheses represent standard deviations; Best results are highlighted in gray \\ \end{tabular}}

\end{tabular}
}
}
\end{table}

\subsection{Impact of incident durations}\label{sec_res_incident_dur}
In this section, we test the model's performance under different incident durations. Note that the actual incident duration is 59 minutes. We assume the demand patterns are the same for all incident scenarios. Results are shown in Table \ref{tab_result_dur}. The IPR model with behavior uncertainty (BU) uses $\epsilon=0, \Gamma=1$.

We find that, across all durations analyzed (30, 59, 90, and 120 minutes), the IPR model considering behavior uncertainty consistently yields the lowest average travel time (ATT), highlighting its effectiveness in mitigating the adverse effects of incidents compared to the other approaches. Moreover, the effect is more prominent in cases with longer incidents. 

\begin{table}[htb]
\centering
\caption{Impact of different incident duration}\label{tab_result_dur}
{
{\renewcommand{\arraystretch}{1} 
\begin{tabular}{@{}cccccc@{}}
\toprule
\multirow{2}{*}{{Incident duration}} & \multirow{2}{*}{{Model}} & \multicolumn{2}{c}{{ATT (all passengers)}} & \multicolumn{2}{c}{{ATT (incident line passengers)}} \\ \cmidrule(l){3-6} 
& & Mean (min)                & Std. (min)                  & Mean (min)              & Std. (min)                       \\ \midrule
\multirow{3}{*}{30 min}           & Capacity-based                     &   26.038 &0.031                      &        30.538 &0.155                                                                  \\
                                      & IPR (w.o. BU)                      &  24.976 &0.024                                                        &   29.753 &0.115                                                                       \\
                                      & \cellcolor{gray!25}IPR (w. BU)                        &  \cellcolor{gray!25}24.802 &\cellcolor{gray!25}0.017                                                              &        \cellcolor{gray!25} 29.539 &\cellcolor{gray!25}0.178                                                                \\ \midrule
\multirow{3}{*}{59 min (actual)}           & Capacity-based                     &   27.609 &0.033 
&        33.848 &0.165                                                                  \\
                                      & IPR (w.o. BU)                      &  26.706 &0.026                                                        &   32.852 &0.122                                                                       \\
                                      & \cellcolor{gray!25}IPR (w. BU)                        &  \cellcolor{gray!25}26.454 & \cellcolor{gray!25} 0.019                                                              &        \cellcolor{gray!25}32.526 & \cellcolor{gray!25} 0.204                                                                \\ \midrule

\multirow{3}{*}{90 min}           & Capacity-based                     &   29.643 &0.037                      &        38.770 &0.170                                                                  \\
                                      & IPR (w.o. BU)                      &  28.721 &0.028                                                        &   36.953 &0.124                                                                       \\
                                      & \cellcolor{gray!25}IPR (w. BU)                        &  \cellcolor{gray!25}28.431 & \cellcolor{gray!25} 0.021                                                              &        \cellcolor{gray!25}36.204 & \cellcolor{gray!25} 0.211                                                                \\ \midrule
\multirow{3}{*}{120 min}           & Capacity-based                     &   32.714 &0.040                      &        44.476 &0.182                                                                  \\
                                      & IPR (w.o. BU)                      &  30.746 &0.031                                                        &   41.473 &0.131                                                                       \\
                                      & \cellcolor{gray!25}IPR (w. BU)                        &  \cellcolor{gray!25}30.325 & \cellcolor{gray!25} 0.025                                                              &        \cellcolor{gray!25}40.588 & \cellcolor{gray!25} 0.228                                                                \\ \bottomrule

\multicolumn{6}{l}{\begin{tabular}[c]{@{}l@{}} \small Numbers in parentheses represent standard deviations; Best results are highlighted in gray \\ 
\end{tabular}}

\end{tabular}
}
}
\end{table}

\section{Synthetic case study}\label{sec_syn_case_study}
In the synthetic case study, we focus on testing the model's generalizability to different network scales and incident situations. 
\subsection{Case study design}
\begin{figure}[htb]
\centering
\includegraphics[width = 1.0\linewidth]{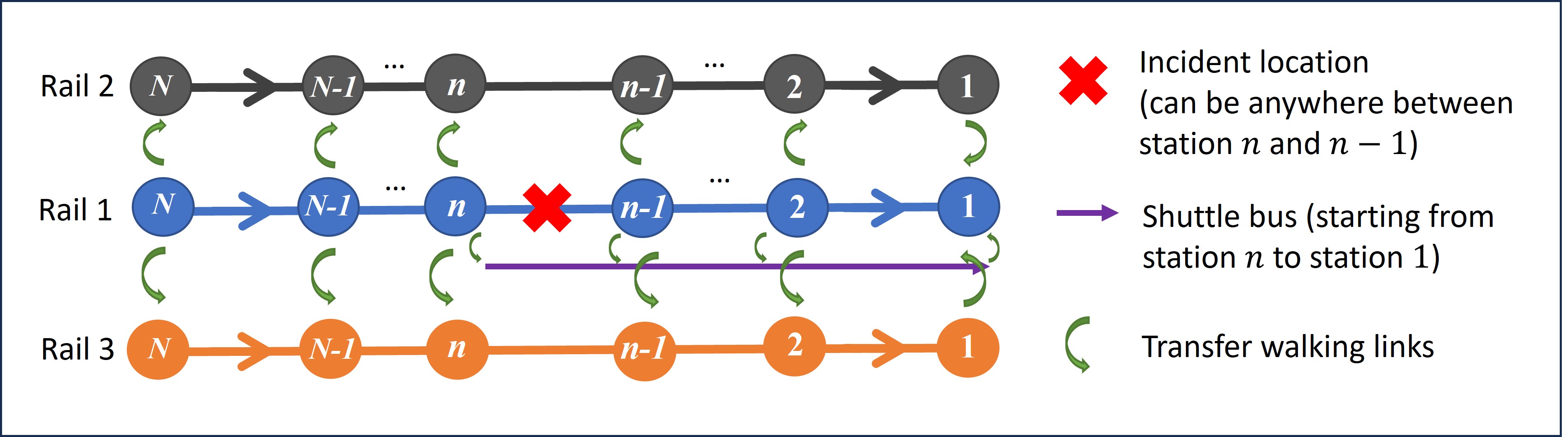}
\caption{Synthetic case study design}
\label{fig_syn_incn}
\end{figure}

Consider a general 3-line synthetic network (Figure \ref{fig_syn_incn}). Each line has $N$ stations, numbering from 1 to $N$. For any stations $n\in[N]$ in rail line 1 (where $[N]=\{1,2,...,N\}$), people can transfer to rail lines 2 or 3 through walking, providing them alternatives during disruptions. For simplicity, we assume station 1 in rail line 1 is the only destination people are willing to go to (e.g., the downtown area for the morning commute). And only rail line 1 has demand. For all OD pairs $(n,1),\; \forall n \in[N]\setminus\{1\}$, the demand $d^{n,1}_t$ is assumed to be the same (i.e., uniform demand). We assume there is an incident happening between stations $\lceil N/2 \rceil+1$ and $\lceil N/2 \rceil$ at rail line 1, lasting for 1 hour (the location of the incident will vary in Section \ref{sec_res_inc_loc}). And there is a shuttle bus serving between stations $n$ and $1$ during the disruptions. All individuals have four options: 1) waiting at rail line 1 until it recovers, 2) transferring to rail line 2, 3) transferring to rail line 3, 4) transferring to the shuttle bus. 

The uniform demand is set as 25 passengers per OD pair per hour. The vehicle capacities for rail lines 1, 2, and 3 and shuttle buses are 500, 300, 300, and 40, respectively. The headway under normal operation times for rail lines 1, 2, and 3 are 10 minutes, 12 minutes, and 13 minutes, respectively. The headway for the shuttle bus during the disruption is 8 minutes. The transfer walking time between rail lines is 10 minutes, and from the rail line to the shuttle bus is 3 minutes. The vehicle travel time between any two stations $n$ and $n-1$ for rail lines 1, 2, and 3 are 5 minutes, 7 minutes, and 8 minutes respectively, and for the shuttle bus is 10 minutes. These settings generally tell us that rail line 1 has the most carrying capacity, followed by rail lines 2, 3, and shuttle buses. The path recommendation hyperparameters are set as $\epsilon=0$, $\Gamma=1$, $\Psi=0$. Other assumptions and setups are the same as the actual case study. In the following testing, we vary $N=2$ (i.e., two-station network) to $N=20$ to test the model performance extensively. 
\subsection{Computational time}
The computational results are summarized in Table \ref{tab_time_comp}. The relative performance of Benders Decomposition (BD) and the Gurobi solver vary across network sizes. For small-scale scenarios (e.g., $N=2$ to $N=10$), Gurobi consistently outperforms BD with lower CPU times. However, as the network size increases beyond $N=12$, the efficiency of BD becomes more apparent. While Gurobi's computational time grows rapidly with increasing problem size—reaching over 450 seconds at $N=20$, BD demonstrates more moderate growth, maintaining a substantial advantage in large-scale instances. This indicates that BD scales more efficiently and is better suited for solving large-scale IPR problems.

\begin{table}[htbp]
\centering
\caption{Computational time comparison on synthetic networks}\label{tab_time_comp}
{
{\renewcommand{\arraystretch}{1} 
\begin{tabular}{@{}ccc|ccc@{}}
\toprule
Scenario & Solver                 & CPU time (sec)   &Scenario& Solver & CPU time (sec)      \\ \midrule
\multirow{2}{*}{$N=2$} & BD & 3.18      &\multirow{2}{*}{$N=12$} & BD & 12.61                     \\
                  & Gurobi                     & 0.51 &    & Gurobi                  & 11.55   \\ \midrule
\multirow{2}{*}{$N=4$} & BD & 3.66      &\multirow{2}{*}{$N=14$} & BD & 16.60                     \\
                  & Gurobi                     & 0.70 &    & Gurobi                    &35.62  \\    \midrule   
\multirow{2}{*}{$N=6$} & BD & 4.12      &\multirow{2}{*}{$N=16$} & BD & 23.45                     \\
                  & Gurobi                     & 0.95 &    & Gurobi                    & 80.79 \\ \midrule
\multirow{2}{*}{$N=8$} & BD & 5.97      &\multirow{2}{*}{$N=18$} & BD & 76.52                     \\ 
                  & Gurobi                     & 1.14 &    & Gurobi                    &  231.83\\ \midrule
\multirow{2}{*}{$N=10$} & BD & 8.18      &\multirow{2}{*}{$N=20$} & BD & 134.42                     \\
                  & Gurobi                     & 2.01 &    & Gurobi                    &  451.32\\\bottomrule
\end{tabular}
}
}
\end{table}

\subsection{Performance}
The results in Table~\ref{tab_net_res} demonstrate that the proposed IPR model consistently outperforms the capacity-based method in terms of average travel time (ATT) across all tested network sizes. Specifically, the IPR model achieves noticeable reductions in ATT, ranging from 13.3\% at $N=4$ to as much as 15.0\% at $N=2$. However, as the network scale increases (e.g., from $N=12$ to $N=20$), the relative improvement gradually diminishes, with the reduction dropping to just 1.8\% at $N=20$. This trend can be attributed to increasing demand and tighter remaining capacity in larger-scale scenarios. Under such conditions, the flexibility of the IPR model becomes more limited, as fewer excess resources are available to reassign passengers optimally. Hence, while the IPR model maintains better performance, its advantage over the capacity-based approach narrows as the system becomes more saturated.

\begin{table}[htb]
\centering
\caption{Average travel time (ATT) comparison for different network scales}\label{tab_net_res}
\resizebox{\textwidth}{!}
{
{\renewcommand{\arraystretch}{1} 
\begin{tabular}{cccc|cccc}
\toprule
\multirow{2}{*}{Scenario} & \multirow{2}{*}{Models} & \multicolumn{2}{c|}{ATT (all passengers)} & \multirow{2}{*}{Scenario} & \multirow{2}{*}{Models} & \multicolumn{2}{c}{ATT (all passengers)} \\ 
                      & & Mean (min)               & Std.  (min)        &&        & Mean (min)                       & Std.  (min)                      \\ \midrule
 \multirow{2}{*}{$N=2$}    &   Capacity-based     & 12.67                      & 0.101                      &  \multirow{2}{*}{$N=12$} &   Capacity-based & 65.81                           & 0.266 \\
  & IPR        & 10.77 (-15.0\%)                   & 0.064               &      & IPR & 60.02 (-8.8\%)                        & 0.145                           \\ \midrule
 \multirow{2}{*}{$N=4$}    &   Capacity-based     & 20.91                      & 0.142                      &  \multirow{2}{*}{$N=14$} &   Capacity-based & 76.76                           & 0.301\\
  & IPR        & 18.14 (-13.3\%)                   & 0.077               &      & IPR & 72.26 (-5.9\%)                        & 0.169                           \\ \midrule
   \multirow{2}{*}{$N=6$}    &   Capacity-based     & 29.64                      & 0.174                      &  \multirow{2}{*}{$N=16$} &   Capacity-based & 86.87                           & 0.331 \\
  & IPR        & 26.95 (-9.1\%)                   & 0.083               &      & IPR & 83.48 (-3.9\%)                        & 0.203                          \\ \midrule
   \multirow{2}{*}{$N=8$}    &   Capacity-based     & 39.79                      & 0.205                     &  \multirow{2}{*}{$N=18$} &   Capacity-based & 93.28                           & 0.379 \\
  & IPR        & 35.59 (-10.6\%)                   & 0.104               &      & IPR & 89.72 (-3.8\%)                        & 0.222                           \\ \midrule
   \multirow{2}{*}{$N=10$}    &   Capacity-based     & 51.40                      & 0.238                      &  \multirow{2}{*}{$N=20$} &   Capacity-based & 96.93                           & 0.403 \\
  & IPR        & 46.39 (-9.7\%)                   & 0.121               &      & IPR & 95.15 (-1.8\%)                        & 0.276                          \\ 
\bottomrule
\multicolumn{8}{l}{\begin{tabular}[c]{@{}l@{}} \small Numbers in parentheses represent percentage travel time reduction compared to the capacity-based method \end{tabular}}

\end{tabular}
}
}
\end{table}

\subsection{Impact of incident location}\label{sec_res_inc_loc}
To evaluate the impact of the incident location, we choose scenarios with $N \in \{8,12,16,20\}$. For each scenario, we consider three different incident locations: Upstream (incident between stations $\lceil 3N/4 \rceil+1$ and $\lceil 3N/4 \rceil$), Middle (incident between stations $\lceil N/2 \rceil+1$ and $\lceil N/2 \rceil$, i.e., same as the above), and Downstream (incident between stations $\lceil N/4 \rceil+1$ and $\lceil N/4 \rceil$). 

The results in Table~\ref{tab_inc_loca} indicate that incident location has no significant impact on the average travel time under either the capacity-based method or the proposed IPR model. Across all scenarios ($N=8, 12, 16, 20$) and for all three incident locations (Upstream, Middle, Downstream), the average travel time values remain relatively consistent, with only minor variations. This outcome is expected, as the model assumes that a service disruption completely blocks the entire line regardless of where the incident occurs (an assumption that reflects common real-world operational protocols). The only difference introduced by varying incident locations lies in how the shuttle bus service is designed. For upstream incidents, the shuttle bus coverage tends to be longer, which can slightly reduce the total system travel time. However, due to the limited capacity of shuttle buses, their influence on the overall passenger travel time remains marginal. 

\begin{table}[htp]
\centering
\caption{Average travel time (ATT) comparison under different incident locations}\label{tab_inc_loca}
\resizebox{\textwidth}{!}
{
{\renewcommand{\arraystretch}{1} 
\begin{tabular}{cccccccccc}
\toprule
\multirow{2}{*}{Scenario} & \multirow{2}{*}{Location} & \multirow{2}{*}{Models} & \multicolumn{2}{c}{ATT (all passengers)} & \multirow{2}{*}{Scenario} & \multirow{2}{*}{Location} & \multirow{2}{*}{Models} & \multicolumn{2}{c}{ATT (all passengers)} \\ 
                     & & & Mean (min)               & Std.  (min)        & &  &      & Mean (min)                       & Std.  (min)                      \\ \midrule
 \multirow{6}{*}{$N=8$}    &  \multirow{2}{*}{Upstream} &  Capacity-based     & 39.55                      & 0.207                      &  \multirow{6}{*}{$N=16$}  &  \multirow{2}{*}{Upstream} & Capacity-based & 86.75                           & 0.328 \\
  & & IPR        & 35.43 (-10.4\%)                   & 0.110               &     & & IPR & 83.29 (-4.0\%)                        & 0.201                           \\ \cmidrule{2-5}  \cmidrule{7-10} 
  &  \multirow{2}{*}{Middle} &  Capacity-based     & 39.79                      & 0.205                      &   &  \multirow{2}{*}{Middle} & Capacity-based & 86.87                           & 0.331 \\
  & & IPR        & 35.59 (-10.6\%)                   & 0.104               &     & & IPR & 83.48 (-3.9\%)                        & 0.203                           \\ \cmidrule{2-5}  \cmidrule{7-10} 
  &  \multirow{2}{*}{Downstream} &  Capacity-based     & 39.91                      & 0.203                      &   &  \multirow{2}{*}{Downstream} & Capacity-based & 86.91                           & 0.332 \\
  & & IPR        & 35.73 (-10.5\%)                   & 0.106               &     & & IPR & 83.62 (-3.8\%)                        & 0.204                           \\ \midrule
 \multirow{6}{*}{$N=12$}    &  \multirow{2}{*}{Upstream} &  Capacity-based     & 65.67                      & 0.263                      &  \multirow{6}{*}{$N=20$}  &  \multirow{2}{*}{Upstream} & Capacity-based & 96.74                           & 0.398 \\
  & & IPR        & 59.84 (-8.9\%)                   & 0.144               &     & & IPR & 95.11 (-1.7\%)                        & 0.275                           \\ \cmidrule{2-5}  \cmidrule{7-10} 
  &  \multirow{2}{*}{Middle} &  Capacity-based     & 65.81                      & 0.266                      &   &  \multirow{2}{*}{Middle} & Capacity-based & 96.93                           & 0.403 \\
  & & IPR        & 60.02 (-8.8\%)                   & 0.145               &     & & IPR & 95.15 (-1.8\%)                        & 0.276                           \\ \cmidrule{2-5}  \cmidrule{7-10} 
  &  \multirow{2}{*}{Downstream} &  Capacity-based     & 65.93                      & 0.266                      &   &  \multirow{2}{*}{Downstream} & Capacity-based & 97.08                           & 0.406 \\
  & & IPR        & 60.52 (-8.2\%)                   & 0.144               &     & & IPR & 95.22 (-1.9\%)                        & 0.273                           \\ 
\bottomrule
\multicolumn{10}{l}{\begin{tabular}[c]{@{}l@{}} \small Numbers in parentheses represent percentage travel time reduction compared to the capacity-based method \end{tabular}}

\end{tabular}
}
}
\end{table}

\section{Conclusion and discussion}\label{sec_conclusion}
This study proposes an individual path recommendation model during PT service disruptions with the objective of minimizing total system travel time and respecting passengers' path choice preferences. Passengers' behavior uncertainty in path choices given recommendations is also considered in the formulation. The original IPR formulation yields a stochastic optimization with decision-dependent distributions. We propose a single-point approximation method to eliminate the expectation operator by introducing two new concepts: $\epsilon$-feasibility and $\Gamma$-concentration. The approximation yields a tractable single-stage mixed integer linear formulation, which can be solved efficiently with Benders decomposition. The approximation gap is proved to be bounded from the above. Additional theoretical analysis shows that $\epsilon$-feasibility and $\Gamma$-concentration are strongly connected to expectation and chance constraints in a typical stochastic optimization formulation, respectively. 

The proposed approach is demonstrated in a real-world case study using data from an urban rail disruption in the CTA system, and a synthetic case study with varied network sizes and incident locations. In the real-world case study, results show that the proposed IPR model reduces the average travel times in the system by 6.6\% compared to the status quo and by 4.2\% compared to a capacity-based benchmark model. In the synthetic case study, the proposed model shows 15.0\% to 1.8\% lower system travel time compared to the capacity-based method, depending on the network sizes and demand situations.

Future studies can be pursued in the following directions. First, it is possible to extend the current framework with more complex recommendation compositions. The challenges in implementing the more general framework stem from the quantification of the posterior path choice probabilities. Future studies may conduct corresponding surveys to calibrate passengers' responses to the recommendations. Second, future studies may consider different sources of uncertainty (including incident duration, in-vehicle time, demand, etc.) for a more realistic modeling framework. Third, future studies could extend beyond individual path recommendations to a co-design framework that integrates personalized routing with service adjustments, such as bus rerouting and shuttle deployment. This approach will optimize both passenger travel experience and operational efficiency by balancing travel time, accessibility, and cost considerations.


\section{Acknowledgement}
The authors would like to thank the Chicago Transit Authority (CTA) for their support and data availability for this research.

\section*{Appendices}
\begin{APPENDICES}

\section{Model extensions}\label{sec_model_exten}
In this section, we discuss several extensions of the model to accommodate more realistic/general scenarios. 
\subsection{Generalization of recommendations}\label{sec_rec_com}
In this study, we assume the information given to passengers is a recommended path. In reality, the recommendation system may provide a bundle of recommended paths with information like estimated in-vehicle time, waiting time, travel cost, etc. The proposed framework can be extended to handle different recommendation typologies. Figure \ref{fig_bh_un_gen} shows an example where the recommendation system will provide a composition of path and travel time information, where each composition can include different paths, different estimated waiting/in-vehicle times, etc. Then, we can change $x_{p,r}$ to $x_{p,c}$, where $x_{p,c}$ indicates whether we will present composition $c$ to passenger $p$. Similarly, each $c$ is associated with a conditional probability $\pi_{p,c}^r$ as shown in Figure \ref{fig_bh_un_gen} (the probability for passenger $p$ to choose path $r$ given that he/she is recommended composition $c$). 
\begin{figure}[htb]
\centering
\includegraphics[width = 1 \linewidth]{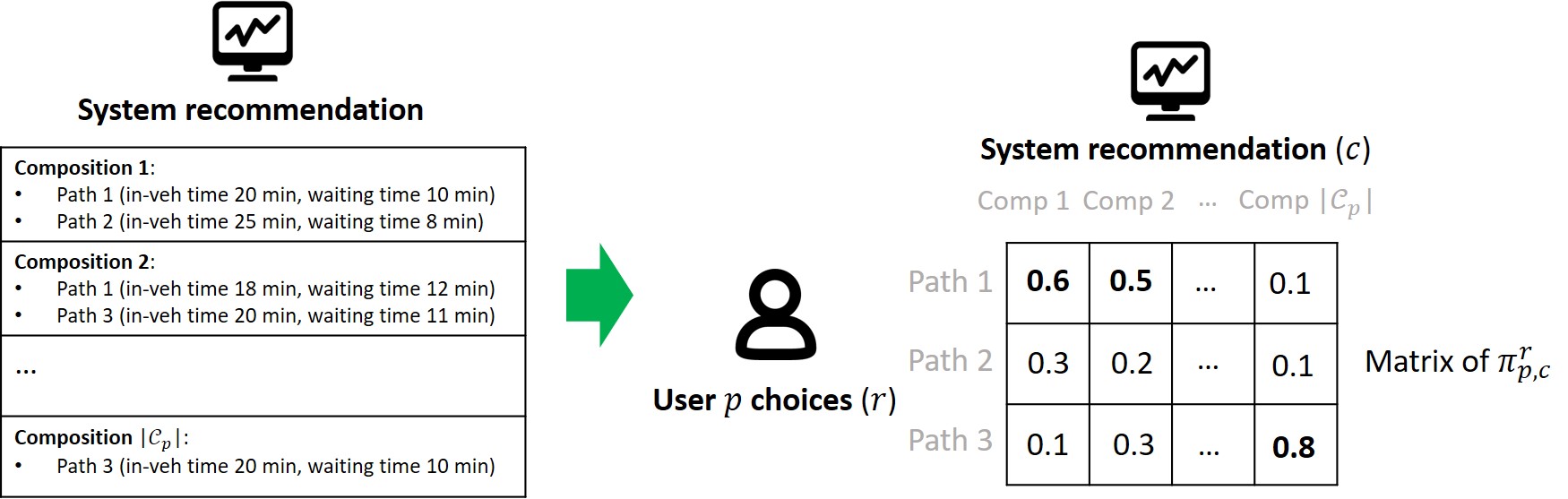}
\caption{Illustration of the generalized recommendation typology. $\Ccal_p$ is the predetermined recommendation composition sets for passenger $p$}
\label{fig_bh_un_gen}
\end{figure}

In this way, we only need to calibrate $\pi_{p,c}^r$ and predetermine the composition set $\mathcal{C}_{p}$ for each passenger $p$. The overall framework proposed above can be easily adapted to the new recommendation typology by replacing $x_{p,r}$ and $\pi_{p,r'}^r$ with $x_{p,c}$ and $\pi_{p,c}^r$, respectively. 

\subsection{Feedback and rolling-horizon}\label{append_roll}
As mentioned in Section \ref{sec_ana_des}, the whole path recommendation problem could be solved in a rolling-horizon manner. At each time interval $t \geq 1$, we update the demand, supply, and system state information, and solve the proposed framework above to get a recommendation strategy $\boldsymbol{x}$. But we only implement the $x_{p,r}$ for $p\in\Pcal^{u,v}_t, \; \forall (u,v)$ (i.e., passengers departing at current time $t$). Specifically, the inputs can be updated as follows. The new supply information can be updated by querying the operators' database. The onboard passengers $\Omega_1$ can be updated based on new tap-in transactions and the simulation from the last time step. The predicted demand without recommendation $f_{t}^{u,v,r}$ can be updated by running the prediction model again with new smart card data transactions. We also need to construct the sets of passengers need recommendations in future time steps: ${\mathcal{P}}_{t+1}^{u,v}, ..., {\mathcal{P}}_{T^D}^{u,v}$. Note that some passengers may already make recommendation requests even if they have not arrived at the system (i.e., part of the passengers in ${\mathcal{P}}_{t+1}^{u,v}, ..., {\mathcal{P}}_{T^D}^{u,v}$ are already known). Then, we can add additional ``predicted individuals'' to match the predicted future demand of each OD pair $(u,v)$ and time $t+h$ in the future with the associated preference matrix. This can be done by sampling passengers from the existing transit application database (usually these applications have passengers' historical travel patterns). We assume all recommendations for passengers who arrived before $t$ are fixed (even if they are not the optimal ones given the newly available information). This can be done by simply adding a new constraint $x_{p,r}=\hat{x}_{p,r}, \;\forall p \in \cup_{t'=1,...,t-1}\mathcal{P}_{t'}^{u,v}$, where $\hat{x}_{p,r}$ is the previous decision. But we can still modify the decisions for passengers in $\mathcal{P}_t^{u,v}, {\mathcal{P}}_{t+1}^{u,v}, ..., {\mathcal{P}}_{T^D}^{u,v}$ as they have not arrived at the system by time $t$. The rolling horizon may also update the associated preference matrix information through passenger feedback when they have not arrived at the transit system. For example, after providing a recommendation, we can ask the passenger to respond whether he/she will actually use it or not. This feedback can be used to update the associated preference matrix.

\section{Proof of Proposition \ref{prop_gamma_linear}}\label{proof_prop_gamma_conc}

From the triangular inequality, we have:
\begin{align}
     \underbrace{\left|Q_t^{u,v,r} - {q}_{t}^{u,v,r}\right|}_{\text{LHS}} & \leq  \left|Q_t^{u,v,r} - {\mu}_{t}^{u,v,r}(\boldsymbol{x})\right| +  \left| {\mu}_{t}^{u,v,r}(\boldsymbol{x}) - {q}_{t}^{u,v,r} \right| \notag \\
     &\leq \underbrace{\left|Q_t^{u,v,r} - {\mu}_{t}^{u,v,r}(\boldsymbol{x})\right| +\epsilon_{t}^{u,v,r} }_{\text{RHS}} .
\end{align}
As $\text{LHS} \leq \text{RHS}$, the probability measure satisfies (for all $a > \epsilon_{t}^{u,v,r}$):
\begin{align}
    \mathbb{P}[\text{LHS} \geq a] \leq \mathbb{P}[\text{RHS} \geq a] \label{eq_lhs_rhs}.
\end{align}
Notice that
\begin{align}
    \mathbb{P}[\text{RHS} \geq a] = \mathbb{P}\left[\left| Q_t^{u,v,r} - {\mu}_{t}^{u,v,r}(\boldsymbol{x})\right| \geq a - \epsilon_{t}^{u,v,r} \right] \leq \frac{({\sigma}_{t}^{u,v,r}(\boldsymbol{x}))^2}{(a - \epsilon_{t}^{u,v,r})^2} .
    \label{eq_rhs}
\end{align}
Equation (\ref{eq_rhs}) is based on Chebyshev's inequality. Therefore,
\begin{align}
\mathbb{P}[\text{LHS} \geq a] = \mathbb{P}[\left|Q_t^{u,v,r} - {q}_{t}^{u,v,r}\right|\geq a] \leq \frac{({\sigma}_{t}^{u,v,r}(\boldsymbol{x}))^2}{(a - \epsilon_{t}^{u,v,r})^2} \label{eq_proof2}.
\end{align}
Comparing (\ref{eq_proof2}) and (\ref{eq_gamma_concen}), we know that to satisfy (\ref{eq_gamma_concen}), we only need ${\sigma}_{t}^{u,v,r}(\boldsymbol{x}) \leq \Gamma_{t}^{u,v,r}$, which completes the proof.

\section{Proof of Lemma \ref{lemma_berge}}\label{app_proof_lemma}

$STT({{\boldsymbol{q}}})$ is obtained by solving linear programming. ${{\boldsymbol{q}}}$ is a parameter in the constraints. For every ${{\boldsymbol{q}}} \geq 0$, the problem is feasible because of the physical meaning of the problem (i.e., assigning flows to the network) as long as the system has enough capacity (i.e., dispatching enough vehicles). Hence, the lemma directly follows Theorem 1 in \citet{martin1975continuity}, which implements Berge's Maximum Theorem in parametric linear programming.

\section{Proof of Proposition \ref{eq_prop_bound_stt_diff}}\label{append_proof1}
According to Lemma \ref{lemma_berge}: 
\begin{align}
    & \left| \mathbb{E}_{\boldsymbol{Q}({\boldsymbol{x}^*})}[STT(\boldsymbol{Q}({\boldsymbol{x}^*}))]- SST(\boldsymbol{q}^*)\right| = \sum_{\hat{\boldsymbol{q}}\in\mathcal{Q}({\boldsymbol{x}^*})} |SST(\hat{\boldsymbol{q}}) - SST(\boldsymbol{q}^*)|\cdot \mathbb{P}_{\boldsymbol{Q}({\boldsymbol{x}^*})}(\hat{\boldsymbol{q}}) \nonumber \\
    & \leq \sum_{\hat{\boldsymbol{q}}\in\mathcal{Q}({\boldsymbol{x}^*})} L\cdot\norm{\hat{\boldsymbol{q}} - \boldsymbol{q}^*}_{1} \cdot \mathbb{P}_{\boldsymbol{Q}({\boldsymbol{x}^*})}(\hat{\boldsymbol{q}})  .
\end{align}
Let us divide the support of the random variable $\boldsymbol{Q}({\boldsymbol{x}^*})$ as three mutually exclusive subsets:
\begin{align}
 &   \mathcal{Q}({\boldsymbol{x}^*})^{\text{Leq}} = \mathcal{Q}({\boldsymbol{x}^*})\cap\{\hat{\boldsymbol{q}}:\; 0 \leq \hat{\boldsymbol{q}} < \mathbb{E}[\boldsymbol{Q}({\boldsymbol{x}^*})] - \boldsymbol{\epsilon}\},\\
&    \mathcal{Q}({\boldsymbol{x}^*})^{\text{Mid}} = \mathcal{Q}({\boldsymbol{x}^*})\cap\{\hat{\boldsymbol{q}}:\; \mathbb{E}[\boldsymbol{Q}({\boldsymbol{x}^*})] - \boldsymbol{\epsilon} \leq \hat{\boldsymbol{q}} \leq \mathbb{E}[\boldsymbol{Q}({\boldsymbol{x}^*})] + \boldsymbol{\epsilon}\},\\
&    \mathcal{Q}({\boldsymbol{x}^*})^{\text{Geq}} = \mathcal{Q}({\boldsymbol{x}^*})\cap\{\hat{\boldsymbol{q}}:\; \mathbb{E}[\boldsymbol{Q}({\boldsymbol{x}^*})] + \boldsymbol{\epsilon} < \hat{\boldsymbol{q}} \leq \boldsymbol{q}^{\text{Max}}\},
\end{align}
where $ \mathcal{Q}({\boldsymbol{x}^*}) = \mathcal{Q}({\boldsymbol{x}^*})^{\text{Leq}} \cup\mathcal{Q}({\boldsymbol{x}^*})^{\text{Mid}} \cup\mathcal{Q}({\boldsymbol{x}^*})^{\text{Geq}} $.

We can calculate the summation over these three subsets separately:

\noindent \textbf{(1) Bounds on the summation over} $\mathcal{Q}({\boldsymbol{x}^*})^{\text{Leq}}$:
\begin{align}
    \sum_{\hat{\boldsymbol{q}}\in\mathcal{Q}({\boldsymbol{x}^*})^{\text{Leq}}} L\cdot\norm{\hat{\boldsymbol{q}} - \boldsymbol{q}^*}_{1} \cdot \mathbb{P}_{\boldsymbol{Q}({\boldsymbol{x}^*})}(\hat{\boldsymbol{q}}) 
    & \leq \sum_{\hat{\boldsymbol{q}}\in\mathcal{Q}({\boldsymbol{x}^*})^{\text{Leq}}} L\cdot\left(\norm{\hat{\boldsymbol{q}} - \mathbb{E}[\boldsymbol{Q}({\boldsymbol{x}^*})]}_{1} + \norm{\mathbb{E}[\boldsymbol{Q}({\boldsymbol{x}^*})] - \boldsymbol{q}^*}_{1}\right) \cdot \mathbb{P}_{\boldsymbol{Q}({\boldsymbol{x}^*})}(\hat{\boldsymbol{q}}) ,
\end{align}
which is followed by the triangle inequality. Notice that
\begin{align}
& \sum_{\hat{\boldsymbol{q}}\in\mathcal{Q}({\boldsymbol{x}^*})^{\text{Leq}}} L\cdot \norm{\hat{\boldsymbol{q}} - \mathbb{E}[\boldsymbol{Q}({\boldsymbol{x}^*})] }_{1} \cdot \mathbb{P}_{\boldsymbol{Q}({\boldsymbol{x}^*})}(\hat{\boldsymbol{q}})  \leq L\cdot\sum_{\hat{\boldsymbol{q}}\in\mathcal{Q}({\boldsymbol{x}^*})^{\text{Leq}}} 
\norm{\mathbb{E}[\boldsymbol{Q}({\boldsymbol{x}^*})]}_1 \cdot \mathbb{P}_{\boldsymbol{Q}({\boldsymbol{x}^*})}(\hat{\boldsymbol{q}}) \nonumber \\
&=
L\cdot\norm{\mathbb{E}[\boldsymbol{Q}({\boldsymbol{x}^*})]}_1\cdot\mathbb{P}\big[ \boldsymbol{Q}({\boldsymbol{x}^*}) \leq \mathbb{E}[\boldsymbol{Q}({\boldsymbol{x}^*})]  - \boldsymbol{\epsilon}\big]\leq L\cdot\norm{\mathbb{E}[\boldsymbol{Q}({\boldsymbol{x}^*})]}_1\cdot \norm{\boldsymbol{\Gamma}}_2^2,
\label{eq_prop2_leq_0}
\end{align}
where the last inequality is the result of the following:
\begin{align}
   & \mathbb{P}\big[ \boldsymbol{Q}({\boldsymbol{x}^*}) \leq \mathbb{E}[\boldsymbol{Q}({\boldsymbol{x}^*})]  - \boldsymbol{\epsilon}\big]\leq \mathbb{P}\big[ |\boldsymbol{Q}({\boldsymbol{x}^*}) - \mathbb{E}[\boldsymbol{Q}({\boldsymbol{x}^*})]| \geq \boldsymbol{\epsilon}\big] \nonumber \\
    &\leq \sum_{i\in\mathcal{F}\times\mathcal{T}} \mathbb{P}[|{Q_i}|_{\boldsymbol{x}^*} - \mathbb{E}[{Q_i}|_{\boldsymbol{x}^*}]| \geq \epsilon_i] \leq \sum_{i\in\mathcal{F}\times\mathcal{T}} (\Gamma_i)^2 = \norm{\boldsymbol{\Gamma}}_2^2,
\end{align}
where the inequality is followed by the union bound and the $\Gamma$-concentration property. Similarly, we have
\begin{align}
& \sum_{\hat{\boldsymbol{q}}\in\mathcal{Q}({\boldsymbol{x}^*})^{\text{Leq}}} L\cdot \norm{\mathbb{E}[\boldsymbol{Q}({\boldsymbol{x}^*})] - {\boldsymbol{q}}^*}_{1} \cdot \mathbb{P}_{\boldsymbol{Q}({\boldsymbol{x}^*})}(\hat{\boldsymbol{q}})  \leq
\sum_{\hat{\boldsymbol{q}}\in\mathcal{Q}({\boldsymbol{x}^*})^{\text{Leq}}} L\cdot\norm{\boldsymbol{\epsilon}}_{1}\cdot \mathbb{P}_{\boldsymbol{Q}({\boldsymbol{x}^*})}(\hat{\boldsymbol{q}}) \nonumber \\
&=
L\cdot\norm{\boldsymbol{\epsilon}}_{1}\cdot\mathbb{P}\big[ \boldsymbol{Q}({\boldsymbol{x}^*}) \leq \mathbb{E}[\boldsymbol{Q}({\boldsymbol{x}^*})]  - \boldsymbol{\epsilon}\big]\leq L\cdot\norm{\boldsymbol{\epsilon}}_{1}\cdot \norm{\boldsymbol{\Gamma}}_2^2,
\label{eq_prop2_leq_1}
\end{align}
where the first inequality is due to the $\epsilon$-feasibility.

Therefore, combining (\ref{eq_prop2_leq_0}) and (\ref{eq_prop2_leq_1}) leads to 
\begin{align}
    \sum_{\hat{\boldsymbol{q}}\in\mathcal{Q}({\boldsymbol{x}^*})^{\text{Leq}}} L\cdot\norm{\hat{\boldsymbol{q}} - \boldsymbol{q}^*}_{1} \cdot \mathbb{P}_{\boldsymbol{Q}({\boldsymbol{x}^*})}(\hat{\boldsymbol{q}})  \leq L\cdot\big(\norm{\mathbb{E}[\boldsymbol{Q}({\boldsymbol{x}^*})]}_1 + \norm{\boldsymbol{\epsilon}}_{1}\big)\cdot \norm{\boldsymbol{\Gamma}}_2^2.
\end{align}

\noindent \textbf{(2) Bounds on the summation over} $\mathcal{Q}({\boldsymbol{x}^*})^{\text{Mid}}$:
\begin{align}
    \sum_{\hat{\boldsymbol{q}}\in\mathcal{Q}({\boldsymbol{x}^*})^{\text{Mid}}} L\cdot\norm{\hat{\boldsymbol{q}} - \boldsymbol{q}^*}_{1} \cdot \mathbb{P}_{\boldsymbol{Q}({\boldsymbol{x}^*})}(\hat{\boldsymbol{q}}) 
    & \leq \sum_{\hat{\boldsymbol{q}}\in\mathcal{Q}({\boldsymbol{x}^*})^{\text{Leq}}} L\cdot\left(\norm{\hat{\boldsymbol{q}} - \mathbb{E}[\boldsymbol{Q}({\boldsymbol{x}^*})]}_{1} + \norm{\mathbb{E}[\boldsymbol{Q}({\boldsymbol{x}^*})] - \boldsymbol{q}^*}_{1}\right) \cdot \mathbb{P}_{\boldsymbol{Q}({\boldsymbol{x}^*})}(\hat{\boldsymbol{q}}) \nonumber\\
    & \leq \sum_{\hat{\boldsymbol{q}}\in\mathcal{Q}({\boldsymbol{x}^*})^{\text{Mid}}} L\cdot(\norm{\boldsymbol{\epsilon}}_1+\norm{\boldsymbol{\epsilon}}_1 )\cdot \mathbb{P}_{\boldsymbol{Q}({\boldsymbol{x}^*})}(\hat{\boldsymbol{q}}) \leq 2L\cdot\norm{\boldsymbol{\epsilon}}_1 .
\end{align}

\noindent \textbf{(3) Bounds on the summation over} $\mathcal{Q}({\boldsymbol{x}^*})^{\text{Geq}}$:

Similar to the proof of $\mathcal{Q}({\boldsymbol{x}^*})^{\text{Leq}}$, notice that
\begin{align}
& \sum_{\hat{\boldsymbol{q}}\in\mathcal{Q}({\boldsymbol{x}^*})^{\text{Geq}}} L\cdot \norm{\hat{\boldsymbol{q}} - \mathbb{E}[\boldsymbol{Q}({\boldsymbol{x}^*})] }_{1} \cdot \mathbb{P}_{\boldsymbol{Q}({\boldsymbol{x}^*})}(\hat{\boldsymbol{q}})  \leq L\cdot\sum_{\hat{\boldsymbol{q}}\in\mathcal{Q}({\boldsymbol{x}^*})^{\text{Geq}}} 
\norm{\boldsymbol{q}^{\text{Max}}}_1 \cdot \mathbb{P}_{\boldsymbol{Q}({\boldsymbol{x}^*})}(\hat{\boldsymbol{q}}) \nonumber \\
&=
L\cdot\norm{\boldsymbol{q}^{\text{Max}}}_1\cdot\mathbb{P}\big[ \boldsymbol{Q}({\boldsymbol{x}^*}) \leq \mathbb{E}[\boldsymbol{Q}({\boldsymbol{x}^*})]  - \boldsymbol{\epsilon}\big]\leq L\cdot\norm{\boldsymbol{q}^{\text{Max}}}_1\cdot \norm{\boldsymbol{\Gamma}}_2^2
\label{eq_prop2_leq_3}.
\end{align}
Combining (\ref{eq_prop2_leq_3}) and (\ref{eq_prop2_leq_1}), we have
\begin{align}
    \sum_{\hat{\boldsymbol{q}}\in\mathcal{Q}({\boldsymbol{x}^*})^{\text{Geq}}} L\cdot\norm{\hat{\boldsymbol{q}} - \boldsymbol{q}^*}_{1} \cdot \mathbb{P}_{\boldsymbol{Q}({\boldsymbol{x}^*})}(\hat{\boldsymbol{q}}) 
    & \leq \sum_{\hat{\boldsymbol{q}}\in\mathcal{Q}({\boldsymbol{x}^*})^{\text{Geq}}} L\cdot\left(\norm{\hat{\boldsymbol{q}} - \mathbb{E}[\boldsymbol{Q}({\boldsymbol{x}^*})]}_{1} + \norm{\mathbb{E}[\boldsymbol{Q}({\boldsymbol{x}^*})] - \boldsymbol{q}^*}_{1}\right) \cdot \mathbb{P}_{\boldsymbol{Q}({\boldsymbol{x}^*})}(\hat{\boldsymbol{q}}) \nonumber\\
    &\leq L\cdot\big(\norm{\boldsymbol{q}^{\text{Max}}}_1 + \norm{\boldsymbol{\epsilon}}_{1}\big)\cdot \norm{\boldsymbol{\Gamma}}_2^2.
\end{align}

In summary, combining the summation over three mutually exclusive sets, we have:
\begin{align}
    & \left| \mathbb{E}_{\boldsymbol{Q}({\boldsymbol{x}^*})}[STT(\boldsymbol{Q}({\boldsymbol{x}^*}))]- SST(\boldsymbol{q}^*)\right| \leq 2L\cdot\norm{\boldsymbol{\epsilon}}_1  + L\cdot\big(\norm{\mathbb{E}[\boldsymbol{Q}({\boldsymbol{x}^*})]}_1  + \norm{\boldsymbol{q}^{\text{Max}}}_1 + 2\norm{\boldsymbol{\epsilon}}_{1}\big)\cdot \norm{\boldsymbol{\Gamma}}_2^2.
\end{align}

\section{Proof of Proposition \ref{prop_eplison_lowerbound}}\label{append_proof2}
When $\boldsymbol{\epsilon} = 0$, we have $\boldsymbol{q} = \mathbb{E}[\boldsymbol{Q}(\boldsymbol{x})]$. Then:
\begin{align}
    G_{\text{EP}}(\boldsymbol{\epsilon} = 0)  = \min_{\boldsymbol{z}}\{g(\mathbb{E}[\boldsymbol{Q}(\boldsymbol{x})], \boldsymbol{z}):\; h_j(\mathbb{E}[\boldsymbol{Q}(\boldsymbol{x})], \boldsymbol{z}) \leq 0 \} \label{eq_our_ep}.
\end{align}
According to Jensen's inequality, we have:
\begin{align}
g(\mathbb{E}[\boldsymbol{Q}(\boldsymbol{x})], \boldsymbol{z}) \leq \mathbb{E}[ g(\boldsymbol{Q}(\boldsymbol{x}), \boldsymbol{z})], \quad h_j(\mathbb{E}[\boldsymbol{Q}(\boldsymbol{x})], \boldsymbol{x}) \leq \mathbb{E}[ h_j(\boldsymbol{Q}(\boldsymbol{x}), \boldsymbol{z})],\;\forall j\in\mathcal{J}.
\end{align}
Therefore, the proposed approach has a smaller objective function and a larger feasible space compared to the stochastic optimization formulation (\ref{eq_stoch_ep}), which makes it a lower bound of (\ref{eq_stoch_ep}).  

\section{Proof of Proposition \ref{prop_gamma_app}}\label{append_proof3}

\textbf{Step 1:} We first show that if $\text{Var}[\boldsymbol{Q}(\boldsymbol{x})]$ is bounded, then $\text{Var}[h_j(\boldsymbol{Q}(\boldsymbol{x}), \boldsymbol{z})]$ is also bounded.

Notice that for any random variable ${X}$, we have $\text{Var}[{X}] = \mathbb{E}[{X}^2] - (\mathbb{E}[{X}])^2 \leq \mathbb{E}[{X}^2]$. Hence, if we take ${X} = h_j(\boldsymbol{Q}(\boldsymbol{x}), \boldsymbol{z}) - h_j(\mathbb{E}[\boldsymbol{Q}(\boldsymbol{x})], \boldsymbol{z})$, we get
\begin{align}
    \text{Var}[h_j(\boldsymbol{Q}(\boldsymbol{x}), \boldsymbol{z}) - h_j(\mathbb{E}[\boldsymbol{Q}(\boldsymbol{x})], \boldsymbol{z})] = \text{Var}[h_j(\boldsymbol{Q}(\boldsymbol{x}), \boldsymbol{z})] \leq \mathbb{E}[(h_j(\boldsymbol{Q}(\boldsymbol{x}), \boldsymbol{z}) - h_j(\mathbb{E}[\boldsymbol{Q}(\boldsymbol{x})], \boldsymbol{z}))^2] \label{eq_sto_gamma1}.
\end{align}
From the Lipschitz continuity of $h_j(\cdot)$, we have
\begin{align}
    |h_j(\boldsymbol{Q}(\boldsymbol{x}), \boldsymbol{z}) - h_j(\mathbb{E}[\boldsymbol{Q}(\boldsymbol{x})], \boldsymbol{z})| \leq C\norm{\boldsymbol{Q}(\boldsymbol{x}) - \mathbb{E}[\boldsymbol{Q}(\boldsymbol{x})]}_2,
\end{align}
which further yields:
\begin{align}
    \mathbb{E}[(h_j(\boldsymbol{Q}(\boldsymbol{x}), \boldsymbol{z}) - h_j(\mathbb{E}[\boldsymbol{Q}(\boldsymbol{x})], \boldsymbol{z}))^2] \leq C^2 \mathbb{E}[\norm{\boldsymbol{Q}(\boldsymbol{x}) - \mathbb{E}[\boldsymbol{Q}(\boldsymbol{x})]}^2_2] = C^2\cdot \sum_{i=1}^n\mathbb{E}\left[(Y_i - \mathbb{E}[Y_i])^2\right] = C^2\sum_{i=1}^n\text{Var}[Y_i].
\end{align}
Combining with (\ref{eq_sto_gamma1}), we have
\begin{align}
    \text{Var}[h_j(\boldsymbol{Q}(\boldsymbol{x}), \boldsymbol{z})] \leq C^2 \cdot \sum_{i=1}^n\text{Var}[Y_i] \leq C^2\cdot \norm{\boldsymbol{\Gamma}}_2^2.
\end{align}

\textbf{Step 2:} We then show that $\text{Var}[h_j(\boldsymbol{Q}(\boldsymbol{x}), \boldsymbol{z})] \leq  C^2\cdot \norm{\boldsymbol{\Gamma}}_2^2$ can deduce a weaker version of the chance constraints. 

Consider Chebyshev's inequality, for a given positive number $a > 0$:
\begin{align}
    \mathbb{P}[|h_j(\boldsymbol{Q}(\boldsymbol{x}), \boldsymbol{z}) - \mathbb{E}[h_j(\boldsymbol{Q}(\boldsymbol{x}), \boldsymbol{z})]| > a] \leq \frac{\text{Var}[h_j(\boldsymbol{Q}(\boldsymbol{x}), \boldsymbol{z})]}{a^2} \label{eq_cheb2}.
\end{align}
Equation (\ref{eq_cheb2}) implies
\begin{align}
    & \mathbb{P}[|h_j(\boldsymbol{Q}(\boldsymbol{x}), \boldsymbol{z}) - \mathbb{E}[h_j(\boldsymbol{Q}(\boldsymbol{x}), \boldsymbol{z})]| \leq a] \geq 1 - \frac{\text{Var}[h_j(\boldsymbol{Q}(\boldsymbol{x}), \boldsymbol{z})]}{a^2} \notag  \\ 
     \Rightarrow \;  &  \mathbb{P}[h_j(\boldsymbol{Q}(\boldsymbol{x}), \boldsymbol{z})  \leq a +\mathbb{E}[h_j(\boldsymbol{Q}(\boldsymbol{x}), \boldsymbol{z})]] \geq 1 - \frac{\text{Var}[h_j(\boldsymbol{Q}(\boldsymbol{x}), \boldsymbol{z})]}{a^2} \label{eq_cheb3}.
\end{align}
Since we know that $\mathbb{E}[h_j(\boldsymbol{Q}(\boldsymbol{x}), \boldsymbol{z})] \leq 0$, (\ref{eq_cheb3}) yields:
\begin{align}
 \mathbb{P}[h_j(\boldsymbol{Q}(\boldsymbol{x}), \boldsymbol{z})  \leq a ] \geq 1 - \frac{\text{Var}[h_j(\boldsymbol{Q}(\boldsymbol{x}), \boldsymbol{z})]}{a^2} \geq 1 - \left(\frac{C}{a}\right)^2\cdot\norm{\boldsymbol{\Gamma}}_2^2 \label{eq_cheb4}.
\end{align}
Let us pick $a = \frac{C\norm{\boldsymbol{\Gamma}}_2}{\sqrt{1-\eta}}$, we have 
\begin{align}
 \mathbb{P}\left[h_j(\boldsymbol{Q}(\boldsymbol{x}), \boldsymbol{z})  \leq \frac{C\norm{\boldsymbol{\Gamma}}_2}{\sqrt{1-\eta}} \right] \geq \eta .
 \label{eq_prop_gamma_approx}
\end{align}
Combining the two steps finishes the proof.

\section{Formulation of the Benders Decomposition}\label{append_benders}
\subsection{Subproblem}
The subproblem is derived by fixing the decision variables $\boldsymbol{x}$, and only considering the components including $\boldsymbol{q}$ and $\boldsymbol{z}$.
\begin{subequations}
\label{eq_sp}
\begin{align}
    [SP(\boldsymbol{x})] \quad \min_{\boldsymbol{q}, \boldsymbol{z}} \quad &  WT(\boldsymbol{q},\boldsymbol{z}) + IVT(\boldsymbol{z})   \\
    \text{s.t.} \quad
    & \text{Constraints } (\ref{const_OF1}) - (\ref{const_OF2}),\\
    & \text{Constraints } (\ref{eq_feasibility1}) - (\ref{eq_feasibility2}).
\end{align}
\end{subequations}
The objective of the dual problem of (\ref{eq_sp}) is
\begin{align}
    D(\boldsymbol{\alpha},\boldsymbol{\beta},\boldsymbol{\gamma},\boldsymbol{\iota},\boldsymbol{\kappa},\boldsymbol{\rho};\boldsymbol{x}) = &\sum_{l \in \mathcal{L}}\sum_{t \in \Tcal} \sum_{t' = t}^{T_{l,t}} K_{l,t}\alpha_{l,t,t'} +\sum_{(u,v,r) \in \mathcal{F}}\sum_{t \in \Tcal} \sum_{t' = t^{\text{min}}}^t f_{t'}^{u,v,r}\beta^{u,v,r}_t \notag\\
    & + \sum_{(u,v,r,i,t) \in \Omega_1}\hat{z}_{t}^{u,v,r,i}\gamma_{t}^{u,v,r,i} + \sum_{(u,v)\in\Wcal}\sum_{t\in\Tcal} d^{u,v}_t \iota^{u,v}_t \notag\\
    & +\sum_{(u,v,r) \in \mathcal{F}}\sum_{t \in \Tcal}\kappa_{t}^{u,v,r} \cdot (1-\epsilon)\sum_{p\in \mathcal{P}^{u,v}_{t}} \sum_{r' \in \mathcal{R}^{u,v}}x_{p,r'} \cdot \pi_{p,r'}^r \notag\\
    & +\sum_{(u,v,r) \in \mathcal{F}}\sum_{t \in \Tcal}\rho_{t}^{u,v,r} \cdot (1+\epsilon)\sum_{p\in \mathcal{P}^{u,v}_{t}} \sum_{r' \in \mathcal{R}^{u,v}}x_{p,r'} \cdot \pi_{p,r'}^r,
\end{align}
where $\boldsymbol{\alpha}, \boldsymbol{\beta}, \boldsymbol{\gamma}, \boldsymbol{\iota}$ are the dual variables associated with constraints (\ref{eq_cap_const}), (\ref{eq_flow_cons_1}), (\ref{eq_flow_cons_2}), (\ref{eq_demand}), respectively. $\boldsymbol{\kappa}$, $\boldsymbol{\rho}$ are the dual variables associated with constraints (\ref{eq_feasibility1}) and (\ref{eq_feasibility2}). Let $\boldsymbol{x} := (\boldsymbol{\alpha}, \boldsymbol{\beta}, \boldsymbol{\gamma}, \boldsymbol{\iota}, \boldsymbol{\kappa}$, $\boldsymbol{\rho})$.
If the dual problem of (\ref{eq_sp}) is feasible and bounded with a solution $\boldsymbol{x}^*$, the following optimality cut is added to the master problem:
\begin{align}
Z \geq D(\boldsymbol{x}^*;\boldsymbol{x})  ,
\end{align}
where $Z$ is a decision variable in the master problem. If the dual problem of (\ref{eq_sp}) is unbounded, and $\boldsymbol{x}^*$ is an optimal extreme ray of the dual, the following feasibility cut is added to the master problem:
\begin{align}
D(\boldsymbol{x}^*;\boldsymbol{x}) \leq 0.
\end{align}

\subsection{Master problem}
Let $\Acal^{\text{O}}$ be the set of solutions $\boldsymbol{x}^* $ of optimality cuts and $\Acal^{\text{F}}$ be the set of solutions $\boldsymbol{x}^*$ of feasibility cuts. At each iteration of the BD, a cut based on the solution of the subproblem is added to the respective set, and the corresponding master problem is defined as follows:
\begin{subequations}
\begin{align}
    [MP(\Acal^{\text{O}}, \Acal^{\text{F}})] \quad \min_{\boldsymbol{x}\in\mathcal{X}, Z} & \quad  \Psi \sum_{p\in\mathcal{P}} \sum_{r\in\mathcal{R}_p} -x_{p,r} \cdot V_{p,r}  + Z  \\
    \text{s.t.} \quad
    &  Z \geq D(\boldsymbol{x}^*;\boldsymbol{x})   \quad \forall \boldsymbol{x}^* \in \Acal^{\text{O}} ,\\
    & D(\boldsymbol{x}^*;\boldsymbol{x}) \leq 0 \quad \forall \boldsymbol{x}^* \in \Acal^{\text{F}},\\
    & \text{Constraints } (\ref{eq_conce}).
\end{align}
\end{subequations}
Note that the master problem has a smaller scale compared to the original problem (because there are no $\boldsymbol{z}$ and $\boldsymbol{q}$), which can be solved efficiently. 

\subsection{Convergence}
Let $(\boldsymbol{x}^{(k)}, Z^{(k)})$ and $(\boldsymbol{q}^{(k)}, \boldsymbol{z}^{(k)})$ be the solutions of the master problem and subproblem, respectively, in the $k$-th iteration. Then, the upper ($UB^{(k)}$) and lower ($LB^{(k)}$) bounds at the $k$-th iteration are given by:
\begin{align}
    UB^{(k)} &= \Psi \sum_{p\in\mathcal{P}} \sum_{r\in\mathcal{R}_p} -x_{p,r}^{(k)} \cdot V_{p,r}  + WT(\boldsymbol{q}^{(k)},\boldsymbol{z}^{(k)}) + IVT(\boldsymbol{z}^{(k)}),\\
    LB^{(k)} &= \Psi \sum_{p\in\mathcal{P}} \sum_{r\in\mathcal{R}_p} -x_{p,r}^{(k)} \cdot V_{p,r}  + Z^{(k)} .
\end{align}
$LB^{(k)}$ will keep increasing as $k$ increases because more cuts are added to the master problem. $UB^{(k)}$ does not necessarily decrease at every iteration. The convergence criterion is 
\begin{align}
    \text{Gap}^{(k)} = \frac{UB^{(k)} - LB^{(k)}}{LB^{(k)}} \leq \text{Predetermined threshold}.
\end{align}

\section{Inference of status quo choices}\label{app_infer_status_quo}
The status quo path choice inference method is based on our previous study \citep{mo2021inferring}, which is also similar to the trip-train method used for destination inference in open public transit systems (i.e., no tap-out). 

\textbf{[In the system when the incident happens]}: Consider a passenger $p \in \mathcal{P}$ with an incident line tap-in record before the end of the incident, meaning that he/she was in the transit system when the incident happens. We then track his/her next tap-in record. If his/her next tap-in is a transfer at a nearby bus or rail station, we can identify his/her chosen path based on the transfer station. We can also identify the waiting passenger if he/she continues to use the incident line to his/her intended destination inferred by his/her next tap-in records. 

\textbf{[Out of the system when the incident happens]}: For a passenger $p \in \mathcal{P}$ with only a tap-in record in nearby bus or rail stations. He/she may be affected by the incident to change the tap-in station, or just use the service as a normal commute. To identify whether he/she was affected, we extracted his/her travel histories on previous days without incidents to get the normal commute trajectories. If his/her tap-in time and location on the incident day have never appeared in the historical records before, we treat him/her as a passenger affected by the incident and identify his/her chosen path based on the tap-in station. 

For passengers in $\mathcal{P}$ without next tap-in records or travel histories, we randomly assign him/her a status quo path based on the proportion of inferred passengers. The final resulted path shares are shown in Figure \ref{fig_path_share_status_quo}. Around 49\% of the passengers chose to wait. 22\% of them chose the parallel bus lines. Others either took NS or WE bus routes and transferred to rails.  
\begin{figure}[htb]
\centering
\includegraphics[width = 0.5 \linewidth]{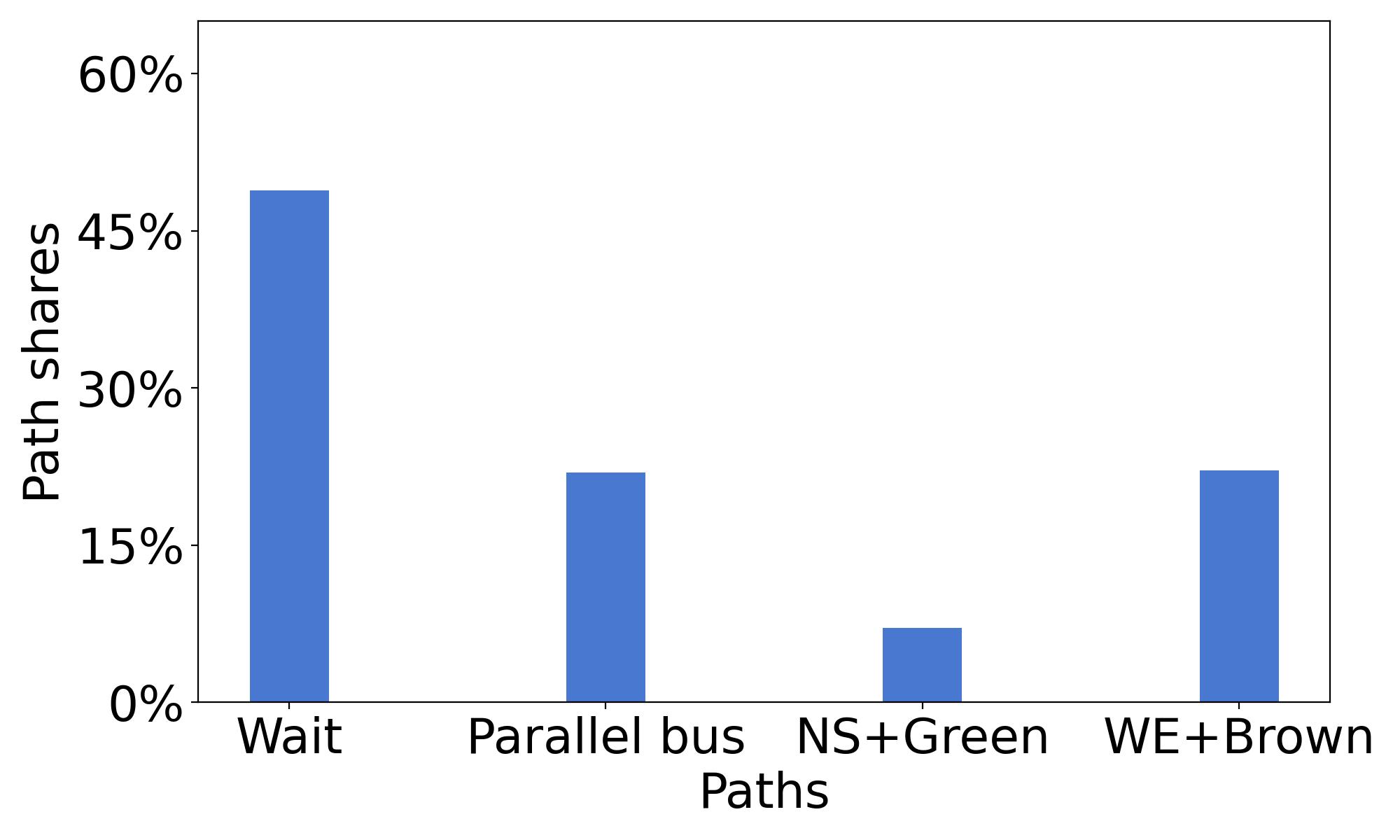}
\caption{Path shares of the inferred status quo choices}
\label{fig_path_share_status_quo}
\end{figure}

\newpage
\section{Individual travel time comparison}\label{append_ind_tt}
It is worth noting that (\ref{eq_of}) does not explicitly output the travel time of passengers using different paths. The travel time of passengers using path $(u,v,r)$ for trips departing at time $t$ (denoted as $TT^{u,v,r}_t$) has to be obtained from the network flow patterns \textbf{after} solving (\ref{eq_of}). Specifically, consider the group of passengers using path $(u,v,r)$ and departing at time $t$. Their arrival time at the destination given realized passenger flow $\hat{\boldsymbol{q}}$ (denoted as $AT^{u,v,r}_t(\hat{\boldsymbol{q}})$) can be calculated as 
\begin{align}
AT^{u,v,r}_t(\hat{\boldsymbol{q}})  = \min \Bigg\{ \tilde{t} \in \mathcal{T}^{u,v,r}_t  :  \sum_{t' = t^{\text{min}}}^{t} \left(f_{t'}^{u,v,r} + \hat{q}_{t'}^{u,v,r}\right) \leq    \sum_{t^{\text{min}} \leq t' +  \delta^{u,v,r,|\Ical^{u,v,r}|}_{t'} \leq {\tilde{t}} } z^{u,v,r,|\Ical^{u,v,r}|}_{t'}  
\Bigg\} \quad \notag
\\ \qquad \qquad \qquad \forall  t\in \mathcal{T}, (u,v,r) \in \mathcal{F}
\label{eq_tt_cal1},
\end{align}
where $\mathcal{T}^{u,v,r}_t$ is the set of possible arrival time indices, defined as $\mathcal{T}^{u,v,r}_t = \{t': t \leq t' \leq T\}$. Equation 
(\ref{eq_tt_cal1}) represents the travel time calculation with cumulative demand curves at origins and destinations. $ \sum_{t' = t^{\text{min}}}^{t} \left(f_{t'}^{u,v,r} + \hat{q}_{t'}^{u,v,r}\right)$ is the cumulative demand up to time $t$ at the origin. $\sum_{t^{\text{min}} \leq t' +  \delta^{u,v,r,|\Ical^{u,v,r}|}_{t'} \leq {\tilde{t}} } z^{u,v,r,|\Ical^{u,v,r}|}_{t'}$ is the cumulative passengers arriving at the destination up to time $t'$. When the cumulative arrivals at the destination are greater or equal to the cumulative demand at the origin (up to time $t$), all passengers finish the trip. So taking the minimum over $t'$ gives the arrival time for passengers departing at $t$. The travel time is then simply:
\begin{align}
    TT^{u,v,r}_t = AT^{u,v,r}_t - t \quad \forall t\in \mathcal{T}, (u,v,r) \in \mathcal{F}.
\label{eq_tt_cal}
\end{align}
Figure \ref{fig_ex_tt_cal} illustrates the travel time calculation.
\begin{figure}[H]
\centering
\includegraphics[width = 0.5 \linewidth]{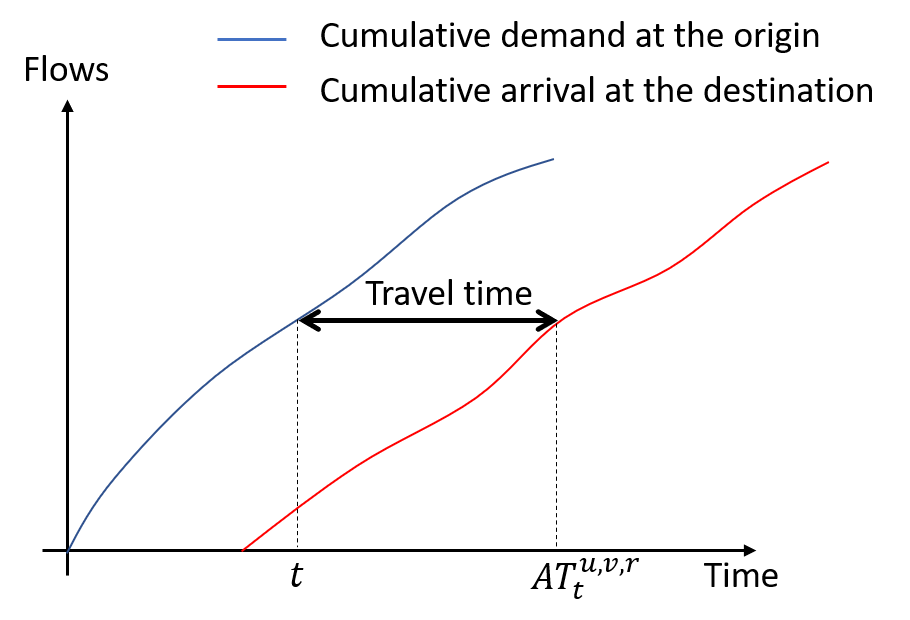}
\caption{Illustration of travel time calculation}
\label{fig_ex_tt_cal}
\end{figure}

Based on the above formulas, we calculate the travel time saving for each individual in the system. The distribution is shown in Figure \ref{fig_compare_tt}, where the negative values imply that the proposed model (IPR) has a lower travel time. For the comparison between IPR and status quo, we saw a sizable population has lower travel time under the recommendation, the maximum saving can be 50 minutes. For the comparison between the IPR and the capacity-based method, the most time saving is within 10 minutes. Note that there are some passengers experiencing higher travel time. This may be due to two reasons. First, there may be multiple ways to recommend passengers that achieve the same system travel time. As the model has no sense of its original travel time, some passengers may be worse off. Second, in order to reach the system optimal, some passengers may need to switch to a worse route in order to make the system better. In order to address this issue, future studies may impose some equity-related constraints to ensure no passenger is provided with a worse path. 
\begin{figure}[htb]
\centering
\subfloat[IPR v.s. Status quo]{\includegraphics[width=0.4\textwidth]{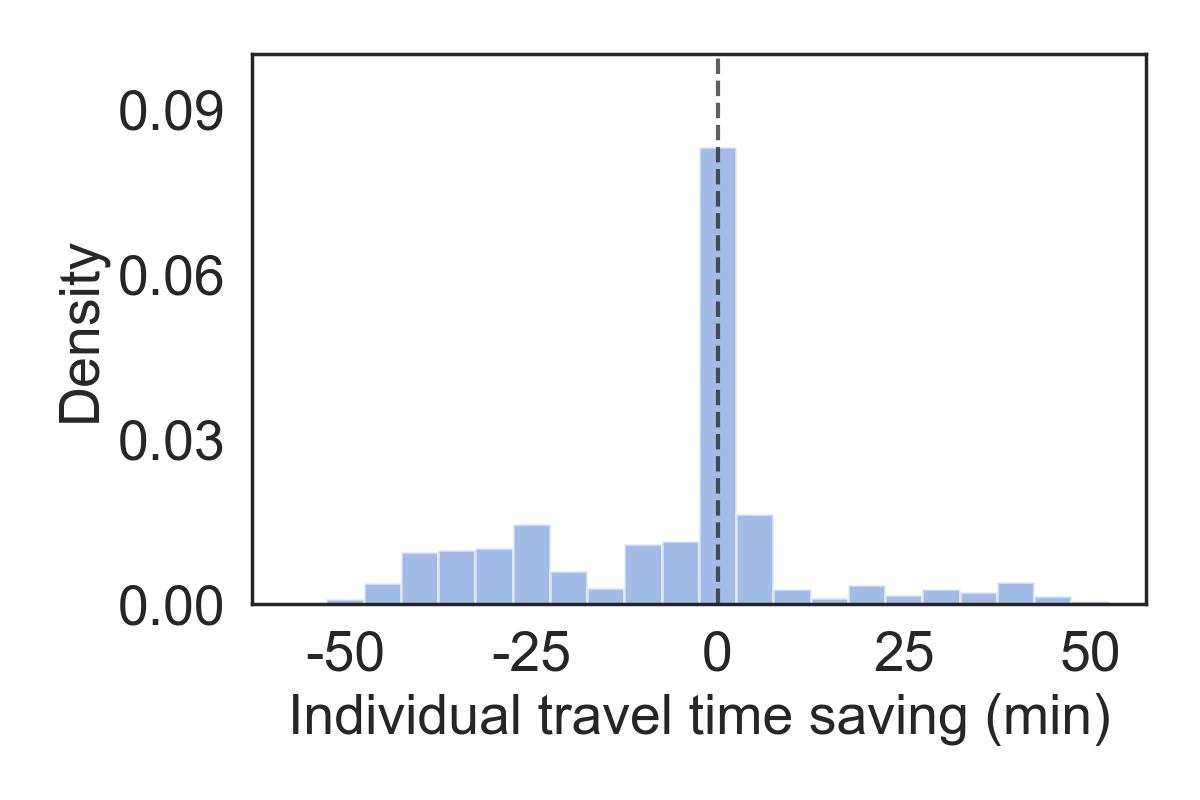}\label{fig_time_save1}}
\hfil
\subfloat[IPR v.s. Capacity-based]{\includegraphics[width=0.4\textwidth]{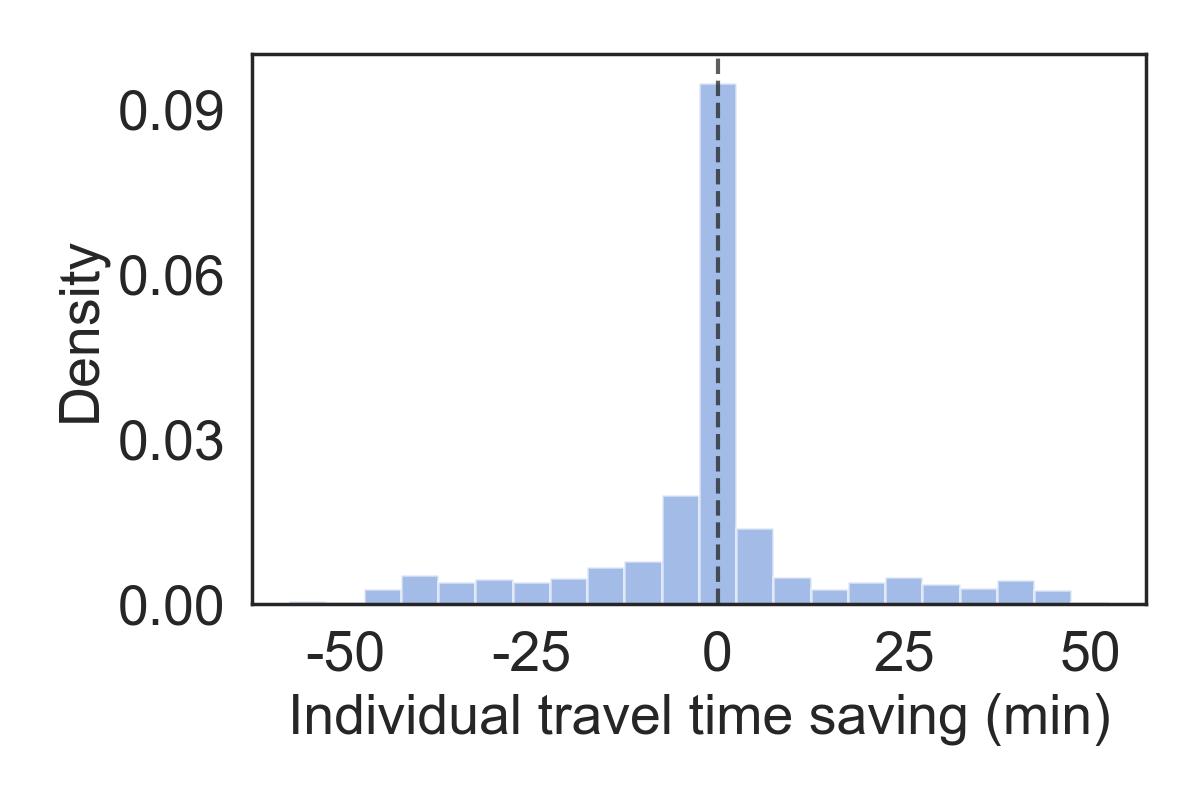}\label{fig_time_save2}}
\caption{Distribution of travel time saving (new travel time minus the old travel time, negative values imply lower travel time).}
\label{fig_compare_tt}
\end{figure}

\end{APPENDICES}

\bibliography{mybibfile}

\end{document}